\documentclass[12pt, reqno]{amsart}
\usepackage{amsmath, amsthm, amscd, amsfonts, amssymb, graphicx, color, mathtools}
\usepackage[bookmarksnumbered, colorlinks, plainpages]{hyperref}
\usepackage{enumerate}
\usepackage{setspace}
\usepackage{multicol}
\usepackage[margin=1in]{geometry}
\usepackage{comment}

\usepackage[normalem]{ulem} % DO NOT FORGET THE OPTION normalem
\newcommand\tsout{\bgroup\markoverwith{\textcolor{red}{\rule[0.5ex]{2pt}{1.4pt}}}\ULon}
\newcommand{\stkout}[1]{\ifmmode\text{\tsout{\ensuremath{#1}}}\else\tsout{#1}\fi}
\allowdisplaybreaks

\theoremstyle{definition}
\newtheorem{theorem}{Theorem}[section]
\newtheorem{lemma}[theorem]{Lemma}
\newtheorem{proposition}[theorem]{Proposition}
\newtheorem{corollary}[theorem]{Corollary}
\newtheorem{definition}[theorem]{Definition}

\theoremstyle{remark}
\newtheorem{remark}[theorem]{Remark}
\numberwithin{equation}{section}

\newcommand{\bff}{\boldsymbol}
\newcommand{\bb}{\mathbb}

\newcommand{\dt}{\mathrm{d}t}
\newcommand{\ddt}{\frac{\mathrm{d}}{\mathrm{d}t}}
\newcommand{\dx}{\mathrm{d}x}
\newcommand{\ds}{\mathrm{d}s}

\newcommand{\dtau}{\mathrm{d}\tau}

\newcommand{\norm}[2]{\left\|{#1}\right\|_{#2}}
\newcommand{\inpro}[2]{\left\langle#1,#2\right\rangle}
\newcommand{\binpro}[2]{\Big\langle#1,#2\Big\rangle}
\newcommand{\abs}[1]{\left|{#1}\right|}

\begin{document}
\setcounter{page}{1}

\title[Global solutions of LLBar equation]{{Global solutions
of the Landau--Lifshitz--Baryakhtar equation}}

\author[Agus L. Soenjaya]{Agus L. Soenjaya}
\address{School of Mathematics and Statistics, The University of New South Wales, Sydney 2052, Australia}
\email{\textcolor[rgb]{0.00,0.00,0.84}{a.soenjaya@unsw.edu.au}}

\author[Thanh Tran]{Thanh Tran}
\address{School of Mathematics and Statistics, The University of New South Wales, Sydney 2052, Australia}
\email{\textcolor[rgb]{0.00,0.00,0.84}{thanh.tran@unsw.edu.au}}

\date{\today}


\begin{abstract}
The Landau--Lifshitz--Baryakhtar (LLBar) equation is a generalisation of the
Landau--Lifshitz--Gilbert and the Landau--Lifshitz--Bloch equations
which takes into account contributions from nonlocal damping and is valid at
moderate temperature below the Curie temperature. Therefore,
it is used to explain some discrepancies between the experimental observations
and the known theories in various problems on magnonics and magnetic domain-wall
dynamics. In this paper, the existence and uniqueness of {global} weak,
{strong}, and regular solutions to LLBar equation are proven. H\"older
continuity of the solution is also discussed.
\end{abstract}
\maketitle

\section{Introduction}
{The theory of micromagnetism deals with the study of magnetic phenomena occurring
in ferromagnetic materials at sub-micrometre length scales.
A better understanding of the magnetisation dynamics at elevated temperature
would contribute to the development of ultrahigh-density storage technology
based on heat-assisted magnetic recording; see \cite{Col-etal23, MeoPanetal20,
Ran-etal22, StrRutEvaCha20} and the references therein.
}

A widely-studied equation which describes the
evolution of magnetic spin field in ferromagnetic material is the
Landau--Lifshitz--Gilbert (LLG) equation \cite{Gil55,LanLif35}. According to this
theory, the magnetisation of a magnetic body
$\Omega\subset\bb{R}^d$, $d\in\{1,2,3\}$, denoted by~$\bff{u}(t,\bff{x}) \in
\bb{R}^3$ for $t>0$ and $\bff{x}\in\Omega$, is described by
\begin{equation}\label{equ:LLE}
\frac{\partial \bff{u}}{\partial t}= -\gamma \bff{u} \times \bff{H}_{\text{eff}}
- \lambda \bff{u} \times (\bff{u}\times \bff{H}_{\text{eff}}),
\end{equation}
where $\gamma>0$ and $\lambda>0$ are the gyromagnetic ratio and a phenomenological damping parameter, respectively, and
$\bff{H}_{\text{eff}}$ is the effective field (consisting of the exchange field,
demagnetising field, external magnetic field and others). It is known that far below the Curie temperature, the
magnetisation of a ferromagnetic material preserves its magnitude. This
property is reflected in equation~\eqref{equ:LLE} (by taking the dot product of
both sides of the equation with $\bff{u}$). 

Mathematically, the LLG equation
has been extensively studied either on bounded or unbounded domains
where various existence, uniqueness and regularity properties were discussed.
A non-exhaustive list includes
\cite{AloSoy92,CarFab01,CimExist07,DiFInnPra20,FeiTra17b,GuoDin08,GuoHon93,GutdeL19,Vis85}.
Since then, various generalisations and improvements to the LLG equation have
been made in the physical and mathematical literatures. A widely-used physical
model for micromagnetism above the Curie temperature is
the Landau--Lifshitz--Bloch (LLB) equation \cite{Gar97}.
%which has also been studied mathematically \cite{Le16,LiGuoLiuLiu20}. 
This equation interpolates between the LLG
equation at low temperatures and the Ginzburg--Landau theory of phase
transitions, and is known not to preserve the magnitude of the magnetisation.
Mathematically, the existence and regularity properties for LLB
equation have been studied \cite{Le16,LiGuoLiuLiu21}.

The LLG and LLB equations, nevertheless, cannot account for some experimental data
and microscopic calculations. These include the nonlocal damping in magnetic
metals and crystals \cite{DvoVanVan13, Wei-etal14}, or the higher-than-expected spin wave
decrement for short-wave magnons~\cite{BarIvaSukMelYu97}. The Landau--Lifshitz--Baryakhtar
(LLBar) equation proposed by Baryakhtar \cite{Bar84,BarDan13,BarIvaSukMelYu97}
is based on Onsager's relations and generalises the LLG and LLB equations
\cite{DvoVanVan13,DvoVanVan14,Wan-etal15}. This equation has also been
implemented on several commonly-used micromagnetic simulation software, such as
\textsc{MuMax} \cite{Au-etal13,Lel-etal18} and \textsc{Fidimag} \cite{Wan15,Wan-etal15}.  Moreover,
various micromagnetic simulations provide evidence that the LLBar equation
agrees with some of the observed experimental findings in micromagnetics,
especially those related to ultrafast magnetisation at an elevated temperature;
see \cite{Au-etal13, DvoVanVan13, LiCheBerJia17, Wan15, Wan-etal15, Wei-etal14}
and the references therein.

The LLBar equation in its most general form \cite{BarIvaSukMelYu97, Wan-etal15} reads
\begin{equation*}
\frac{\partial\bff{u}}{\partial t}= -\gamma \bff{u}\times \bff{H}_{\text{eff}} + \bff{\Lambda}_r \cdot \bff{H}_{\text{eff}} - \bff{\Lambda}_{e, ij} \frac{\partial^2 \bff{H}_{\text{eff}}}{\partial \bff{x}_i \partial \bff{x}_j}, 
\end{equation*}
where $\bff{u}$ represents the magnetisation vector, $\bff{\Lambda}_r$ and
$\bff{\Lambda}_e$ denote the relaxation tensor and the exchange tensor,
respectively. Here, Einstein's summation notation is used. For a
polycrystalline, amorphous soft magnetic materials and magnetic metals at
moderate temperature (where nonlocal damping and longitudinal relaxation are
significant), this equation simplifies \cite{DvoVanVan13, Wan-etal15} to
\begin{equation*}
\frac{\partial\bff{u}}{\partial t} = -\gamma \bff{u} \times \bff{H}_{\text{eff}} + \lambda_r \bff{H}_{\text{eff}}- \lambda_e \Delta \bff{H}_{\text{eff}}.
\end{equation*}
where the positive scalars $\gamma$, $\lambda_r$, and $\lambda_e$ are the
electron gyromagnetic ratio, relativistic damping constant, and exchange damping
constant, respectively. The effective field $\bff{H}_{\text{eff}}$ is given by
\begin{equation*}
\bff{H}_{\text{eff}}= \Delta \bff{u}+ \frac{1}{2\chi}(1-|\bff{u}|^2)\bff{u} + \text{lower order terms},
\end{equation*}
with $\chi>0$ being the magnetic susceptibility of the material.

If the exchange interaction is dominant (as is the case for ordinary
ferromagnetic material), then $\bff{u}:[0,T]\times\Omega\to\bb{R}^3$ solves the following
problem:
\begin{subequations}\label{llbar}
%\left\{
\begin{alignat}{2}
&\frac{\partial\bff{u}}{\partial t} - \beta_1 \Delta \bff{u} + \beta_2 \Delta^2 \bff{u} 
   = \beta_3 (1-|\bff{u}|^2)\bff{u}-\beta_4 \bff{u}\times \Delta \bff{u} 
   + \beta_5 \Delta(|\bff{u}|^2 \bff{u})
	& \ &\text{in $(0,T)\times\Omega$},
	\label{llbar a}
	\\[1ex]
	& \bff{u}(0,\cdot) = \bff{u_0} &\ &\text{in $\Omega$},
	\label{llbar b}
	\\[1ex]
& \frac{\partial \bff{u}}{\partial \bff{n}} = 
  \frac{\partial(\Delta\bff{u})}{\partial \bff{n}}= \bff{0} 
	& \ &\text{on $(0,T)\times\partial\Omega$},
	\label{llbar c}
\end{alignat}
%\right.
\end{subequations}
where $\beta_1 = \lambda_r - \lambda_e/(2\chi)$ is a real constant (which may be
positive or negative), while $\beta_2, \ldots, \beta_5$ are
positive constants. Here, $\partial\Omega$ is the boundary of $\Omega$ with
exterior unit normal vector denoted by $\bff{n}$.

Typically, the constant $\beta_1$ will be positive since $\lambda_e/(2\chi)$ is
much smaller than $\lambda_r$. However, in certain situations occurring in
spintronics or magnonics where the wavelength of the magnons is approaching the
exchange length of the ferromagnetic material, $\lambda_e$ can be significant
\cite{DvoVanVan13}. Therefore, we allow $\beta_1$ to take positive or negative values in
\eqref{llbar}.

To the best of our knowledge, mathematical analysis of the LLBar equation
does not exist in the literature. In this paper, we prove the existence,
uniqueness, and regularity of a weak {and strong} solution to problem~\eqref{llbar} in one,
two and three spatial dimensions (see Theorem~\ref{the:weakexist}), by using
the Faedo--Galerkin approximation and compactness method. We also prove H\"older
continuity properties of the solution (Theorem~\ref{the:Hol}).
This gives a mathematical foundation for the rigorous theory of LLBar equation
which is not currently available in the literature.

Another advantage of studying the LLBar equation is 
for a given initial data $\bff{u}_0$, the weak solution to
the LLBar equation generally has better regularity compared to that of the LLG
or the LLB equation.  Moreover, it is known that the existence
of global solutions to the LLG equation in 2-D is only guaranteed for
sufficiently small initial data \cite{CarFab01, FeiTra17b}, whereas for general initial data,
solutions in 2-D could blow-up in finite time \cite{Har04}. As we show in this
paper, the solution to the LLBar equation exists globally.

{We note that a related model of magnetisation dynamics in the framework of
frustrated magnets (which takes into account local and nonlocal interactions)
has recently been explored in \cite{DorMel22}. The model is based on the LLG
equation involving the bilaplacian operator, where the magnitude of the
magnetisation is conserved in that case.}

The paper is organised as follows. In Section \ref{formula}, we introduce some
notations and formulate the main results. In Section \ref{faedo}, we establish some a priori
estimates that are needed for the proof of the main theorems.
Section~\ref{weaksol} is devoted to the proof of the main results.
Finally, we collect in the appendix some essential mathematical facts that are used throughout
the paper.

\section{Formulation of the main results}\label{formula}
\subsection{Notation}
We begin by defining some notations used in this paper. The function space $\bb{L}^p := \bb{L}^p(\Omega; \bb{R}^3)$ denotes the space of $p$-th integrable functions taking values in $\bb{R}^3$ and $\bb{W}^{k,p} := \bb{W}^{k,p}(\Omega; \bb{R}^3)$ denotes the Sobolev space of 
functions on $\Omega \subset \bb{R}^d$ taking values in $\bb{R}^3$. Also, we
write $\bb{H}^k := \bb{W}^{k,2}$. Here, $\Omega\subset \bb{R}^d$ for $d=1,2,3$
is an open domain with smooth boundary. The operator $\Delta$ denotes the Neumann Laplacian. The partial derivative
$\partial/\partial x_i$ will be written by $\partial_i$ for short.

If $X$ is a normed vector space, the spaces $L^p(0,T, X)$ and $W^{k,p}(0,T,X)$ denote respectively the usual Lebesgue and Sobolev spaces of functions on $(0,T)$ taking values in $X$. The space $C([0,T],X)$ denotes the space of continuous function on $[0,T]$ taking values in $X$. Throughout this paper, we denote the scalar product in a Hilbert space $H$ by $\langle \cdot, \cdot\rangle_H$ and its corresponding norm by $\|\cdot\|_H$. We will not distinguish between the scalar product of $\bb{L}^2$ vector-valued functions taking values in $\bb{R}^3$ and the scalar product of $\bb{L}^2$ matrix-valued functions taking values in $\bb{R}^{3\times 3}$, and still denote them by $\langle\cdot,\cdot\rangle_{\bb{L}^2}$.

%{
%(To me a vector in $\bb{R}^3$ is a column matrix and is written as
%\[
%	\bff{u} =
%	\begin{bmatrix}
%	u_1 \\ u_2 \\ u_3
%	\end{bmatrix}
%	= (u_1,u_2,u_3)
%	=
%	\begin{bmatrix}
%		u_1 & u_2 & u_3
%	\end{bmatrix}^\top.
%\]
%Note the round brackets and commas. Sometimes, I also write matrices with round
%brackets. However, I disagree with this notation, which some people use,
%\[
%	\bff{u} =
%	(u_1, u_2, u_3)^\top,
%\]
%because the $\top$ symbol indicates that
%$(u_1,u_2,u_3)$ is a matrix, but we don't use commas in matrices. I want to make
%clear of these symbols here.)
%}

The following frequently-used notations are collected here for the reader's
convenience. Firstly, for any vector $\bff{z}\in\bb{R}^3$ and
matrices~$\bff{A}$, $\bff{B}\in\bb{R}^{3\times d}$, we define
\begin{equation}\label{equ:dot cro 2}
	\begin{alignedat}{2}
	\bff{z}\cdot\bff{A} &:=
	\begin{bmatrix}
		\bff{z}\cdot\bff{A}^{(1)} & \cdots & \bff{z}\cdot\bff{A}^{(d)}
	\end{bmatrix}
	\in\bb{R}^{1\times d},
	\quad&
	\bff{A} \cdot \bff{B}
	&:=
	\sum_{j=1}^d \bff{A}^{(j)} \cdot \bff{B}^{(j)}
	\in\bb{R},
	\\
	\bff{z}\times\bff{A}
	&:=
	\begin{bmatrix}
		\bff{z}\times\bff{A}^{(1)} & \cdots & \bff{z}\times\bff{A}^{(d)}
	\end{bmatrix}
	\in\bb{R}^{3\times d},
	\qquad&
	\bff{A} \times \bff{B}
	&:=
	\sum_{j=1}^d \bff{A}^{(j)} \times \bff{B}^{(j)}
	\in\bb{R}^3,
	\end{alignedat}
\end{equation}
where~$\bff{A}^{(j)}$ and~$\bff{B}^{(j)}$ denote the $j^{\rm th}$ column of~$\bff{A}$
and~$\bff{B}$, respectively. 

Next, for any vector-valued
function~$\bff{v}=(v_1,v_2,v_3):\Omega\subset\bb{R}^d\to\bb{R}^3$, we define
\begin{equation}\label{equ:nab del}
	\left\{
\begin{alignedat}{2}
	&\nabla\bff{v} : \Omega\to\bb{R}^{3\times d}
	&\quad\text{by}\quad
	&\nabla\bff{v} :=
\begin{bmatrix}
	\partial_1 \bff{v} \cdots \partial_d \bff{v}
\end{bmatrix}
=
\begin{bmatrix}
	\partial_1 v_1 & \cdots & \partial_d v_1
	\\
	\partial_1 v_2 & \cdots & \partial_d v_2
	\\
	\partial_1 v_3 & \cdots & \partial_d v_3
\end{bmatrix},
\\
	& \frac{\partial\bff{v}}{\partial\bff{n}} :
	\partial\Omega\to\bb{R}^{3\times1}
	&\quad\text{by}\quad
	& \frac{\partial\bff{v}}{\partial\bff{n}} 
	:=
	(\nabla\bff{v})\bff{n}
	=
	\begin{bmatrix}
		\displaystyle 	\frac{\partial v_1}{\partial\bff{n}} &
		\displaystyle 	\frac{\partial v_2}{\partial\bff{n}} &
		\displaystyle 	\frac{\partial v_3}{\partial\bff{n}}
	\end{bmatrix}^\top,
\\
	&\Delta\bff{v} : \Omega\to\bb{R}^{3\times 1}
	&\quad\text{by}\quad
	&\Delta\bff{v} :=
\begin{bmatrix}
	\Delta v_1 & \Delta v_2 & \Delta v_3
\end{bmatrix}^\top,
\\
	&\Delta\nabla\bff{v} : \Omega\to\bb{R}^{3\times d}
	&\quad\text{by}\quad
	&\Delta\nabla\bff{v} :=
\begin{bmatrix}
	\Delta \partial_1 v_1 & \cdots & \Delta \partial_d v_1
\\ 
	\Delta\partial_1  v_2 & \cdots & \Delta \partial_d v_2
\\ 
	\Delta\partial_1  v_3 & \cdots & \Delta \partial_d v_3
\end{bmatrix}
=
\nabla \Delta \bff{v}
.
\end{alignedat}
\right.
\end{equation}
As a consequence, if~$\bff{u}$ and~$\bff{v}$ satisfy suitable assumptions 
and~$\partial \bff{u}/\partial\bff{n}=0$ (where~$\bff{n}$ is the outward normal vector
to~$\partial D$), then
\begin{equation}\label{equ:vec Gre}
\begin{aligned}
-\inpro{\Delta\bff{u}}{\bff{v}}_{\bb{L}^2}
&=
-\sum_{i=1}^3
\inpro{\Delta u_i}{v_i}_{L^2}
=
\sum_{i=1}^3
\inpro{\nabla u_i}{\nabla v_i}_{\bb{L}^2}
=
\inpro{\nabla\bff{u}}{\nabla\bff{v}}_{\bb{L}^2},
\\
\inpro{\bff{u}\times\Delta\bff{u}}{\bff{v}}_{\bb{L}^2}
&=
\inpro{\Delta\bff{u}}{\bff{v}\times\bff{u}}_{\bb{L}^2}
=
%-
%\inpro{\nabla\bff{u}}{\nabla(\bff{v}\times\bff{u}}_{\bb{L}^2}
%=
-
\inpro{\bff{u}\times\nabla\bff{u}}{\nabla\bff{v}}_{\bb{L}^2}.
\end{aligned}
\end{equation}

Finally, throughout this paper, the constant $C$ in the estimate denotes a
generic constant which takes different values at different occurrences. If
the dependence of $C$ on some variable, e.g.~$T$, is highlighted, we often write
$C(T)$. The notation $A \lesssim B$ means $A \le C B$ where the specific form of
the constant $C$ is not important to clarify.
%The positive constant $\varepsilon$ denotes a constant which is made
%sufficiently small in the estimates, but whose exact value is not important.

\subsection{Main results}
In the following, we define the notion of weak solutions to \eqref{llbar}. We
first multiply~\eqref{llbar a} (dot product) with a test function~$\bff{\phi}$,
integrate over~$\Omega$, and (formally) use integration by parts,
noting~\eqref{llbar c}, to obtain
\begin{align}\label{equ:formal} 
	& \binpro{\frac{\partial\bff{u}(t)}{\partial t}}{\bff{\phi}}_{\bb{L}^2}
	+
	\beta_1 \inpro{\nabla\bff{u}(t)}{\nabla\bff{\phi}}_{\bb{L}^2}
	+
	\beta_2 \inpro{\Delta\bff{u}(t)}{\Delta\bff{\phi}}_{\bb{L}^2}
	\nonumber\\
	&=
	\beta_3 \inpro{(1-|\bff{u}(t)|^2)\bff{u}(t)}{\bff{\phi}}_{\bb{L}^2}
	+
	\beta_4 \inpro{\bff{u}(t) \times
	\nabla\bff{u}(t)}{\nabla\bff{\phi}}_{\bb{L}^2}
	-
	\beta_5 \inpro{ \nabla\big(|\bff{u}(t)|^2
	\bff{u}(t)\big)}{\nabla\bff{\phi}}_{\bb{L}^2}.
\end{align}
We next find sufficient conditions for the terms on the right-hand
side to be well-defined {for $d=1,2,3$}. If $\bff{u}\in\bb{H}^1$, then
$\bff{u}\in\bb{L}^4$ so that
$|\bff{u}|^2 \bff{u} \in\bb{L}^{4/3}$. Therefore, the
term $ \inpro{|\bff{u}|^2\bff{u}}{\bff{\phi}}_{\bb{L}^2}$ is well defined
if~$\bff{u}\in\bb{H}^1$ and
$\bff{\phi}\in\bb{H}^1$. Moreover, if $\bff{u}\in\bb{H}^2$, then
$\bff{u}\in\bb{L}^{\infty}$. Thus, the second term~$\inpro{\bff{u}(t) \times
\nabla\bff{u}(t)}{\nabla\bff{\phi}}_{\bb{L}^2}$ and the third term
\[
\inpro{ \nabla\big(|\bff{u}(t)|^2 \bff{u}(t)\big)}{\nabla\bff{\phi}}_{\bb{L}^2}
=
\sum_{i=1}^{d} \sum_{j=1}^{3} 
\inpro{\partial_i\big(|\bff{v}|^2 v_j\big)}{\partial_i\phi_j}_{L^2}
\]
on the right-hand side of the above equation are also well defined if
$\bff{u}\in\bb{H}^2$ and $\bff{\phi}\in\bb{H}^1$. This motivates
the following definition of solutions to problem~\eqref{llbar}.

{
\begin{definition}\label{def:wea str sol}
	Given $T>0$ and $\bff{u}_0\in\bb{L}^2(\Omega)$, a function
	$\bff{u}:[0,T]\to \bb{H}^2$ is a \emph{weak solution} to the problem
	\eqref{llbar} if $\bff{u}$ belongs to $C([0,T]; \bb{L}^2) \cap L^2(0,T; \bb{H}^2)$ and satisfies
	\begin{align}\label{weakform}
		\inpro{ \bff{u}(t)}{\bff{\phi}}_{\bb{L}^2} 
		&+ 
		\beta_1 \int_0^t \inpro{\nabla\bff{u}(s)}{\nabla\bff{\phi}}_{\bb{L}^2} \ds 
		+ \beta_2 \int_0^t \inpro{\Delta \bff{u}(s)}{\Delta \bff{\phi}}_{\bb{L}^2} \ds  
		\nonumber\\
		&= \inpro{\bff{u_0}}{ \bff{\phi}}_{\bb{L}^2} 
		+ \beta_3 \int_0^t \inpro{ (1-|\bff{u}(s)|^2) \bff{u}(s)}{\bff{\phi} }_{\bb{L}^2} \ds 
		\nonumber\\
		&\quad + \beta_4 \int_0^t \inpro{ \bff{u}(s)\times \nabla \bff{u}(s)}{ \nabla\bff{\phi} }_{\bb{L}^2} \ds 
		- \beta_5 \int_0^t \inpro{ \nabla\left(|\bff{u}(s)|^2 \bff{u}(s)\right)}{
			\nabla\bff{\phi} }_{\bb{L}^2} \ds,
	\end{align}
	for all $\bff{\phi}\in \bb{H}^2$ and $t\in [0,T]$.
	
	A weak solution $\bff{u}: [0,T] \to \bb{H}^2$ is called a \emph{strong
	solution} if it belongs to $C([0,T];\bb{H}^2) \cap L^2(0,T; \bb{H}^4)$.
	In this case, it satisfies~\eqref{llbar} almost everywhere
	in~$[0,T]\times\Omega$.
\end{definition}
}

We now state the main theorems of the paper, the proofs of which will be given in
Section~\ref{weaksol} and Section~\ref{sec:Hol}. The
first theorem gives the existence, uniqueness, and regularity of the solution.

{
\begin{theorem}\label{the:weakexist}
	Let $\Omega\subset \bb{R}^d$, $d=1,2,3$, be a bounded domain with
	$C^{r+1,1}$-boundary and let $\bff{u}_0\in \bb{H}^r$, $r\in \{0,1,2,3\}$, be a
	given initial data. For any $T>0$, there exists a global
	weak solution to \eqref{llbar} such that
	\begin{equation}\label{equ:wea sol smo}
	\bff{u} \in C([0,T]; \bb{H}^r) \cap L^2(0,T; \bb{H}^{r+2}).
	\end{equation}
	Furthermore, this solution depends continuously on the
	$\bb{H}^r$-norm of the initial data, which implies uniqueness. More
	precisely, if $\bff{u}$ and $\bff{v}$ are solutions corresponding to the
	initial data~$\bff{u}_0$ and~$\bff{v}_0$, respectively, then the
	following estimates hold
	\begin{align}\label{equ:cont dep u-v L2}
		\norm{\bff{u}(t) - \bff{v}(t)}{\bb{L}^2}^2
		\lesssim
		\norm{\bff{u}_0- \bff{v}_0}{\bb{L}^2}^2 \,
		\exp \left(\int_0^t \Big(1+ \norm{\bff{u}(s)}{\bb{L}^\infty}^4
		+ \norm{\bff{v}(s)}{\bb{L}^\infty}^4 \Big) \ds \right)
	\end{align}
	and
	\begin{equation}\label{equ:uv uv0 Hr}
		\norm{\bff{u}(t)-\bff{v}(t)}{\bb{H}^r}
		\lesssim
		\norm{\bff{u}_0 - \bff{v}_0}{\bb{H}^r},
	\end{equation}
	where the constant depends on $T$.
	In particular, if the solution $\bff{u}$ to problem~\eqref{llbar} satisfies
	\begin{align}\label{equ:add reg unique}
		\int_0^T \norm{\bff{u}(s)}{\bb{L}^\infty}^4 \ds < \infty,
	\end{align}
	then it is unique.

	Moreover, if $r=2,3$ then the solution is a strong solution in the
	sense of Definition~\ref{def:wea str sol}.
\end{theorem}
}

The next theorem shows that the strong solution is H\"older continuous in time.

\begin{theorem}\label{the:Hol}
Let $T>0$ and $\bff{u}$ be the unique strong solution of \eqref{llbar} with initial data $\bff{u}_0\in \bb{H}^2$. Then
\[
	\bff{u} \in C^{0,\alpha}(0,T; \bb{L}^2) \cap C^{0,\beta}(0,T; \bb{L}^\infty),
\]
where $\alpha \in \left(0, \frac{1}{2}\right]$ and $\beta \in \left(0, \frac{1}{2} - \frac{d}{8} \right]$.
\end{theorem}
{
\begin{remark}\label{rem:rem assump}
It can be seen that assumption~\eqref{equ:add reg unique} holds in many
different cases.
\begin{enumerate}
\renewcommand{\labelenumi}{\theenumi}
\renewcommand{\theenumi}{{\rm (\roman{enumi})}}
\item If $r=0$, assumption~\eqref{equ:add reg unique} holds for $d=1$ or $d=2$.
	Indeed, in case $d=1$, we have by using the Gagliardo--Nirenberg
	inequality (Theorem~\ref{the:Gal Nir}) and the fact that $\bff{u}\in
	L^2(0,T;\bb{H}^2)$,
\begin{align*}
	\int_0^t \norm{\bff{u}(s)}{\bb{L}^\infty}^4 \ds 
	\lesssim
	\int_0^t \norm{\bff{u}(s)}{\bb{L}^2}^2 
	\norm{\bff{u}(s)}{\bb{H}^1}^2 \ds
	\leq
	\int_0^t \norm{\bff{u}(s)}{\bb{H}^1}^2 \ds \lesssim 1,
\end{align*}
and in case $d=2$,
\begin{align*}
	\int_0^t \norm{\bff{u}(s)}{\bb{L}^\infty}^4 \ds 
	\lesssim
	\int_0^t \norm{\bff{u}(s)}{\bb{L}^2}^2 
	\norm{\bff{u}(s)}{\bb{H}^2}^2 \ds
	\leq
	\int_0^t \norm{\bff{u}(s)}{\bb{H}^2}^2 \ds \lesssim 1.
\end{align*}
\item If $r\in \{1,2,3\}$, assumption~\eqref{equ:add reg unique} holds for
	$d=1,2,3$. In this case, with $d=3$, we have
\[
	\int_0^t \norm{\bff{u}(s)}{\bb{L}^\infty}^4 \ds
	\lesssim
	\int_0^t \norm{\bff{u}(s)}{\bb{H}^1}^2 \norm{\bff{u}(s)}{\bb{H}^2}^2 \ds
	\leq 
	\int_0^t \norm{\bff{u}(s)}{\bb{H}^2}^2 \ds
	\lesssim 1.
\]
\end{enumerate}
\end{remark}
}

\section{Faedo--Galerkin Approximation}\label{faedo}

Let $\{\bff{e_i}\}_{i=1}^\infty$ denote an orthonormal basis of $\bb{L}^2$
consisting of eigenvectors for $-\Delta$ such that
\begin{align*}
-\Delta \bff{e_i} =\lambda_i \bff{e_i} \ \text{ in $\Omega$}
\quad\text{and}\quad
 \frac{\partial \bff{e_i}}{\partial \bff{n}}= \bff{0} \ \text{ on } \partial \Omega,
\end{align*}
where $\lambda_i>0$ are the eigenvalues of $-\Delta$, associated with
$\bff{e}_i$. By elliptic regularity results, $\bff{e_i}$ is smooth up to the
boundary, and we also have
\begin{align*}
\Delta^2 \bff{e_i} = \lambda_i^2 \bff{e_i} \ \text{ in $\Omega$}
\quad\text{and}\quad
 \frac{\partial \bff{e_i}}{\partial \bff{n}}=\frac{\partial \Delta
 \bff{e_i}}{\partial \bff{n}} =\bff{0} \ \text{ on } \partial \Omega.
\end{align*}
Let $\bb{V}_n:= \text{span}\{\bff{e_1},\ldots,\bff{e_n}\}$ and $\Pi_n: \bb{L}^2\to
\bb{V}_n$ be the orthogonal projection defined by
\begin{align*}
\inpro{ \Pi_n \bff{v}}{\bff{\phi}}_{\bb{L}^2}
= \inpro{ \bff{v}}{\bff{\phi}
}_{\bb{L}^2} \quad \text{for all } 
\bff{\phi}\in \bb{V}_n \text{ and all }
\bff{v}\in\bb{L}^2.
\end{align*}
Note that $\Pi_n$ is self-adjoint and satisfies
\begin{equation}\label{equ:Pin}
\norm{\Pi_n \bff{v}}{\bb{L}^2} \leq \norm{\bff{v}}{\bb{L}^2} \quad \text{for all }
\bff{v}\in \bb{L}^2.
\end{equation}
%\begin{align*}
%\langle \bff{v}, \Pi_n \bff{w}\rangle_{\bb{L}^2} = \langle \Pi_n \bff{v}, \Pi_n \bff{w}\rangle_{\bb{L}^2}=  \langle \Pi_n \bff{v}, \bff{w}\rangle_{\bb{L}^2} \quad \text{for all } \bff{v,w}\in \bb{L}^2
%\end{align*}
%and {satisfies}

To prove the existence of a weak solution to \eqref{llbar}, we will use the
Faedo--Galerkin method. We first prove the following two lemmas.

\begin{lemma}\label{lem:der u2u}
	For any vector-valued function $\bff{v}:\Omega\to\bb{R}^3$, we have
	\begin{align}
	\nabla (|\bff{v}|^2 \bff{v}) 
	&= 
	2 \bff{v} \ (\bff{v}\cdot \nabla\bff{v}) 
	+ |\bff{v}|^2 \nabla \bff{v},
	\label{equ:nab un2}
	\\
	\frac{\partial\big(|\bff{v}|^2\bff{v}\big)}{\partial\bff{n}}
	&=
	2 \bff{v} \Big(\bff{v}\cdot \frac{\partial\bff{v}}{\partial\bff{n}}\Big)
	+
	|\bff{v}|^2 \frac{\partial\bff{v}}{\partial\bff{n}},
	\label{equ:nor der v2v}
	\\
	\Delta (|\bff{v}|^2 \bff{v}) 
	&= 
	2|\nabla \bff{v}|^2 \bff{v} 
	+ 2(\bff{v}\cdot \Delta \bff{v})\bff{v} 
	+ 4 \nabla \bff{v} \ (\bff{v}\cdot \nabla\bff{v})^\top
	+ |\bff{v}|^2 \Delta \bff{v}, 
	\label{equ:del v2v}
	\end{align}
	provided that the partial derivatives are well defined.
\end{lemma}
\begin{proof}
Recall the notations introduced in \eqref{equ:dot cro 2} and \eqref{equ:nab
del}. Also note that
\[
	\nabla(|\bff{v}|^2) = 2(\bff{v}\cdot\nabla\bff{v})^\top
	\quad\text{and}\quad
	\Delta(|\bff{v}|^2) = 2 |\nabla\bff{v}|^2 + 2 \bff{v}\cdot\Delta\bff{v}.
\]
Hence, it follows from the product rule that
\[
	\nabla(|\bff{v}|^2 \bff{v})
	=
	\bff{v} \big( \nabla(|\bff{v}|^2) \big)^{\top}
	+
	|\bff{v}|^2 \nabla\bff{v}
	=
	2 \bff{v} (\bff{v}\cdot\nabla\bff{v}) 
	+
	|\bff{v}|^2 \nabla\bff{v},
\]
proving \eqref{equ:nab un2}. Identity~\eqref{equ:nor der v2v} then follows
from~\eqref{equ:nab un2} and the definition of normal derivatives. 

Finally, the product rule gives
\begin{align*} 
	\Delta(|\bff{v}|^2 \bff{v})
	&=
	\Delta(|\bff{v}|^2) \bff{v} 
	+ 2 \nabla\bff{v} \ \nabla(|\bff{v}|^2) 
	+ |\bff{v}|^2 \Delta\bff{v}
	\\
	&=
	2 |\nabla\bff{v}|^2 \bff{v}
	+
	2 (\bff{v}\cdot\Delta\bff{v}) \bff{v}
	+
	4 \nabla\bff{v} \ (\bff{v}\cdot\nabla\bff{v})^\top
	+
	|\bff{v}|^2 \Delta\bff{v},
\end{align*}
proving \eqref{equ:del v2v}.
\end{proof}

\begin{lemma}\label{lem:Lip}
	For each $n\in\bb{N}$ and $\bff{v}\in\bb{V}_n$, define
	%, define $F_n^j : \bb{V}_n \to \bb{V}_n$, $j=1,\ldots,5$, by
\begin{align*}
	F_n^1(\bff{v}) &= \Delta \bff{v},
	\\
	F_n^2(\bff{v}) &= \Delta^2 \bff{v},
	\\
	F_n^3(\bff{v}) &= \Pi_n (|\bff{v}|^2 \bff{v}),
	\\
	F_n^4(\bff{v}) &= \Pi_n (\bff{v}\times \Delta \bff{v}),
	\\
	F_n^5(\bff{v}) &= \Pi_n \Delta (|\bff{v}|^2 \bff{v}).
\end{align*}
Then $F_n^j$, $j=1,\ldots,5$, are well-defined mappings from $\bb{V}_n$
into~$\bb{V}_n$. Moreover, $F_n^1$ and $F_n^2$ are globally Lipschitz while
$F_n^3$, $F_n^4$, and $F_n^5$ are locally Lipschitz.
\end{lemma}
\begin{proof}
For any $\bff{v}\in \bb{V}_n$, since $\bff{v} = \sum_{i=1}^{n}
\inpro{\bff{v}}{\bff{e}_i}_{\bb{L}^2} \bff{e}_i$, we have
\begin{align*}
-\Delta\bff{v} = \sum_{i=1}^n \lambda_i \langle \bff{v},
\bff{e}_i\rangle_{\bb{L}^2}\, \bff{e}_i \in \bb{V}_n
\quad \text{and} \quad 
\Delta^2 \bff{v} = \sum_{i=1}^n \lambda_i^2 \langle \bff{v},
\bff{e}_i\rangle_{\bb{L}^2}\, \bff{e}_i \in \bb{V}_n.
\end{align*}
Therefore, $F_n^1$ and $F_n^2$ map $\bb{V}_n$ into $\bb{V}_n$. Moreover, if the
boundary of $\Omega$ is sufficiently smooth, then the
eigenfunctions~$\bff{e}_i$, $i\in\bb{N}$, are smooth functions, and so is
$\bff{v}\in\bb{V}_n$. This implies that $|\bff{v}|^2\bff{v}$,
$\Delta(|\bff{v}|^2\bff{v})$, and $\bff{v}\times\Delta\bff{v}$ all belong to
$\bb{L}^2(\Omega)$, so that $F_n^3$, $F_n^4$, and $F_n^5$ are well defined.

We now prove the Lipschitz property of these mappings.
Using the triangle inequality, the orthonormality of $\{\bff{e}_i\}$ and
H\"{o}lder's inequality, for any $\bff{v},\bff{w}\in \bb{V}_n$ and for $j=1,2$, we have
\begin{align*} 
	\norm{F_n^j(\bff{v}) - F_n^j(\bff{w})}{\bb{L}^2}^2
	&=
	\norm{ \sum_{i=1}^{n} \lambda_i^j
		\inpro{\bff{v}-\bff{w}}{\bff{e}_i}_{\bb{L}^2} \bff{e}_i
	}{\bb{L}^2}^2
	\\
	&=
	\sum_{i=1}^{n} \lambda_i^{2j}
		\left|\inpro{\bff{v}-\bff{w}}{\bff{e}_i}_{\bb{L}^2} \right|^2
	\le 
	\Big(\sum_{i=1}^{n} \lambda_i^{2j}\Big)
        \norm{\bff{v}-\bff{w}}{\bb{L}^2}^2.
\end{align*}
Hence, $F_n^1$ and $F_n^2$ are globally Lipschitz.

Next, it follows from \eqref{equ:Pin} that
\begin{align*}
\|F_n^3 (\bff{v}) -F_n^3 (\bff{w})\|_{\bb{L}^2} &\le  \||\bff{v}|^2 \bff{v}- |\bff{w}|^2 \bff{w}\|_{\bb{L}^2} \\
&\leq \| |\bff{v}|^2 (\bff{v}-\bff{w})\|_{\bb{L}^2} + \|(\bff{v}-\bff{w})\cdot (\bff{v}+\bff{w})\bff{w}\|_{\bb{L}^2} \\
&\leq \left(\|\bff{v}\|_{\bb{L}^\infty}^2 + \|\bff{v}+\bff{w}\|_{\bb{L}^\infty} \|\bff{w}\|_{\bb{L}^\infty} \right) \|\bff{v}-\bff{w}\|_{\bb{L}^2},
\end{align*}
where we used the fact that all norms are equivalent in the finite dimensional
subspace $\bb{V}_n$. This shows that $F_n^3$ is locally Lipschitz.

Similarly, it follows from \eqref{equ:Pin} that
\begin{align*}
\|F_n^4 (\bff{v}) -F_n^4 (\bff{w})\|_{\bb{L}^2} &\le \|\bff{v}\times \Delta \bff{v}- \bff{w}\times \Delta\bff{w} \|_{\bb{L}^2} \\
&\leq \|\bff{v} \times (\Delta \bff{v}-\Delta \bff{w})\|_{\bb{L}^2} + \|(\bff{v}-\bff{w}) \times \Delta \bff{w}\|_{\bb{L}^2}\\
&\leq \|\bff{v}\|_{\bb{L}^\infty} \|F_n^2(\bff{v})-F_n^2(\bff{w})\|_{\bb{L}^2} + \|\bff{v}-\bff{w}\|_{\bb{L}^2} \|\Delta \bff{w}\|_{\bb{L}^\infty}.
\end{align*}
Since $F_n^2$ is Lipschitz, we deduce that $F_n^4$ is locally Lipschitz.

Finally, note that if $\bff{v}\in\bb{V}_n$, then
$\partial\bff{v}/\partial\bff{n}=0$. Thus \eqref{equ:nor der v2v} implies
$\partial\big(|\bff{v}|^2\bff{v}\big)/\partial\bff{n}=0$,
which allows us to use integration by parts to obtain
\[
\inpro{\Delta\big(|\bff{v}|^2\bff{v}\big)}{\bff{w}}
	=
	- \inpro{\nabla\big(|\bff{v}|^2\bff{v}\big)}{\nabla\bff{w}}
	=
	\inpro{|\bff{v}|^2\bff{v}}{\Delta\bff{w}}
	\quad\forall\bff{v}, \ \bff{w} \in \bb{V}_n.
\]
Therefore, for any $\bff{v},\bff{w}\in\bb{V}_n$, we can use the definition of $\Pi_n$
and integration by parts again to have
\begin{align*} 
	\inpro{\Pi_n \Delta\big(|\bff{v}|^2\bff{v}\big)}{\bff{w}}
	&=
	\inpro{\Delta\big(|\bff{v}|^2\bff{v}\big)}{\bff{w}}
	%=
	%- \inpro{\nabla\big(|\bff{v}|^2\bff{v}\big)}{\nabla\bff{w}}
	=
	\inpro{|\bff{v}|^2\bff{v}}{\Delta\bff{w}}
	%\\
	%&=
	=
	\inpro{\Pi_n\big(|\bff{v}|^2\bff{v}\big)}{\Delta\bff{w}}
	=
	\inpro{\Delta\Pi_n\big(|\bff{v}|^2\bff{v}\big)}{\bff{w}}.
\end{align*}
This means $\Delta$ and $\Pi_n$ commute, so that 
$F_n^5(\bff{v})= F_n^1 \circ F_n^3 (\bff{v})$. Since $F_n^1$ is
Lipschitz and $F_n^3$ is locally Lipschitz, we have that $F_n^5$ is locally
Lipschitz as well. 
This completes the proof.
\end{proof}

The Faedo--Galerkin method seeks to approximate the solution to~\eqref{llbar} by
$\bff{u}_n(t) \in \bb{V}_n$ satisfying the equation
%\multlinegap=2cm
\begin{equation}\label{approx}
\left\{
\begin{alignedat}{2}
&\frac{\partial \bff{u}_n}{\partial t} 
- \beta_1 \Delta \bff{u}_n 
+ \beta_2 \Delta^2 \bff{u}_n 
	-\beta_3 \Pi_n((1-|\bff{u_n}|^2) \bff{u}_n) &&
\\
&\quad +\beta_4 \Pi_n (\bff{u}_n \times \Delta \bff{u}_n) 
- \beta_5 \Pi_n(\Delta (|\bff{u}_n|^2 \bff{u}_n))
	= \bff{0} && \quad\text{in $(0,T)\times\Omega$,}
\\
		  &\bff{u}_n(0) = \bff{u}_{0n}
		     && \quad\text{in $\Omega$,}
\end{alignedat}
\right.
\end{equation}
where $\bff{u}_{0n}\in \bb{V}_n$ is an approximation of $\bff{u}_0$ {such
that if $\bff{u}_0 \in \bb{H}^r$, then
\begin{equation}\label{equ:un0}
	\norm{\bff{u}_{0n}}{\bb{H}^r}
	\lesssim
	\norm{\bff{u}_{0}}{\bb{H}^r}
\end{equation}
}

Lemma~\ref{lem:Lip} assures us that all the terms in~\eqref{approx} are well
defined. Moreover, the existence of solutions to the above ordinary differential
equation in $\bb{V}_n$ is guaranteed by this lemma and the Cauchy--Lipschitz
theorem.

We now prove some a priori estimates for the solution of \eqref{approx}. First
we need the following results.

\begin{lemma}\label{equivnorm}
Let $\Omega \subset \bb{R}^d$ be an open bounded domain with $C^\infty$-boundary and $\epsilon>0$ be
given. Then there exists a positive constant $C$ such that the following
inequalities hold:
\begin{enumerate}
\renewcommand{\labelenumi}{\theenumi}
\renewcommand{\theenumi}{{\rm (\roman{enumi})}}
\item for any $\bff{v} \in \bb{L}^2(\Omega)$ such that $\Delta \bff{v}\in \bb{L}^2(\Omega)$ satisfying $\displaystyle
	\frac{\partial \bff{v}}{\partial\bff{n}}=0$ on $\partial \Omega$,
\begin{align}
	\norm{\bff{v}}{\bb{H}^2}^2 &\leq C \left( \norm{\bff{v}}{\bb{L}^2}^2 + \norm{\Delta \bff{v}}{\bb{L}^2}^2 \right), \label{eq1} \\
	\norm{\nabla \bff{v}}{\bb{L}^2}^2 &\leq C \norm{\bff{v}}{\bb{L}^2}^2 + \varepsilon \norm{\Delta \bff{v}}{\bb{L}^2}^2, \label{eq2}
\end{align}

\item for any $\bff{v} \in \bb{H}^1(\Omega)$ such that $\Delta \bff{v}\in \bb{H}^1$ satisfying $\displaystyle
	\frac{\partial \bff{v}}{\partial\bff{n}}=0$ on $\partial \Omega$,
\begin{align} 
	\norm{\Delta \bff{v}}{\bb{L}^2}^2 
	&\leq 
	\norm{\nabla \bff{v}}{\bb{L}^2} 
	\norm{\nabla \Delta \bff{v}}{\bb{L}^2}, 
	\label{eq3} 
	\\
	\norm{\bff{v}}{\bb{H}^3}^2 
	&\leq 
	C \left( 
	\norm{\bff{v}}{\bb{L}^2}^2 
	+ 
	\norm{\nabla \bff{v}}{\bb{L}^2}^2 
	+ 
	\norm{\nabla \Delta \bff{v}}{\bb{L}^2}^2 
	\right), \label{eq5}
\end{align}

\item for any $\bff{v} \in \bb{H}^2(\Omega)$ such that $\Delta^2 \bff{v} \in \bb{L}^2(\Omega)$ satisfying $\displaystyle
	\frac{\partial \bff{v}}{\partial\bff{n}}=\frac{\partial \Delta \bff{v}}{\partial
	\bff{n}}=0$ on $\partial \Omega$,
\begin{align}
	\norm{\nabla \Delta \bff{v}}{\bb{L}^2}^2 
	&\leq 
	\norm{\Delta \bff{v}}{\bb{L}^2}
	\norm{\Delta^2 \bff{v}}{\bb{L}^2},
	\label{eq4}
	\\
	\norm{\bff{v}}{\bb{H}^4}^2
	&\leq
	C \left(\norm{\bff{v}}{\bb{L}^2}^2
	+
	\norm{\Delta \bff{v}}{\bb{L}^2}^2
	+
	\norm{\Delta^2 \bff{v}}{\bb{L}^2}^2 \right),
	\label{eq6}
\end{align}

\item for any $\bff{v} \in \bb{H}^3(\Omega)$ such that $\Delta^2 \bff{v} \in \bb{H}^1(\Omega)$ satisfying $\displaystyle
\frac{\partial \bff{v}}{\partial\bff{n}}=\frac{\partial \Delta \bff{v}}{\partial
	\bff{n}}=0$ on $\partial \Omega$,
\begin{align}
	\norm{\bff{v}}{\bb{H}^5}^2
	&\leq
	C \left(\norm{\bff{v}}{\bb{L}^2}^2
	+
	\norm{\nabla \bff{v}}{\bb{L}^2}^2
	+
	\norm{\nabla \Delta \bff{v}}{\bb{L}^2}^2
	+
	\norm{\nabla \Delta^2 \bff{v}}{\bb{L}^2}^2 \right).
	\label{eq8}
\end{align}

\end{enumerate}
\end{lemma}

\begin{proof}
Inequality \eqref{eq1} follows from the standard elliptic regularity result with Neumann boundary data \cite[Corollary 2.2.2.6]{Gri11}.
Next, using integration by parts, H\"{o}lder's inequality, and Young's inequality, we obtain
\begin{align*}
	\norm{\nabla \bff{v}}{\bb{L}^2}^2 
	= \inpro{\nabla \bff{v}}{\nabla \bff{v}}_{\bb{L}^2}
	= -\inpro{\bff{v}}{\Delta \bff{v}}
	\leq \norm{\bff{v}}{\bb{L}^2} \norm{\Delta \bff{v}}{\bb{L}^2}
	\leq C \norm{\bff{v}}{\bb{L}^2}^2 + \varepsilon \norm{\Delta \bff{v}}{\bb{L}^2}^2,
\end{align*}
proving \eqref{eq2}.

Similarly, we have
\begin{align*}
	\norm{\Delta \bff{v}}{\bb{L}^2}^2
	= -\inpro{\nabla \Delta \bff{v}}{\nabla \bff{v}}
	\leq \norm{\nabla \Delta \bff{v}}{\bb{L}^2} \norm{\nabla \bff{v}}{\bb{L}^2}
\end{align*}
and
\begin{align*}
	\norm{\nabla \Delta \bff{v}}{\bb{L}^2}^2
	= - \inpro{\Delta \bff{v}}{\Delta^2 \bff{v}}
	\leq \norm{\Delta \bff{v}}{\bb{L}^2} \norm{\Delta^2 \bff{v}}{\bb{L}^2},
\end{align*}
proving \eqref{eq3} and \eqref{eq4}.

Next, applying the higher-order regularity result \cite[Remark 2.5.1.2]{Gri11} to the elliptic operator $(-\Delta + I)$ yields
\begin{align*}
	\norm{\bff{v}}{\bb{H}^3}^2
	&\leq
	C \big( \norm{\bff{v}}{\bb{H}^1}^2
	+ \norm{\Delta \bff{v}}{\bb{H}^1}^2 \big)
	\leq
	C \big( \norm{\bff{v}}{\bb{L}^2}^2
	+
	\norm{\nabla \bff{v}}{\bb{L}^2}^2
	+
	\norm{\nabla \Delta \bff{v}}{\bb{L}^2}^2 \big),
\end{align*}
proving \eqref{eq5}. 
Furthermore, by the same elliptic regularity result,
\begin{align*}
	\norm{\bff{v}}{\bb{H}^4}^2
	\leq
	C \big(\norm{\bff{v}}{\bb{L}^2}^2
	+
	\norm{\Delta \bff{v}}{\bb{H}^2}^2 \big)
	\leq
	C \big( \norm{\bff{v}}{\bb{L}^2}^2
	+
	\norm{\Delta \bff{v}}{\bb{L}^2}^2
	+
	\norm{\Delta^2 \bff{v}}{\bb{L}^2}^2 \big),
\end{align*}
where we applied \eqref{eq1} to $\Delta \bff{v}$, noting that the assumptions are satisfied. This proves \eqref{eq6}.

Similarly, we have
\begin{align*}
	\norm{\bff{v}}{\bb{H}^5}^2 
	\leq
	C \big(\norm{\Delta \bff{v}}{\bb{H}^3}^2 
	+ \norm{\bff{v}}{\bb{H}^3}^2 \big)
	\leq
	C \big( 
	\norm{\bff{v}}{\bb{L}^2}^2
	+ \norm{\nabla \bff{v}}{\bb{L}^2}^2
	+ \norm{\nabla \Delta \bff{v}}{\bb{L}^2}^2
	+ \norm{\nabla \Delta^2 \bff{v}}{\bb{L}^2}^2 \big),
\end{align*}
where we applied \eqref{eq5} to $\Delta \bff{v}$, thus proving \eqref{eq8}.
\end{proof}
{
\begin{remark}
	Elliptic regularity result in $\bb{H}^2$ \eqref{eq1} holds more
	generally for a domain with $C^{1,1}$-boundary \cite[Remark
	2.2.2.6]{Gri11} or a convex polygonal domain \cite[Theorem
	4.3.1.4]{Gri11}. The $\bb{H}^r$-regularity results also hold for a
	domain with $C^{r-1,1}$-boundary (see \cite[Remark 2.5.1.2]{Gri11}).
\end{remark}
}
{
\begin{lemma}
	Let $\bff{k}=(k_1,\cdots,k_d)$ be a multi-index and define the operator
	\begin{align*}
		D^{\bff{k}} :=
		\frac{\partial^{|\bff{k}|}}{\partial_{x_1}^{k_1}\cdots\partial_{x_d}^{k_d}} 
	\end{align*}
	where $|\bff{k}| = k_1 + \cdots + k_d$.
	Then the following inequalities hold:
	\begin{enumerate}
	\renewcommand{\labelenumi}{\theenumi}
	\renewcommand{\theenumi}{{\rm (\roman{enumi})}}
		\item for any $\bff{u}, \bff{v} \in \bb{H}^s(\Omega)$, where $s>d/2$,
			\begin{align} \label{eq7}
				\norm{\abs{\bff{u}} \abs{\bff{v}}}{H^s}
				\leq
				C \norm{\bff{u}}{\bb{H}^s}
				\norm{\bff{v}}{\bb{H}^s}.	
			\end{align}
		\item for any $\bff{u}$, $\bff{v} \in \bb{H}^{2} \cap
			\bb{H}^{|\bff{k}|}$,
			\begin{align}\label{equ:D k u2u-v2v}
				\norm{D^{\bff{k}} (|\bff{u}|^2 \bff{u}-
				|\bff{v}|^2 \bff{v})}{\bb{L}^2}
				\lesssim
				\begin{cases}
					\big( \norm{\bff{u}}{\bb{L}^\infty}^2 + \norm{\bff{v}}{\bb{L}^\infty}^2 \big) 
					\norm{\bff{u}-\bff{v}}{\bb{L}^2}
					\ &
					\text{if $|\bff{k}|=0$},
					\\[1ex]
					\big( \norm{\bff{u}}{\bb{H}^1} + \norm{\bff{v}}{\bb{H}^1} \big)
					\big( \norm{\bff{u}}{\bb{H}^2} + \norm{\bff{v}}{\bb{H}^2} \big) 
					\norm{\bff{u}-\bff{v}}{\bb{H}^1} 
					\ &
					\text{if $|\bff{k}|=1$},
					\\[1ex]
					\big( \norm{\bff{u}}{\bb{H}^{|\bff{k}|}}^2 +
					\norm{\bff{v}}{\bb{H}^{|\bff{k}|}}^2 \big)
					\norm{\bff{u}-\bff{v}}{\bb{H}^{|\bff{k}|}}
					\ &
					\text{if $|\bff{k}|\geq 2$}.
				\end{cases}
			\end{align}
		\item for any $\bff{u}$, $\bff{v} \in \bb{H}^{|\bff{k}|} \cap
			\bb{L}^\infty$,
			\begin{align}\label{equ:u cross Dk u}
				\norm{\bff{u} \times D^{\bff{k}} \bff{u} - \bff{v} \times D^{\bff{k}} \bff{v} }{\bb{L}^2}
				\lesssim
					\norm{\bff{u}}{\bb{L}^\infty} \norm{D^{\bff{k}} (\bff{u}-\bff{v})}{\bb{L}^2}
					+
					\norm{\big(\bff{u}-\bff{v}\big)D^{\bff{k}}
					\bff{v}}{\bb{L}^2}.
			\end{align}
	\end{enumerate}
\end{lemma}

\begin{proof}
Firstly, \eqref{eq7} follows from %Theorem 4.39 in \cite{AdaFou03}.
$\norm{|\bff{v}||\bff{w}|}{H^s} \lesssim \norm{\bff{v}}{\bb{L}^\infty}\norm{\bff{w}}{\bb{H}^s}
+ \norm{\bff{v}}{\bb{H}^s}\norm{\bff{w}}{\bb{L}^\infty}$
and the Sobolev embedding.
	Next, we prove \eqref{equ:D k u2u-v2v}. For the case $|\bff{k}|=0$, by H\"{o}lder's inequality we have
	\begin{align*}
		\norm{ |\bff{u}|^2 \bff{u} - |\bff{v}|^2 \bff{v}}{\bb{L}^2} 
		&\leq
		\norm{|\bff{u}|^2 (\bff{u}-\bff{v})}{\bb{L}^2}
		+
		\norm{|\bff{u}+\bff{v}||\bff{u}-\bff{v}|| \bff{v}}{\bb{L}^2}
		\\
		&\leq 
		\big( \norm{ \bff{u}}{\bb{L}^\infty}^2 
		+ 
		\norm{\bff{v}}{\bb{L}^\infty}
		\norm{\bff{u} + \bff{v} }{\bb{L}^\infty}
		\big) 
		\norm{\bff{u}-\bff{v}}{\bb{L}^2}
		\\
		&\lesssim 
		\big( \norm{\bff{u}}{\bb{L}^\infty}^2 + \norm{\bff{v}}{\bb{L}^\infty}^2 \big)
		\norm{\bff{u}-\bff{v}}{\bb{L}^2}.
	\end{align*}
	For the case $|\bff{k}|=1$, note that
	\begin{align*} 
		\big|D^{\bff{k}} (|\bff{u}|^2 \bff{u}- |\bff{v}|^2 \bff{v})\big|
		&\leq
		\big|D^{\bff{k}}(|\bff{u}|^2) (\bff{u}-\bff{v})\big|
		+ \big||\bff{u}|^2 D^{\bff{k}} (\bff{u}-\bff{v})\big|
		+ \big|D^{\bff{k}}(|\bff{u}|^2 - |\bff{v}|^2) \bff{v}\big|
		\\
		&\quad
		+ \big|(|\bff{u}|^2-|\bff{v}|^2) D^{\bff{k}} \bff{v}\big|
		\\
		&\leq
		2|\bff{u}| \big|D^{\bff{k}}(\bff{u})\big| \big|\bff{u}-\bff{v}\big|
		+ 
		|\bff{u}|^2 \big|D^{\bff{k}} (\bff{u}-\bff{v})\big|
		+
		\big|D^{\bff{k}}\bff{u}+D^{\bff{k}}\bff{v}\big|
		|\bff{u}-\bff{v}| |\bff{v}|
		\\
		&\quad
		+
		|\bff{u}+\bff{v}| 
		\big|D^{\bff{k}}\bff{u}-D^{\bff{k}}\bff{v}\big|
		|\bff{v}|
		+
		|\bff{u}+\bff{v}| 
		|\bff{u}+\bff{v}| 
		\big|D^{\bff{k}}\bff{v}\big|.
	\end{align*}
	Using this, we have
	\begin{align*}
		&\norm{D^{\bff{k}}(|\bff{u}|^2 \bff{u}- |\bff{v}|^2 \bff{v})}{\bb{L}^2}
		\\
		&\le
		2\norm{|\bff{u}| \big|D^{\bff{k}}(\bff{u})\big| |\bff{u}-\bff{v}|}{\bb{L}^2}
		+ 
		\norm{|\bff{u}|^2 D^{\bff{k}} (\bff{u}-\bff{v})}{\bb{L}^2}
		+
		\norm{|D^{\bff{k}} \bff{u}+ D^{\bff{k}} \bff{v}| \ |\bff{u}-\bff{v}| \ |\bff{v}|}{L^2} 
		\\
		&\quad
		+
		\norm{ |\bff{u}+\bff{v}| \ |D^{\bff{k}}\bff{u}-D^{\bff{k}} \bff{v}| \ |\bff{v}|}{L^2}
		+
		\norm{ |\bff{u}+\bff{v}| \ |\bff{u}- \bff{v}| \ |D^{\bff{k}}\bff{v}|}{L^2}
		\\
		&\lesssim
		\big( \norm{\bff{u}}{\bb{L}^6} \norm{D^{\bff{k}}
			\bff{u}}{\bb{L}^6} 
		+ \norm{\bff{u}}{\bb{L}^6} \norm{D^{\bff{k}}
		\bff{v}}{\bb{L}^6}
		+ \norm{\bff{v}}{\bb{L}^6} \norm{D^{\bff{k}}
		\bff{u}}{\bb{L}^6}
		+ \norm{\bff{v}}{\bb{L}^6} \norm{D^{\bff{k}}
	\bff{v}}{\bb{L}^6} \big)
		\norm{\bff{u}-\bff{v}}{\bb{L}^2} 
		\\
		&\quad
		+
		\big( \norm{\bff{u}}{\bb{L}^\infty}^2 
		+ \norm{\bff{v}}{\bb{L}^\infty}^2 \big)
		\norm{D^{\bff{k}} \bff{u}-D^{\bff{k}} \bff{v}}{\bb{L}^2} 
		\\
		&\lesssim
		\big( \norm{\bff{u}}{\bb{H}^1} \norm{D^{\bff{k}} \bff{u}}{\bb{H}^1} 
		+ \norm{\bff{u}}{\bb{H}^1} \norm{D^{\bff{k}} \bff{v}}{\bb{H}^1}
		+ \norm{\bff{v}}{\bb{H}^1} \norm{D^{\bff{k}} \bff{u}}{\bb{H}^1}
		+ \norm{\bff{v}}{\bb{H}^1} \norm{D^{\bff{k}} \bff{v}}{\bb{H}^1} \big)
		\norm{\bff{u}-\bff{v}}{\bb{H}^1}
		\\
		&\quad
		+
		\big( \norm{\bff{u}}{\bb{H}^1} \norm{\bff{u}}{\bb{H}^2}
		+ \norm{\bff{v}}{\bb{H}^1} \norm{\bff{v}}{\bb{H}^2} \big)
		\norm{\bff{u}-\bff{v}}{\bb{H}^1}
		\\
		&\lesssim
		\big( \norm{\bff{u}}{\bb{H}^1} \norm{\bff{u}}{\bb{H}^2} 
		+ \norm{\bff{u}}{\bb{H}^1} \norm{\bff{v}}{\bb{H}^2} 
		+ \norm{\bff{v}}{\bb{H}^1} \norm{\bff{u}}{\bb{H}^2}
		+ \norm{\bff{v}}{\bb{H}^1} \norm{\bff{v}}{\bb{H}^2}
		\big) 
		\norm{\bff{u}-\bff{v}}{\bb{H}^1}
		\\
		&=
		\big( \norm{\bff{u}}{\bb{H}^1} + \norm{\bff{v}}{\bb{H}^1} \big)
		\big( \norm{\bff{u}}{\bb{H}^2} + \norm{\bff{v}}{\bb{H}^2} \big) 
		\norm{\bff{u}-\bff{v}}{\bb{H}^1},
	\end{align*}
where in the third step we used Gagliardo--Nirenberg's inequality
$\norm{\bff{u}}{\bb{L}^\infty}^2 \lesssim \norm{\bff{u}}{\bb{H}^1}
\norm{\bff{v}}{\bb{H}^2}$ (valid for $d\in \{1,2,3\}$) and the Sobolev embedding
$\bb{H}^1 \subset \bb{L}^6$.

For the case $|\bff{k}|\geq 2$, we use H\"{o}lder's inequality and \eqref{eq7} to obtain
	\begin{align*}
		\norm{D^{\bff{k}} (|\bff{u}|^2 \bff{u}- |\bff{v}|^2 \bff{v})}{\bb{L}^2} 
		&\leq
		\norm{D^{\bff{k}} (|\bff{u}|^2 \bff{u}- |\bff{v}|^2 \bff{v})}{\bb{L}^2}
		\leq
		\norm{|\bff{u}|^2 \bff{u}- |\bff{v}|^2 \bff{v}}{\bb{H}^{|\bff{k}|}}
		\\
		&\leq 
		\big( \norm{ \bff{u}}{\bb{H}^{|\bff{k}|}}^2 
		+ 
		\norm{\bff{v}}{\bb{H}^{|\bff{k}|}}
		\norm{\bff{u} + \bff{v} }{\bb{H}^{|\bff{k}|}}
		\big) 
		\norm{\bff{u}-\bff{v}}{\bb{H}^{|\bff{k}|}}
		\\
		&\lesssim 
		\big( 
		\norm{\bff{u}}{\bb{H}^{|\bff{k}|}}^2 
		+ \norm{\bff{v}}{\bb{H}^{|\bff{k}|}}^2 
		\big)
		\norm{\bff{u}-\bff{v}}{\bb{H}^{|\bff{k}|}}.
	\end{align*}
This completes the proof of \eqref{equ:D k u2u-v2v}. The proof of~\eqref{equ:u
cross Dk u} is obvious and is omitted.
\end{proof}
}
We now use the above lemmas to derive a priori estimates on the Galerkin
solution~$\bff{u}_n$.

\begin{proposition}\label{pro:Delta un}
	Let $T>0$ be arbitrary {and assume that $\bff{u}_0\in\bb{L}^2$}. For each
	$n\in \bb{N}$ and all $t\in [0,T]$,
\begin{align}\label{equ:un nab del}
\norm{\bff{u}_n(t)}{\bb{L}^2}^2 
&+ 
{
\int_0^t \norm{\nabla \bff{u}_n(s)}{\bb{L}^2}^2 \ds 
}
+
\int_0^t \norm{\Delta \bff{u}_n(s)}{\bb{L}^2}^2 \ds 
+  
\int_0^t \norm{\bff{u}_n(s)}{\bb{L}^4}^4 \ds 
\nonumber\\
&+ 
\int_0^t \norm{ |\bff{u}_n(s)| |\nabla\bff{u}_n(s)| }{\bb{L}^2}^2 \ds 
+  
\int_0^t \norm{\bff{u}_n(s)\cdot \nabla \bff{u}_n(s)}{\bb{L}^2}^2 \ds 
\lesssim 
{
\norm{\bff{u}_0}{\bb{L}^2}^2,
}
\end{align}
where the constant depends on $T$, but is independent of $n$.
\end{proposition}

\begin{proof}
Taking the inner product of \eqref{approx} with $\bff{u}_n(t)$,
integrating by parts with respect to $\bff{x}$ (noting \eqref{equ:nor der
v2v}) and using the identity~\eqref{equ:nab un2}, we obtain, for any $\epsilon > 0$,
\begin{align*}
\frac{1}{2} \ddt \norm{\bff{u}_n}{\bb{L}^2}^2 
&+\beta_2 \norm{\Delta \bff{u}_n}{\bb{L}^2}^2 
+\beta_3 \norm{\bff{u}_n}{\bb{L}^4}^4 
+ 2 \beta_5 \norm{\bff{u}_n \cdot \nabla \bff{u}_n}{\bb{L}^2}^2 
+ \beta_5 \norm{ |\bff{u}_n| |\nabla \bff{u}_n| }{\bb{L}^2}^2 
\\
&\le \beta_3 \norm{\bff{u}_n}{\bb{L}^2}^2 
+ |\beta_1| \norm{\nabla \bff{u}_n}{\bb{L}^2}^2
\le
C \norm{\bff{u}_n}{\bb{L}^2}^2 
+ \epsilon \norm{\Delta \bff{u}_n}{\bb{L}^2}^2,
\end{align*}
where in the last step we used \eqref{eq2}. Rearranging the above
equation, choosing sufficiently small $\epsilon$, and integrating over $(0,t)$,
we deduce
\begin{multline*}
	\norm{\bff{u}_n(t)}{\bb{L}^2}^2 
	+ 
	\int_0^t \norm{\Delta \bff{u}_n(s)}{\bb{L}^2}^2 \ds 
	+
	\int_0^t \norm{\bff{u}_n(s)}{\bb{L}^4}^4 \ds 
	+ 
	\int_0^t \norm{ |\bff{u}_n(s)| |\nabla \bff{u}_n(s)| }{\bb{L}^2}^2 \ds 
	\\
	+
	\int_0^t \norm{\bff{u}_n(s) \cdot \nabla \bff{u}_n(s)}{\bb{L}^2}^2 \ds 
	\lesssim 
	\norm{\bff{u}_n(0)}{\bb{L}^2}^2 
	+ 
	\int_0^t \norm{\bff{u}_n(s)}{\bb{L}^2}^2 \ds.
\end{multline*}
{
Invoking Gronwall's inequality yields the required estimate for all the terms on
the left-hand side of~\eqref{equ:un nab del}, except $\int_0^t \norm{\nabla
\bff{u}_n(s)}{\bb{L}^2}^2 \ds$. The bound for this term follows
from~\eqref{eq2}, completing the proof of the proposition.
}
\end{proof}

\begin{proposition}\label{pro:dt un H-2}
	{Under the assumption of Proposition~\ref{pro:Delta un}}, we have
	\begin{align*}
			\norm{\partial_t \bff{u}_n}{L^2(0,T;\bb{H}^{-2})} 
			\lesssim
			\norm{\bff{u}_0}{\bb{L}^2}^2,
	\end{align*}
	where the constant depends on $T$ but is independent of $n$.
\end{proposition}
	
\begin{proof}
	Taking the inner product of \eqref{approx} with $\bff{\varphi}\in \bb{H}^2$ such that $\norm{\bff{\varphi}}{\bb{H}^2} \leq 1$ and integrating by parts with respect to $\bff{x}$, we have
		\begin{align*}
			\inpro{\partial_t \bff{u}_n}{\bff{\varphi}}_{\bb{L}^2}
			&=
			-\beta_1 \inpro{\nabla \bff{u}_n}{\nabla \bff{\varphi}}_{\bb{L}^2}
			-
			\beta_2 \inpro{\Delta \bff{u}_n}{\Delta \bff{\varphi}}_{\bb{L}^2}
			+
			\beta_3 \inpro{\bff{u}_n}{\bff{\varphi}}_{\bb{L}^2}
			-
			\beta_3 \inpro{|\bff{u}_n|^2 \bff{u}_n}{\bff{\varphi}}_{\bb{L}^2}
			\\
			&\quad
			+
			\beta_4 \inpro{\bff{u}_n \times \nabla \bff{u}_n}{\nabla \bff{\varphi}}_{\bb{L}^2}
			+
			\beta_5 \inpro{\nabla(|\bff{u}_n|^2 \bff{u}_n)}{\nabla \bff{\varphi}}_{\bb{L}^2}.
		\end{align*}
	It follows from this equation and H\"{o}lder's inequality that
		\begin{align*}
			\big| \inpro{\partial_t \bff{u}_n}{\bff{\varphi}}_{\bb{L}^2} \big|
			&\lesssim
			\norm{\bff{u}_n}{\bb{H}^1} \norm{\bff{\varphi}}{\bb{H}^1}
			+
			\norm{\bff{u}_n}{\bb{H}^2} \norm{\bff{\varphi}}{\bb{H}^2}
			+
			\norm{\bff{u}_n}{\bb{L}^2} \norm{\bff{\varphi}}{\bb{L}^2}
			+
			\norm{\bff{u}_n}{\bb{L}^2} \norm{\bff{u}_n}{\bb{L}^6}^2 \norm{\bff{\varphi}}{\bb{L}^6}
			\\
			&\quad
			+
			\norm{\bff{u}_n}{\bb{L}^4} \norm{\nabla \bff{u}_n}{\bb{L}^2}
			\norm{\nabla \bff{\varphi}}{\bb{L}^4}
			+
			\norm{\bff{u}_n}{\bb{L}^2} \norm{\bff{u}_n}{\bb{L}^6}
			\norm{\nabla \bff{u}_n}{\bb{L}^6} \norm{\nabla \bff{\varphi}}{\bb{L}^6},
		\end{align*}
	where we used the Sobolev embedding $\bb{H}^1 \subset \bb{L}^6$. Therefore, integrating over $(0,t)$ and noting that $\norm{\bff{\varphi}}{\bb{H}^2} \leq 1$, we obtain
		\[
		\int_0^t \big| \inpro{\partial_s \bff{u}_n(s)}{\bff{\varphi}}_{\bb{L}^2} \big|^2 \,\ds
		\lesssim
		\int_0^t \norm{\bff{u}_n}{\bb{H}^2}^2 \ds \lesssim 1
		\]
	by Proposition \ref{pro:Delta un}. Taking supremum over the set
	$\{\varphi\in \bb{H}^2: \norm{\bff{\varphi}}{\bb{H}^2} \leq 1\}$ and
	noting that
	\[
		\norm{\partial_t\bff{u}_n}{\bb{H}^{-2}}
		\le
		\sup_{\bff{\varphi}\in\bb{H}^2}
		\frac{\big|{\inpro{\partial_t\bff{u}_n}{\bff{\varphi}}_{\bb{L}^2}}\big|}
		{\norm{\bff{\varphi}}{\bb{H}^2}}, 
	\]
	we obtain the required estimate.
	\end{proof}

\begin{proposition}\label{pro:nab delta un}
Let $T>0$ be arbitrary {and assume that $\bff{u}_0\in\bb{H}^1$.} Then there exists
$T^\ast >0$ such that for $n\in \bb{N}$ and $t\in [0,T^\ast]$, we have
\begin{align*}
	\norm{\nabla \bff{u}_n(t)}{\bb{L}^2}^2 
	&+ 
	\int_0^t \norm{\nabla \Delta \bff{u}_n(s)}{\bb{L}^2}^2 \ds  
	+ 
	\int_0^t \norm{ |\bff{u}_n(s)| |\Delta\bff{u}_n(s)| }{\bb{L}^2}^2 \ds 
	+ 
	\int_0^t \norm{\bff{u}_n(s)\cdot \Delta \bff{u}_n(s)}{\bb{L}^2}^2 \ds 
	\\
	&\lesssim
	\norm{\bff{u}_n(0)}{\bb{H}^1}^2
	{\lesssim
	\norm{\bff{u}_0}{\bb{H}^1}^2.}
\end{align*}
The constant depends on $T^\ast$, but is independent of $n$. Here,
\[
\begin{cases}
T^\ast = T \quad & \text{for $d=1,2$,}
\\
T^\ast \le T \quad & \text{for $d=3$,}
\end{cases}
\]
{where $T^\ast = T^\ast(\norm{\bff{u}_0}{\bb{H}^1})$.}
\end{proposition}

\begin{proof}
Taking the inner product of \eqref{approx} with $-\Delta \bff{u}_n(t)$ and
integrating by parts with respect to $\bff{x}$, we have
\begin{multline}\label{equ:Del un}
	\frac{1}{2} \ddt \norm{\nabla \bff{u}_n}{\bb{L}^2}^2 
	+\beta_1 \norm{\Delta \bff{u}_n}{\bb{L}^2}^2 
	+ \beta_2 \norm{\nabla \Delta \bff{u}_n}{\bb{L}^2}^2  
	- \beta_3 \norm{\nabla \bff{u}_n}{\bb{L}^2}^2
	\\
	+\beta_3 \inpro{\nabla (|\bff{u}_n|^2 \bff{u}_n)}{\nabla \bff{u}_n}_{\bb{L}^2} +
	\beta_5 \inpro{\Delta (|\bff{u}_n|^2 \bff{u}_n)}{\Delta \bff{u}_n}_{\bb{L}^2} =
	0.
\end{multline}
It follows from \eqref{equ:nab un2} and \eqref{equ:del v2v} that
\[
	\inpro{\nabla (|\bff{u}_n|^2 \bff{u}_n)}{\nabla \bff{u}_n}_{\bb{L}^2}
	=
	2 \norm{\bff{u}_n\cdot\nabla\bff{u}_n}{\bb{L}^2}^2
	+
	\norm{|\bff{u}_n||\nabla\bff{u}_n|}{\bb{L}^2}^2
\]
and
\begin{align*} 
\inpro{\Delta (|\bff{u}_n|^2 \bff{u}_n)}{\Delta \bff{u}_n}_{\bb{L}^2}
&=
2 \inpro{|\nabla\bff{u}_n|^2\bff{u}_n}{\Delta\bff{u}_n}_{\bb{L}^2}
+
2 \norm{\bff{u}_n\cdot\Delta\bff{u}_n}{\bb{L}^2}^2
%\inpro{(\bff{u}_n\cdot\Delta\bff{u}_n)\bff{u}_n}{\Delta\bff{u}_n}_{\bb{L}^2}
\\
&\quad+
4
\inpro{\nabla\bff{u}_n(\bff{u}_n\cdot\nabla\bff{u}_n)^\top}{\Delta\bff{u}_n}_{\bb{L}^2}
+
\norm{|\bff{u}_n||\Delta\bff{u}_n|}{\bb{L}^2}^2.
%\inpro{|\bff{u}_n|^2 \Delta\bff{u}_n}{\Delta\bff{u}_n}_{\bb{L}^2}
\end{align*}
Therefore, after rearranging the terms in \eqref{equ:Del un}, we obtain
\begin{align}\label{equ:mul Del un}
& \frac{1}{2} \ddt \norm{\nabla \bff{u}_n}{\bb{L}^2}^2 
+ \beta_2 \norm{\nabla \Delta \bff{u}_n}{\bb{L}^2}^2 
+ 2\beta_3 \norm{\bff{u}_n \cdot \nabla \bff{u}_n}{\bb{L}^2}^2 
+ \beta_3 \norm{ |\bff{u}_n| |\nabla \bff{u}_n| }{\bb{L}^2}^2 
\nonumber\\
&\qquad + 2\beta_5 \norm{\bff{u}_n\cdot\Delta\bff{u}_n}{\bb{L}^2}^2
+ \beta_5 \norm{|\bff{u}_n||\Delta\bff{u}_n|}{\bb{L}^2}^2
\nonumber\\
&= 
-\beta_1 \norm{\Delta \bff{u}_n}{\bb{L}^2}^2
+ \beta_3 \norm{\nabla \bff{u}_n}{\bb{L}^2}^2 
- 2\beta_5 \inpro{|\nabla\bff{u}_n|^2\bff{u}_n}{\Delta\bff{u}_n}_{\bb{L}^2}
- 4\beta_5
\inpro{\nabla\bff{u}_n(\bff{u}_n\cdot\nabla\bff{u}_n)^\top}{\Delta\bff{u}_n}_{\bb{L}^2}
\nonumber\\
&\le
|\beta_1| \norm{\Delta \bff{u}_n}{\bb{L}^2}^2
+ \beta_3 \norm{\nabla \bff{u}_n}{\bb{L}^2}^2 
+ 6\beta_5
\int_{\Omega} |\bff{u}_n| |\Delta \bff{u}_n| |\nabla \bff{u}_n|^2 \,\dx.
\end{align}
Using H\"older's inequality and Young's inequality, we can estimate the last term
on the right-hand side by
\begin{align*}
6\beta_5 \int_\Omega |\bff{u}_n| |\Delta \bff{u}_n| |\nabla \bff{u}_n|^2 \,\dx 
&\leq 
6\beta_5 \norm{ |\bff{u}_n| |\Delta \bff{u}_n| }{\bb{L}^2} 
\norm{\nabla \bff{u}_n}{\bb{L}^4}^2
\\
&\leq 
\epsilon \norm{ |\bff{u}_n| |\Delta \bff{u}_n| }{\bb{L}^2}^2 
+ C \norm{\nabla \bff{u}_n}{\bb{L}^4}^4,
\end{align*}
where $\varepsilon>0$ is sufficiently small. This inequality together
with~\eqref{equ:mul Del un} yields
\begin{align} \label{equ:nab un L4} 
	\ddt \norm{\nabla \bff{u}_n}{\bb{L}^2}^2 
	&+ \norm{\nabla \Delta \bff{u}_n}{\bb{L}^2}^2 
	+ \norm{\bff{u}_n \cdot \nabla \bff{u}_n}{\bb{L}^2}^2 
	+ \norm{ |\bff{u}_n| |\nabla \bff{u}_n| }{\bb{L}^2}^2 
	+ \norm{\bff{u}_n\cdot\Delta\bff{u}_n}{\bb{L}^2}^2
	+ \norm{|\bff{u}_n||\Delta\bff{u}_n|}{\bb{L}^2}^2
	\nonumber\\
	&\lesssim 
	\norm{\nabla \bff{u}_n}{\bb{L}^2}^2
	+ \norm{\Delta \bff{u}_n}{\bb{L}^2}^2
	+ \norm{\nabla \bff{u}_n}{\bb{L}^4}^4.
\end{align}
We estimate the last term on the right-hand side of \eqref{equ:nab un L4} by
invoking Gagliardo--Nirenberg's inequality \eqref{gagliardo}. 

\medskip
\noindent
\underline{Case 1: $d=1$}. Applying inequality~\eqref{gagliardo}
with~$\bff{v}=\bff{u}_n$, $q=4$, $r=1$, $s_1=0$, and $s_2=3$ gives
\begin{align*} 
	\norm{\nabla\bff{u}_n}{\bb{L}^4}^4
	&\lesssim
	\norm{\bff{u}_n}{\bb{L}^2}^{7/3}
	\norm{\bff{u}_n}{\bb{H}^3}^{5/3}
	\lesssim
	\norm{\bff{u}_n}{\bb{H}^3}^{5/3}
\end{align*}
where in the last step we used Proposition~\ref{pro:Delta un} and the assumption
that $\bff{u}_n(0)=\bff{u}_{0n}\in\bb{V}_n$ which approximates~$\bff{u}_0$. 
Young's inequality implies, for any $\epsilon>0$,
\[
	\norm{\nabla\bff{u}_n}{\bb{L}^4}^4
	\le C + \epsilon\norm{\bff{u}_n}{\bb{H}^3}^{2}
	\lesssim 
	1 
	+ \norm{\nabla\bff{u}_n}{\bb{L}^2}^2
	+ \epsilon\norm{\nabla\Delta\bff{u}_n}{\bb{L}^2}^2,
\]
where in the last step we used \eqref{eq5}. Therefore, by choosing $\epsilon>0$
sufficiently small, we deduce from \eqref{equ:nab un L4}
\begin{align*}
\ddt\|\nabla \bff{u}_n(t)\|_{\bb{L}^2}^2 
\lesssim
1 + \norm{\nabla\bff{u}_n(t)}{\bb{L}^2}^2
+ \norm{\Delta\bff{u}_n(t)}{\bb{L}^2}^2.
%&+ 
%\int_{0}^{t} \|\nabla \Delta \bff{u}_n(s)\|_{\bb{L}^2}^2 \ds
%+ \int_{0}^{t} \|\bff{u}_n(s) \cdot \nabla \bff{u}_n(s)\|_{\bb{L}^2}^2 \ds
%+ \int_{0}^{t} \| |\bff{u}_n(s)| |\nabla \bff{u}_n(s)| \|_{\bb{L}^2}^2 \ds
%\\
%&\quad
%+ \int_{0}^{t} \norm{\bff{u}_n(s) \cdot\Delta\bff{u}_n(s) }{\bb{L}^2}^2 \ds
%+ \int_{0}^{t} \norm{|\bff{u}_n(s)||\Delta\bff{u}_n(s)|}{\bb{L}^2}^2 \ds
%\nonumber\\
%&\le
%\|\nabla \bff{u}_n(0)\|_{\bb{L}^2}^2 + C
%+ C\int_{0}^{t} \|\Delta \bff{u}_n(s)\|_{\bb{L}^2}^2 \ds
%+ C\int_{0}^{t} \|\nabla \bff{u}_n(s)\|_{\bb{L}^2}^2 \ds
%\\
%&\le
%\|\nabla \bff{u}_n(0)\|_{\bb{L}^2}^2 + C
%+ C\int_{0}^{t} \|\nabla \bff{u}_n(s)\|_{\bb{L}^2}^2 \ds,
\end{align*}
Integrating over $(0,t)$ and using  Proposition~\ref{pro:Delta un}, we obtain
\[
	\norm{\nabla\bff{u}_n(t)}{\bb{L}^2}^2
	\le
	\norm{\nabla\bff{u}_n(0)}{\bb{L}^2}^2
	+
	C + C \int_{0}^{t} \norm{\nabla\bff{u}_n(s)}{\bb{L}^2}^2 \ds.
\]
Gronwall's inequality yields the required estimate.

\medskip
\noindent
\underline{Case 2: $d=2$}. Applying inequality~\eqref{gagliardo}
with~$\bff{v}=\nabla\bff{u}_n$, $q=4$, $r=0$, $s_1=0$, and $s_2=1$ gives
\begin{align*} 
	\norm{\nabla\bff{u}_n}{\bb{L}^4}^4
	&\lesssim
	\norm{\nabla\bff{u}_n}{\bb{L}^2}^{2}
	\norm{\nabla\bff{u}_n}{\bb{H}^1}^{2}
	\lesssim
	\big( 1 + \norm{\Delta\bff{u}_n}{\bb{L}^2}^2 \big)
	\norm{\nabla\bff{u}_n}{\bb{L}^2}^2,
\end{align*}
where in the last step we used \eqref{eq1} and Proposition~\ref{pro:Delta un}.
Therefore, inequality~\eqref{equ:nab un L4} gives
\[
	\ddt\norm{\nabla\bff{u}_n(t)}{\bb{L}^2}^2
	\lesssim
	\norm{\Delta\bff{u}_n(t)}{\bb{L}^2}^2
	+
	\big( 1 + \norm{\Delta\bff{u}_n(t)}{\bb{L}^2}^2 \big)
	\norm{\nabla\bff{u}_n(t)}{\bb{L}^2}^2.
\]
Integrating over $(0,t)$, using Gronwall's inequality and
Proposition~\ref{pro:Delta un}, we deduce
\[
\norm{\nabla\bff{u}_n(t)}{\bb{L}^2}^2
	\lesssim
	1 + \exp\Big( \int_{0}^{T} 
	\big( 1 + \norm{\Delta\bff{u}_n(t)}{\bb{L}^2}^2 \big) \dt \Big)
	\lesssim 1,
\]
proving the result for $d=2$.

\medskip
\noindent
\underline{Case 3: $d=3$}. Applying inequality~\eqref{gagliardo}
with~$\bff{v}=\nabla\bff{u}_n$, $q=4$, $r=0$, $s_1=0$, $s_2=2$, and
using~\eqref{eq5}, we infer
\begin{align*} 
	\norm{\nabla\bff{u}_n}{\bb{L}^4}^4
	&\lesssim
	\norm{\nabla\bff{u}_n}{\bb{L}^2}^{5/2}
	\norm{\nabla\bff{u}_n}{\bb{H}^2}^{3/2}
	\lesssim
	\norm{\nabla\bff{u}_n}{\bb{L}^2}^{5/2}
	\big( 1 + \norm{\nabla\bff{u}_n}{\bb{L}^2}^{2}
	+ \norm{\nabla\Delta\bff{u}_n}{\bb{L}^2}^{2} \big)^{3/4}
	\\
	&\lesssim
	\norm{\nabla\bff{u}_n}{\bb{L}^2}^{5/2}
	+
	\norm{\nabla\bff{u}_n}{\bb{L}^2}^{4}
	+
	\norm{\nabla\bff{u}_n}{\bb{L}^2}^{5/2}
	\norm{\nabla\Delta\bff{u}_n}{\bb{L}^2}^{3/2}.
\end{align*}
Young's inequality yields, for $\epsilon>0$ sufficiently small,
\begin{align*} 
	\norm{\nabla\bff{u}_n}{\bb{L}^4}^4
	\le
	C \Big(
	\norm{\nabla\bff{u}_n}{\bb{L}^2}^{5/2}
	+
	\norm{\nabla\bff{u}_n}{\bb{L}^2}^{4}
	+
	\norm{\nabla\bff{u}_n}{\bb{L}^2}^{10}
	\Big)
	+
	\epsilon \norm{\nabla\Delta\bff{u}_n}{\bb{L}^2}^{2}.
\end{align*}
Inserting this estimate into \eqref{equ:nab un L4} and rearranging the terms, we
deduce
\begin{align*} 
	\ddt \norm{\nabla\bff{u}_n}{\bb{L}^2}^{2}
	&\lesssim
	\norm{\nabla\bff{u}_n}{\bb{L}^2}^{2}
	+
	\norm{\Delta\bff{u}_n}{\bb{L}^2}^{2}
	+
	\norm{\nabla\bff{u}_n}{\bb{L}^2}^{5/2}
	+
	\norm{\nabla\bff{u}_n}{\bb{L}^2}^{4}
	+
	\norm{\nabla\bff{u}_n}{\bb{L}^2}^{10}
	\\
    &\lesssim
	\norm{\nabla\bff{u}_n}{\bb{L}^2}^{2}
	+
	\norm{\Delta\bff{u}_n}{\bb{L}^2}^{2}
	+
	\norm{\nabla\bff{u}_n}{\bb{L}^2}^{10}.
\end{align*}
Integrating over $(0,t)$ and using Proposition~\ref{pro:Delta un} give
\begin{align*} 
	\norm{\nabla\bff{u}_n(t)}{\bb{L}^2}^{2}
	&\le
	\norm{\nabla\bff{u}_n(0)}{\bb{L}^2}^{2} + C
	+
	C \int_{0}^{t}
	\norm{\nabla\bff{u}_n(s)}{\bb{L}^2}^{10} \ds.
\end{align*}
{By using Gronwall--Bihari's inequality (Theorem~\ref{the:bihari}) with $f(x)=x^5$ (so
that $F(x)=x^{-4}/4$ and $F^{-1}(x)=x^{-1/4}/\sqrt{2}$), and
noting~\eqref{equ:un0}, we obtain the required estimate for any $t\in
[0,T^\ast]$, where $T^\ast = \big( \norm{\nabla\bff{u}_n(0)}{\bb{L}^2}^{2} + C
\big)^{-4}/4$.  This completes the proof of the proposition.
}
\end{proof}

{
\begin{remark}
Since~$\bff{u}_n(0) = \bff{u}_{0n}$ approximates $\bff{u}_0$, we have
$T^{\ast} \approx \norm{\nabla\bff{u}_0}{\bb{L}^2}^{-8}$.
\end{remark}
}

\begin{proposition}\label{pro:dt un H-1}
	{Under the assumption of Proposition \ref{pro:nab delta un}}, we have
	\begin{align*}
		\norm{\partial_t \bff{u}_n}{L^2(0,T^\ast; \bb{H}^{-1})}^2
		\lesssim
		\norm{\bff{u}_0}{\bb{H}^1}^2,
	\end{align*}
where the constant depends on $T^\ast$ but is independent of $n$.
\end{proposition}

\begin{proof}
	Taking the inner product of \eqref{approx} with $\bff{\varphi}\in \bb{H}^1$ such that $\norm{\bff{\varphi}}{\bb{H}^1} \leq 1$ and integrating by parts with respect to $\bff{x}$, we have
	\begin{align*}
		\inpro{\partial_t \bff{u}_n}{\bff{\varphi}}_{\bb{L}^2}
		&=
		-\beta_1 \inpro{\nabla \bff{u}_n}{\nabla \bff{\varphi}}_{\bb{L}^2}
		+
		\beta_2 \inpro{\nabla \Delta \bff{u}_n}{\nabla \bff{\varphi}}_{\bb{L}^2}
		+
		\beta_3 \inpro{\bff{u}_n}{\bff{\varphi}}_{\bb{L}^2}
		-
		\beta_3 \inpro{|\bff{u}_n|^2 \bff{u}_n}{\bff{\varphi}}_{\bb{L}^2}
		\\
		&\quad
		+
		\beta_4 \inpro{\bff{u}_n \times \nabla \bff{u}_n}{\nabla \bff{\varphi}}_{\bb{L}^2}
		+
		\beta_5 \inpro{\nabla(|\bff{u}_n|^2 \bff{u}_n)}{\nabla \bff{\varphi}}_{\bb{L}^2}.
	\end{align*}
	It follows from this equation and H\"{o}lder's inequality that
	\begin{align*}
		\big| \inpro{\partial_t \bff{u}_n}{\bff{\varphi}}_{\bb{L}^2} \big|
		&\lesssim
		\norm{\bff{u}_n}{\bb{H}^1} \norm{\bff{\varphi}}{\bb{H}^1}
		+
		\norm{\bff{u}_n}{\bb{H}^3} \norm{\bff{\varphi}}{\bb{H}^1}
		+
		\norm{\bff{u}_n}{\bb{L}^2} \norm{\bff{\varphi}}{\bb{L}^2}
		+
		\norm{\bff{u}_n}{\bb{L}^2} \norm{\bff{u}_n}{\bb{L}^6}^2 \norm{\bff{\varphi}}{\bb{L}^6}
		\\
		&\quad
		+
		\norm{\bff{u}_n}{\bb{L}^4} \norm{\nabla \bff{u}_n}{\bb{L}^4}
		\norm{\nabla \bff{\varphi}}{\bb{L}^2}
		+
		\norm{\bff{u}_n}{\bb{L}^6} \norm{\bff{u}_n}{\bb{L}^6}
		\norm{\nabla \bff{u}_n}{\bb{L}^6} \norm{\nabla \bff{\varphi}}{\bb{L}^2}.
	\end{align*}
	Therefore, integrating over $(0,t)$, using the Sobolev embedding $\bb{H}^1 \subset \bb{L}^6$, and noting that $\norm{\bff{\varphi}}{\bb{H}^1} \leq 1$, we obtain
	\[
	\int_0^t \big| \inpro{\partial_s \bff{u}_n(s)}{\bff{\varphi}}_{\bb{L}^2} \big|^2 \,\ds
	\lesssim
	\int_0^t \norm{\bff{u}_n}{\bb{H}^3}^2 \ds \lesssim 1
	\]
	by Proposition \ref{pro:Delta un}. Taking supremum over the set
	$\{\varphi\in \bb{H}^1: \norm{\bff{\varphi}}{\bb{H}^1} \leq 1\}$ and
	noting that
	\[
	\norm{\partial_t\bff{u}_n}{\bb{H}^{-1}}
	\le
	\sup_{\bff{\varphi}\in\bb{H}^1}
	\frac{\big|{\inpro{\partial_t\bff{u}_n}{\bff{\varphi}}_{\bb{L}^2}}\big|}
	{\norm{\bff{\varphi}}{\bb{H}^1}}, 
	\]
	we obtain the required estimate.
\end{proof}

\begin{proposition}\label{pro:partial un}
	{Assume that $\bff{u}_0\in\bb{H}^2$.}
	Let $T>0$ be arbitrary and $T^\ast$ be defined as in
	Proposition~\ref{pro:nab delta un}. For each $n\in \bb{N}$ and all $t\in
	[0,T^\ast]$,
\begin{align*} 
	\int_{0}^{t} \norm{\partial_s\bff{u}_n(s)}{\bb{L}^2}^2 \ds
	&+
	\norm{\Delta\bff{u}_n(t)}{\bb{L}^2}^2
	+
	\norm{\bff{u}_n(t)}{\bb{L}^4}^4
	+
	\norm{\bff{u}_n(t)\cdot\nabla\bff{u}_n(t)}{\bb{L}^2}^2
	+
	\norm{|\bff{u}_n(t)| \ |\nabla\bff{u}_n(t)|}{\bb{L}^2}^2
	\\
	 &\lesssim
	\norm{\bff{u}_n(0)}{\bb{H}^2}^2
	{
	\lesssim
	\norm{\bff{u}_0}{\bb{H}^2}^2,
}
\end{align*}
where the constant depends on $T^\ast$, but is independent of $n$. Here,
$\partial_s\bff{u}_n:=\displaystyle\frac{\partial \bff{u}_n}{\partial s}$.
\end{proposition}

\begin{proof}
	Taking the inner product of \eqref{approx} with 
	$\partial_t\bff{u}_n$ and integrating by parts with respect
	to~$\bff{x}$, we obtain
\begin{align*}
	\norm{\partial_t \bff{u}_n}{\bb{L}^2}^2 
	&+ \frac{\beta_1}{2} \ddt \norm{\nabla \bff{u}_n}{\bb{L}^2}^2 
	+ \frac{\beta_2}{2} \ddt \norm{\Delta \bff{u}_n}{\bb{L}^2}^2 
	+ \frac{\beta_3}{4} \ddt \norm{\bff{u}_n}{\bb{L}^4}^4 
	+ \beta_4 \inpro{\bff{u}_n \times
	\Delta\bff{u}_n}{\partial_t\bff{u}_n}_{\bb{L}^2}
	\\
	&+ \beta_5 \inpro{\nabla(|\bff{u}_n|^2 \bff{u}_n)}
	{\partial_t\nabla \bff{u}_n}_{\bb{L}^2} 
	= 
	\frac{\beta_3}{2} \ddt \norm{\bff{u}_n}{\bb{L}^2}^2.
\end{align*}
For the last term on the left-hand side, it follows from \eqref{equ:nab un2}
that
\begin{align*} 
	\beta_5 
	\inpro{\nabla(|\bff{u}_n|^2 \bff{u}_n)}{\partial_t\nabla\bff{u}_n}_{\bb{L}^2} 
	&=
	2\beta_5
	\inpro{\bff{u}_n\big(\bff{u}_n\cdot\nabla\bff{u}_n\big)}%
	{\partial_t\nabla\bff{u}_n}_{\bb{L}^2}
	+
	\beta_5
	\inpro{|\bff{u}_n|^2\nabla\bff{u}_n}{\partial_t\nabla\bff{u}_n}_{\bb{L}^2}
	\\
	&=
	\beta_5
	\ddt \norm{\bff{u}_n\cdot\nabla\bff{u}_n}{\bb{L}^2}^2
	-
	\beta_5
	\inpro{\bff{u}_n\cdot\nabla\bff{u}_n}%
	{\partial_t\bff{u}_n\cdot\nabla\bff{u}_n}_{\bb{L}^2}
	\\
	&\quad
	+
	\frac{\beta_5}{2} 
	\inpro{|\bff{u}_n|^2}{\partial_t\big(|\nabla\bff{u}_n|^2\big)}_{L^2}
	\\
	&=
	\beta_5
	\ddt \norm{\bff{u}_n\cdot\nabla\bff{u}_n}{\bb{L}^2}^2
	-
	\beta_5
	\inpro{\bff{u}_n\cdot\nabla\bff{u}_n}%
	{\partial_t\bff{u}_n\cdot\nabla\bff{u}_n}_{\bb{L}^2}
	\\
	&\quad
	+
	\frac{\beta_5}{2} 
	\ddt \norm{|\bff{u}_n| \ |\nabla\bff{u}_n|}{L^2}^2
	-
	\beta_5
	\inpro{|\nabla\bff{u}_n|^2}{\bff{u}_n\cdot\partial_t\bff{u}_n}_{L^2}.
\end{align*}
Therefore,
\begin{align*} 
	\norm{\partial_t \bff{u}_n}{\bb{L}^2}^2 
	&+ \frac{\beta_1}{2} \ddt \norm{\nabla \bff{u}_n}{\bb{L}^2}^2 
	+ \frac{\beta_2}{2} \ddt \norm{\Delta \bff{u}_n}{\bb{L}^2}^2 
	+ \frac{\beta_3}{4} \ddt \norm{\bff{u}_n}{\bb{L}^4}^4 
	\\
	&\quad
	+ \beta_5
	\ddt \norm{\bff{u}_n\cdot\nabla\bff{u}_n}{\bb{L}^2}^2
	+ \frac{\beta_5}{2} 
	\ddt \norm{|\bff{u}_n| \ |\nabla\bff{u}_n|}{L^2}^2
	\\
	&=
	\frac{\beta_3}{2} \ddt \norm{\bff{u}_n}{\bb{L}^2}^2
	- \beta_4 
	\inpro{\bff{u}_n\times\Delta\bff{u}_n}{\partial_t\bff{u}_n}_{\bb{L}^2}
	\\
	&\quad
	+ \beta_5
	\inpro{\bff{u}_n\cdot\nabla\bff{u}_n}%
	{\partial_t\bff{u}_n\cdot\nabla\bff{u}_n}_{\bb{L}^2}
	+ \beta_5
	\inpro{|\nabla\bff{u}_n|^2}{\bff{u}_n\cdot\partial_t\bff{u}_n}_{L^2}
	\\
	&\le
	\frac{\beta_3}{2} \ddt \norm{\bff{u}_n}{\bb{L}^2}^2
	+ \beta_4 
	\norm{\bff{u}_n}{\bb{L}^\infty} \norm{\Delta\bff{u}_n}{\bb{L}^2}
	\norm{\partial_t \bff{u}_n}{\bb{L}^2}
	+ 2\beta_5
	\norm{\bff{u}_n}{\bb{L}^\infty} \norm{\nabla\bff{u}_n}{\bb{L}^4}^2
	\norm{\partial_t\bff{u}_n}{\bb{L}^2}
	\\
	&\le
	\frac{\beta_3}{2} \ddt \norm{\bff{u}_n}{\bb{L}^2}^2
	+
	C \norm{\bff{u}_n}{\bb{L}^\infty}^2 \norm{\Delta\bff{u}_n}{\bb{L}^2}^2
	+
	C \norm{\bff{u}_n}{\bb{L}^\infty}^2 \norm{\nabla\bff{u}_n}{\bb{L}^4}^4
	+
	\epsilon \norm{\partial_t\bff{u}_n}{\bb{L}^2}^2,
\end{align*}
for any $\epsilon>0$, where in the last step we used Young's inequality.
Rearranging the inequality, we obtain
\begin{align}\label{equ:dtu} 
	\norm{\partial_t \bff{u}_n}{\bb{L}^2}^2 
	&+ \ddt \norm{\nabla \bff{u}_n}{\bb{L}^2}^2 
	+ \ddt \norm{\Delta \bff{u}_n}{\bb{L}^2}^2 
	+ \ddt \norm{\bff{u}_n}{\bb{L}^4}^4 
	+ \ddt \norm{\bff{u}_n\cdot\nabla\bff{u}_n}{\bb{L}^2}^2
	+ \ddt \norm{|\bff{u}_n| \ |\nabla\bff{u}_n|}{L^2}^2
	\nonumber\\
	&\lesssim
	\ddt \norm{\bff{u}_n}{\bb{L}^2}^2
	+
	\norm{\bff{u}_n}{\bb{L}^\infty}^2 \norm{\Delta\bff{u}_n}{\bb{L}^2}^2
	+
	\norm{\bff{u}_n}{\bb{L}^\infty}^2 \norm{\nabla\bff{u}_n}{\bb{L}^4}^4.
\end{align}
We now estimate the last two terms on the right-hand side of \eqref{equ:dtu}.

\medskip
\noindent
\underline{Case 1: $d=1$}. It follows from the Sobolev embedding,
Proposition~\ref{pro:Delta un}, and Proposition~\ref{pro:nab delta un} that
\[
	\norm{\bff{u}_n(t)}{\bb{L}^\infty}^2 
	\lesssim
	\norm{\bff{u}_n(t)}{\bb{H}^1}^2 
	\lesssim
	1, \quad t\in[0,T].
\]
Moreover, the Gagliardo--Nirenberg inequality (Theorem~\ref{the:Gal Nir} with
$\bff{v}=\bff{u}_n$, $q=4$, $r=1$, $s_1=1$, and $s_2=2$) together
with~\eqref{eq1} implies
\begin{align*} 
	\norm{\nabla\bff{u}_n(t)}{\bb{L}^4}^4
	&\lesssim
	\norm{\bff{u}_n(t)}{\bb{H}^1}^{3}
	\norm{\bff{u}_n(t)}{\bb{H}^2}
	\lesssim
	1 + \norm{\Delta\bff{u}_n(t)}{\bb{L}^2}^2,
	\quad t\in[0,T].
\end{align*}
Therefore, inequality~\eqref{equ:dtu} yields the required result, after
integrating over~$(0,t)$ and using Proposition~\ref{pro:Delta un}.

\medskip
\noindent
\underline{Case 2: $d=2$}. The Gagliardo--Nirenberg inequality
(respectively with $\bff{v}=\bff{u}_n$, $q=\infty$, $r=s_1=0$, $s_2=2$, and
with~$\bff{v}=\nabla\bff{u}_n$, $q=4$, $r=s_1=0$, $s_2=1$) implies
\begin{alignat*}{2} 
	\norm{\bff{u}_n(t)}{\bb{L}^\infty}^2
	&\lesssim
	\norm{\bff{u}_n(t)}{\bb{L}^2}
	\norm{\bff{u}_n(t)}{\bb{H}^2}
	\lesssim
	\norm{\bff{u}_n(t)}{\bb{H}^2}
	\lesssim
	1 + \norm{\Delta\bff{u}_n(t)}{\bb{L}^2},
	\quad && t\in[0,T],
	\\
	\norm{\nabla\bff{u}_n(t)}{\bb{L}^4}^4
	&\lesssim
	\norm{\nabla\bff{u}_n(t)}{\bb{L}^2}^2
	\norm{\nabla\bff{u}_n(t)}{\bb{H}^1}^2
	\lesssim
	\norm{\bff{u}_n(t)}{\bb{H}^2}^2
	\lesssim
	1 + \norm{\Delta\bff{u}_n(t)}{\bb{L}^2}^2,
	\quad && t\in[0,T],
\end{alignat*}
where we also used~\eqref{eq1} and Proposition~\ref{pro:nab delta un}. Inserting
these estimates into~\eqref{equ:dtu} and integrating over~$(0,t)$ yield
\begin{align}\label{equ:dsu} 
	\int_{0}^{t} \norm{\partial_s\bff{u}_n(s)}{\bb{L}^2}^2 \ds
	&+
	\norm{\Delta\bff{u}_n(t)}{\bb{L}^2}^2
	+
	\norm{\bff{u}_n(t)}{\bb{L}^4}^4
	+ 
	\norm{\bff{u}_n(t)\cdot\nabla\bff{u}_n(t)}{\bb{L}^2}^2
	+ 
	\norm{|\bff{u}_n(t)| \ |\nabla\bff{u}_n(t)|}{L^2}^2
	\nonumber\\
	&\lesssim
	1 + \int_{0}^{t} \Big(
	\norm{\Delta\bff{u}_n(s)}{\bb{L}^2}
	+
	\norm{\Delta\bff{u}_n(s)}{\bb{L}^2}^2
	+
	\norm{\Delta\bff{u}_n(s)}{\bb{L}^2}^3
	\Big) \ds
	\nonumber\\
	&\lesssim
	1 + \int_{0}^{t} \Big(
	\norm{\Delta\bff{u}_n(s)}{\bb{L}^2}^2
	+
	\norm{\Delta\bff{u}_n(s)}{\bb{L}^2}^3
	\Big) \ds,
\end{align}
where in the last step we used Young's inequality for the
term~$\norm{\Delta\bff{u}_n(s)}{\bb{L}^2}$.
For the last term on the right-hand side, we use~\eqref{eq3} to obtain
\begin{align*} 
	\norm{\Delta\bff{u}_n}{\bb{L}^2}^3
	&\le
	\norm{\nabla\bff{u}_n}{\bb{L}^2}^{3/2}
	\norm{\nabla\Delta\bff{u}_n}{\bb{L}^2}^{3/2}
	\lesssim
	\norm{\nabla\Delta\bff{u}_n}{\bb{L}^2}^{3/2}
	\lesssim
	1 + \norm{\nabla\Delta\bff{u}_n}{\bb{L}^2}^{2},
\end{align*}
where in the penultimate step we used Proposition~\ref{pro:nab delta un}, and in the
last step we used Young's inequality. Therefore the right-hand side
of~\eqref{equ:dsu} is bounded independent of $n$ due to
Proposition~\ref{pro:nab delta un}, proving the proposition for this case.

\medskip
\noindent
\underline{Case 3: $d=3$}. The Gagliardo--Nirenberg inequality
(respectively with $\bff{v}=\bff{u}_n$, $q=\infty$, $r=0$, $s_1=1$, $s_2=2$, and
with $\bff{v}=\nabla\bff{u}_n$, $q=4$, $r=s_1=0$, $s_2=1$) implies
\begin{align*}
	\norm{\bff{u}_n(t)}{\bb{L}^\infty}^2
	&\lesssim
	\norm{\bff{u}_n(t)}{\bb{H}^1}
	\norm{\bff{u}_n(t)}{\bb{H}^2}
	\lesssim
	\norm{\bff{u}_n(t)}{\bb{H}^2}
	\lesssim
	1 + \norm{\Delta\bff{u}_n(t)}{\bb{L}^2},
	\\
	\norm{\nabla\bff{u}_n(t)}{\bb{L}^4}^4
	&\lesssim
	\norm{\nabla\bff{u}_n(t)}{\bb{L}^2}
	\norm{\nabla\bff{u}_n(t)}{\bb{H}^1}^3
	\lesssim
	\norm{\bff{u}_n(t)}{\bb{H}^2}^3
	\lesssim
	1 + \norm{\Delta\bff{u}_n(t)}{\bb{L}^2}^3,
\end{align*}
for all $t\in[0,T^\ast]$ where $T^\ast$ is given in Proposition~\ref{pro:nab
delta un}. {Inserting these estimates into~\eqref{equ:dtu}, integrating
over~$(0,t)$, and using~\eqref{equ:un0} yield}
\begin{align*}
	\int_{0}^{t} \norm{\partial_s\bff{u}_n(s)}{\bb{L}^2}^2 \ds
	&+
	\norm{\Delta\bff{u}_n(t)}{\bb{L}^2}^2
	+
	\norm{\bff{u}_n(t)}{\bb{L}^4}^4
	+ 
	\norm{\bff{u}_n(t)\cdot\nabla\bff{u}_n(t)}{\bb{L}^2}^2
	+ 
	\norm{|\bff{u}_n(t)| \ |\nabla\bff{u}_n(t)|}{L^2}^2
	\nonumber\\
	&\lesssim
	1 + \int_{0}^{t} \Big(
	\norm{\Delta\bff{u}_n(s)}{\bb{L}^2}
	+
	\norm{\Delta\bff{u}_n(s)}{\bb{L}^2}^2
	+
	\norm{\Delta\bff{u}_n(s)}{\bb{L}^2}^3
	+
	\norm{\Delta\bff{u}_n(s)}{\bb{L}^2}^4
	\Big) \ds
	\nonumber\\
	&\lesssim
	1 + \int_{0}^{t} \Big(
	\norm{\Delta\bff{u}_n(s)}{\bb{L}^2}^2
	+
	\norm{\Delta\bff{u}_n(s)}{\bb{L}^2}^4
	\Big) \ds
	\\
	&\lesssim
	1 + \int_{0}^{t} \Big(
	\norm{\Delta\bff{u}_n(s)}{\bb{L}^2}^2
	+
	\norm{\nabla\Delta\bff{u}_n(s)}{\bb{L}^2}^2
	\Big) \ds
	\\
	&\lesssim 1,
\end{align*}
where in the last step we used Proposition~\ref{pro:nab delta un}, completing the proof
of the proposition.
\end{proof}

%\begin{remark}\label{estimate4}
%Proposition \ref{pro:partial un} implies $\partial_t\bff{u}_n \in L^2(0,T;\bb{L}^2)$. %Therefore,
%by \eqref{approx}, we have $\Delta^2 \bff{u}_n \in L^2(0,T;\bb{L}^2)$.
%\end{remark}

\begin{proposition}\label{pro:Del sqr}
	{Assume that $\bff{u}_0\in\bb{H}^2$.}
	Let $T>0$ be arbitrary and $T^\ast$ be defined as in
	Proposition~\ref{pro:nab delta un}. For each $n\in \bb{N}$ and all $t\in
	[0,T^\ast]$,
\begin{align*} 
	\int_{0}^{t} 
	\norm{\Delta^2\bff{u}_n(s)}{\bb{L}^2}^2 \ds
	\lesssim
	\norm{\bff{u}_n(0)}{\bb{H}^2}^2
	{
	\lesssim
	\norm{\bff{u}_0}{\bb{H}^2}^2,
}
\end{align*}
where the constant depends on $T^\ast$ but is independent of $n$.
\end{proposition}

\begin{proof}
	Taking the inner product of \eqref{approx} with 
	$\Delta^2\bff{u}_n$ and integrating by parts with respect
	to~$\bff{x}$, we obtain
\begin{align}
	\label{equ:Del sqr0}
	&\frac{1}{2} \ddt \norm{\Delta \bff{u}_n}{\bb{L}^2}^2 
	+ 
	\beta_1 \norm{\nabla \Delta \bff{u}_n}{\bb{L}^2}^2
	+
	\beta_2 \norm{ \Delta^2 \bff{u}_n}{\bb{L}^2}^2
	\nonumber
	\\
	&= 
	\beta_3 \norm{\Delta \bff{u}_n}{\bb{L}^2}^2 
	-
	\beta_3 \inpro{ \abs{\bff{u}_n}^2 \bff{u}_n}{\Delta^2 \bff{u}_n }_{\bb{L}^2} 
	\nonumber
	\\
	&\quad-
	\beta_4 \inpro{ \bff{u}_n \times \Delta \bff{u}_n}{\Delta^2 \bff{u}_n }_{\bb{L}^2} 
	- 
	\beta_5  \inpro{ \Delta (\abs{\bff{u}_n}^2 \bff{u}_n)}{\Delta^2
	\bff{u}_n}_{\bb{L}^2}.   
\end{align}
Each term on the right-hand side can be estimated as follows. For the first term, by Young's inequality, Sobolev embedding and Proposition \ref{pro:nab delta un},
\begin{align*}
	\nonumber
	\big|\inpro{ |\bff{u}_n|^2 \bff{u}_n}{\Delta^2 \bff{u}_n }_{\bb{L}^2} \big|
	&\leq
	C \norm{|\bff{u}_n|^2 \bff{u}_n}{\bb{L}^2}^2 
	+ 
	\epsilon \norm{\Delta^2 \bff{u}_n}{\bb{L}^2}^2 
	\\
	\nonumber
	&=
	C \norm{\bff{u}_n}{\bb{L}^6}^6
	+
	\epsilon \norm{\Delta^2 \bff{u}_n}{\bb{L}^2}^2
	\\
	\nonumber
	&\leq
	C \norm{\bff{u}_n}{\bb{H}^1}^6
	+
	\epsilon \norm{\Delta^2 \bff{u}_n}{\bb{L}^2}^2
	\\
	&\lesssim
	1 + \epsilon \norm{\Delta^2 \bff{u}_n}{\bb{L}^2}^2  %\label{equ:Del sqr1}
\end{align*}
for any $\epsilon > 0$.
For the second term, by H\"{o}lder's inequality, Young's inequality, Sobolev
embedding, and Proposition \ref{pro:partial un}, we have
\begin{align*}
	\nonumber
	\big| 
	\inpro{ \bff{u}_n \times \Delta \bff{u}_n}{\Delta^2 \bff{u}_n }_{\bb{L}^2}
	\big|
	&\leq
	\norm{\bff{u}_n}{\bb{L}^\infty}
	\norm{\Delta \bff{u}_n}{\bb{L}^2}
	\norm{\Delta^2 \bff{u}_n}{\bb{L}^2}	
	\\
	\nonumber
	&\leq	
	C(\norm{\bff{u}_n}{\bb{H}^2}^2 \norm{\Delta \bff{u}_n}{\bb{L}^2}^2)
	+
	\epsilon \norm{\Delta^2 \bff{u}_n}{\bb{L}^2}^2
	\\
	&\lesssim
	1 + \epsilon \norm{\Delta^2 \bff{u}_n}{\bb{L}^2}^2.  %\label{equ:Del sqr2}
\end{align*}
Finally, by H\"{o}lder's and Young's inequalities, we have
\begin{align}\label{equ:Del sqr3}
	\big|
	\inpro{ \Delta (\abs{\bff{u}_n}^2 \bff{u}_n)}{\Delta^2 \bff{u}_n}
	\big|
	&\leq
	\norm{\Delta(\abs{\bff{u}_n}^2 \bff{u}_n)}{\bb{L}^2}
	\norm{\Delta^2 \bff{u}_n}{\bb{L}^2}
	\nonumber
	\\
	&\le
	C \norm{\Delta(\abs{\bff{u}_n}^2 \bff{u}_n)}{\bb{L}^2}^2
	+
	\epsilon \norm{\Delta^2 \bff{u}_n}{\bb{L}^2}^2.
\end{align}
For the first term on the right-hand side, it follows from~\eqref{equ:del v2v},
H\"older's inequality, and Sobolev embedding that
\begin{align}\label{equ:Del sqr3 b} 
	\norm{\Delta(\abs{\bff{u}_n}^2 \bff{u}_n)}{\bb{L}^2}^2
	&\lesssim
	\norm{\nabla \bff{u}_n}{\bb{L}^6}^4 \norm{\bff{u}_n}{\bb{L}^6}^2 
	+
	\norm{\Delta \bff{u}_n}{\bb{L}^6}^2 \norm{\bff{u}_n}{\bb{L}^6}^4
	\nonumber
	\\
	&\leq
	\norm{\nabla \bff{u}_n}{\bb{H}^1}^4 \norm{\bff{u}_n}{\bb{H}^1}^2
	+ 
	\norm{\Delta \bff{u}_n}{\bb{H}^1}^2 \norm{\bff{u}_n}{\bb{H}^1}^4
	\nonumber
	\\
	&\lesssim
	\norm{\bff{u}_n}{\bb{H}^2}^6
	+
	\norm{\Delta \bff{u}_n}{\bb{H}^1}^2 \norm{\bff{u}_n}{\bb{H}^1}^4
	\nonumber
	\\
	&\lesssim
	1 + \norm{\nabla \Delta \bff{u}_n}{\bb{L}^2}^2,
\end{align}
where in the last step we also used Proposition~\ref{pro:partial un}.
Altogether, we deduce from~\eqref{equ:Del sqr0} after integrating over~$(0,t)$ that
\begin{align*}
	\norm{\Delta \bff{u}_n(t)}{\bb{L}^2}^2 
	+
	\int_0^t \norm{ \Delta^2 \bff{u}_n(s)}{\bb{L}^2}^2 \ds
 	%+
 	%\int_{0}^{t} 
 	%\norm{\bff{u}_n(s)\cdot\nabla\Delta\bff{u}_n(s)}{\bb{L}^2}^2 \ds
 	%\\
 	%&\quad
 	%+
 	%\int_{0}^{t} 
 	%\norm{|\bff{u}_n(s)| \ |\nabla\Delta\bff{u}_n(s)|}{\bb{L}^2}^2 \ds
 	%\\
	\lesssim
	1
	+
	\int_0^t \norm{\nabla \Delta \bff{u}_n(s)}{\bb{L}^2}^2 \ds
	\lesssim
	1,
\end{align*}
where in the last step we used Proposition \ref{pro:nab delta un}. This
completes the proof of the proposition.
\end{proof}

\begin{proposition}\label{pro:nab Del sq}
{Assume that $\bff{u}_0\in\bb{H}^3$.}
Let $T>0$ be arbitrary and $T^\ast$ be defined as in
Proposition~\ref{pro:nab delta un}. For each $n\in \bb{N}$ and all $t\in
[0,T^\ast]$,
\begin{align*}
	\|\nabla\Delta \bff{u}_n(t)\|_{\bb{L}^2}^2 
	+  
	\int_0^t \|\nabla\Delta^2 \bff{u}_n(s)\|_{\bb{L}^2}^2\,\ds 
	\lesssim 
	\|\bff{u}_n(0)\|_{\bb{H}^3}^2
	{\lesssim
	\|\bff{u}_0\|_{\bb{H}^3}^2,}
\end{align*}
where the constant depends on $T^*$ but is independent of $n$.
\end{proposition}

\begin{proof}
Taking the inner product of \eqref{approx} with $\Delta^3 \bff{u}_n$ and
integrating by parts with respect to $\bff{x}$, we have
\begin{align}\label{equ:nab Del sq-1}
	&\frac{1}{2} \ddt \norm{\nabla \Delta \bff{u}_n}{\bb{L}^2}^2 
	+ 
	\beta_1 \norm{\Delta^2 \bff{u}_n}{\bb{L}^2}^2 
	+ 
	\beta_2 \norm{\nabla \Delta^2 \bff{u}_n}{\bb{L}^2}^2 
	-
	\beta_3 \norm{\nabla \Delta \bff{u}_n}{\bb{L}^2}^2 
	\nonumber\\
	&=
	\beta_3 \inpro{ \Delta \big(|\bff{u}_n|^2 \bff{u}_n\big)}{ \Delta^2 \bff{u}_n }_{\bb{L}^2} 
	- 
	2 \beta_4 \inpro{ \nabla \bff{u}_n \times \nabla \Delta \bff{u}_n}{ \Delta^2 \bff{u}_n}_{\bb{L}^2} 
	\nonumber\\
	&\quad
	- 
	\beta_5 \inpro {\nabla\Delta (|\bff{u}_n|^2 \bff{u}_n)}{ \nabla\Delta^2 \bff{u}_n }_{\bb{L}^2}.
\end{align}
Each term on the right-hand side can be estimated as follows. For the first
term, by \eqref{equ:Del sqr3} and~\eqref{equ:Del sqr3 b} we have
\begin{align} \label{equ:nab Del sq-2}
	\big|
	\inpro{ \Delta \big(|\bff{u}_n|^2 \bff{u}_n\big)}{ \Delta^2 \bff{u}_n }_{\bb{L}^2} 
	\big|
	\lesssim
	1 + \norm{\nabla \Delta \bff{u}_n}{\bb{L}^2}^2
	+
	\epsilon
	\norm{\Delta^2 \bff{u}_n}{\bb{L}^2}^2.
\end{align}
For the second term, Sobolev embedding and H\"{o}lder's inequality give
\begin{align} \label{equ:nab Del sq-3}
	\nonumber
	\big|
	\inpro{ \nabla \bff{u}_n \times \nabla \Delta \bff{u}_n}{ \Delta^2 \bff{u}_n}_{\bb{L}^2}
	\big|
	&\leq
	\norm{\nabla \bff{u}_n}{\bb{L}^3}
	\norm{\nabla \Delta \bff{u}_n}{\bb{L}^6}
	\norm{\Delta^2 \bff{u}_n}{\bb{L}^2}
	\\
	\nonumber
	&\leq
	\norm{\bff{u}_n}{\bb{H}^2}
	\norm{\bff{u}_n}{\bb{H}^4}
	\norm{\Delta^2 \bff{u}_n}{\bb{L}^2}
	\\
	&\lesssim
	1 + \norm{\Delta^2 \bff{u}_n}{\bb{L}^2}^2,
\end{align}
where in the last step we used \eqref{eq1}, \eqref{eq6}, and Proposition \ref{pro:partial un}.
For the last term on the right-hand side of~\eqref{equ:nab Del sq-1}, by
H\"{o}lder's and Young's inequalities, and \eqref{equ:nab un2}, we deduce
\begin{align} \label{equ:nab Del sq-4}
	\nonumber
	\big|
	\inpro {\nabla\Delta (|\bff{u}_n|^2 \bff{u}_n)}{ \nabla\Delta^2 \bff{u}_n }_{\bb{L}^2}
	\big|
	&\leq
	\norm{\nabla(|\bff{u}_n|^2 \bff{u}_n)}{\bb{H}^2} \norm{\nabla \Delta^2 \bff{u}_n}{\bb{L}^2}
	\\
	\nonumber
	&\leq 
	C \norm{\bff{u}_n \big( \bff{u}_n \cdot \nabla \bff{u}_n \big)}{\bb{H}^2}^2
	+
	C \norm{|\bff{u}_n|^2 |\nabla \bff{u}_n|}{H^2}^2
	+
	\epsilon \norm{\nabla \Delta^2 \bff{u}_n}{\bb{L}^2}^2
	\\
	\nonumber
	&\leq
	C \norm{\bff{u}_n}{\bb{H}^2}^4 \norm{\nabla \bff{u}_n}{\bb{H}^2}^2
	+
	\epsilon \norm{\nabla \Delta^2 \bff{u}_n}{\bb{L}^2}^2
	\\
	&\lesssim
	1 + \norm{\nabla \Delta \bff{u}_n}{\bb{L}^2}^2
	+ 
	\epsilon \norm{\nabla \Delta^2 \bff{u}_n}{\bb{L}^2}^2
\end{align}
for any $\epsilon >0$, where in the penultimate step we used~\eqref{eq7} and in
the last step we used \eqref{eq5}, Proposition \ref{pro:nab delta un}, and
Proposition \ref{pro:partial un}.
Inserting the estimates \eqref{equ:nab Del sq-2}, \eqref{equ:nab Del sq-3}, and
\eqref{equ:nab Del sq-4} into \eqref{equ:nab Del sq-1} and integrating over
$(0,t)$ yield
\begin{align*}
	\norm{\nabla \Delta \bff{u}_n(t)}{\bb{L}^2}^2 + \int_0^t \norm{\nabla \Delta^2 \bff{u}_n(s)}{\bb{L}^2}^2 \ds
	&\lesssim
	1 + \int_0^t \norm{\nabla \Delta \bff{u}_n(s)}{\bb{L}^2}^2 \ds + \int_0^t \norm{\Delta^2 \bff{u}_n(s)}{\bb{L}^2}^2 \ds
	\\
	&\lesssim
	1,
\end{align*}
where in the last step we used Proposition \ref{pro:nab delta un} and
\ref{pro:Del sqr}. This completes the proof.
\end{proof}

\begin{proposition}\label{pro:dt un H1}
	{Under the assumption of Proposition \ref{pro:nab Del sq}}, we have
	\begin{align*} 
		\norm{\partial_t \bff{u}_n}{L^2(0,T^\ast; \bb{H}^1)}^2
		{\lesssim
		\norm{\bff{u}_0}{\bb{H}^3}^2,}
	\end{align*}
	where the constant depends on $T^\ast$, but is independent of $n$. 
\end{proposition}

\begin{proof}
	Taking the inner product of \eqref{approx} with $-\Delta \partial_t \bff{u}_n$ and integrating by parts with respect to $\bff{x}$, we have
	\begin{align}\label{equ:dt un-1}
		&\norm{\nabla \partial_t \bff{u}_n}{\bb{L}^2}^2 
		+ 
		\beta_1 \ddt \norm{\Delta \bff{u}_n}{\bb{L}^2}^2 
		+ 
		\beta_2 \ddt \norm{\nabla \Delta \bff{u}_n}{\bb{L}^2}^2 
		\nonumber\\
		\nonumber
		&=
		\beta_3 \ddt \norm{\nabla \bff{u}_n}{\bb{L}^2}^2 
		-
		\beta_3 \inpro{ \nabla \big(|\bff{u}_n|^2 \bff{u}_n\big)}{ \nabla \partial_t \bff{u}_n }_{\bb{L}^2} 
		- 
		\beta_4 \inpro{ \nabla \bff{u}_n \times \Delta \bff{u}_n}{ \nabla \partial_t \bff{u}_n}_{\bb{L}^2} 
		\\
		&\quad
		-
		\beta_4 \inpro{\bff{u}_n \times \nabla \Delta \bff{u}_n}{ \nabla \partial_t \bff{u}_n}_{\bb{L}^2} 
		+
		\beta_5 \inpro {\nabla\Delta (|\bff{u}_n|^2 \bff{u}_n)}{ \nabla\partial_t \bff{u}_n }_{\bb{L}^2}.
	\end{align}
	Each inner product on the right-hand side can be estimated as follows. For the first
	inner product, by H\"{o}lder's inequality, Proposition \ref{pro:Delta un}, \ref{pro:nab delta un} and \ref{pro:partial un}, and Sobolev embedding $\bb{H}^1~\subset~\bb{L}^6$, we have
	\begin{align} \label{equ:dt un inpro 1}
		\nonumber
		\beta_3 \big|
		\inpro{ \nabla \big(|\bff{u}_n|^2 \bff{u}_n\big)}{ \nabla \partial_t \bff{u}_n }_{\bb{L}^2} 
		\big|
		&\lesssim
		\norm{\nabla \bff{u}_n}{\bb{L}^6}
		\norm{\bff{u}_n}{\bb{L}^6}^2
		\norm{\nabla \partial_t \bff{u}_n}{\bb{L}^2}
		\\
		\nonumber
		&\leq
		\norm{\bff{u}_n}{\bb{H}^2}
		\norm{\bff{u}_n}{\bb{H}^1}^2
		\norm{\nabla \partial_t \bff{u}_n}{\bb{L}^2}
		\\
		&\lesssim
		1
		+
		\epsilon
		\norm{\nabla \partial_t \bff{u}_n}{\bb{L}^2}^2,
	\end{align}
	for any $\epsilon >0$, where in the last step we used Young's inequality.
	For the second inner product, H\"{o}lder's inequality, Proposition \ref{pro:Delta un}, \ref{pro:nab delta un}, \ref{pro:nab Del sq}, and Sobolev embedding $\bb{H}^2\subset \bb{L}^\infty$ give
	\begin{align} \label{equ:dt un inpro 2}
		\nonumber
		\beta_4 \big|
		\inpro{ \nabla \bff{u}_n \times \Delta \bff{u}_n}{ \nabla \partial_t \bff{u}_n}_{\bb{L}^2} 
		\big|
		&\lesssim
		\norm{\nabla \bff{u}_n}{\bb{L}^\infty}
		\norm{\Delta \bff{u}_n}{\bb{L}^2}
		\norm{\nabla \partial_t \bff{u}_n}{\bb{L}^2}
		\\
		\nonumber
		&\leq
		\norm{\nabla \bff{u}_n}{\bb{H}^2}
		\norm{\Delta \bff{u}_n}{\bb{L}^2}
		\norm{\nabla \partial_t \bff{u}_n}{\bb{L}^2}
		\\
		&\lesssim
		1 + \epsilon \norm{\nabla \partial_t \bff{u}_n}{\bb{L}^2}^2,
	\end{align}
	for any $\epsilon>0$. Similarly, for the third inner product, we have
	\begin{align}\label{equ:dt un inpro 3}
		\nonumber
		\beta_4 \big|
		\inpro{\bff{u}_n \times \nabla \Delta \bff{u}_n}{ \nabla \partial_t \bff{u}_n}_{\bb{L}^2}
		\big|
		&\lesssim
		\norm{\bff{u}_n}{\bb{L}^\infty}
		\norm{\nabla \Delta \bff{u}_n}{\bb{L}^2}
		\norm{\nabla \partial_t \bff{u}_n}{\bb{L}^2}
		\\
		\nonumber
		&\leq
		\norm{\bff{u}_n}{\bb{H}^2}
		\norm{\nabla \Delta \bff{u}_n}{\bb{L}^2}
		\norm{\nabla \partial_t \bff{u}_n}{\bb{L}^2}
		\\
		&\lesssim
		1 + \epsilon \norm{\nabla \partial_t \bff{u}_n}{\bb{L}^2}^2.
	\end{align}
	For the last inner product on \eqref{equ:dt un-1}, by
	H\"{o}lder's and Young's inequality, we obtain
	\begin{align}\label{equ:dt un inpro 4}
		\nonumber
		\beta_5 \big|
		\inpro {\nabla\Delta (|\bff{u}_n|^2 \bff{u}_n)}{ \nabla\partial_t \bff{u}_n }_{\bb{L}^2}
		\big|
		&\leq
		\norm{|\bff{u}_n|^2 \bff{u}_n}{\bb{H}^3} 
		\norm{\nabla \partial_t \bff{u}_n}{\bb{L}^2}
		\\
		\nonumber
		&\leq
		C \norm{\bff{u}_n}{\bb{H}^3}^3
		\norm{\nabla \partial_t \bff{u}_n}{\bb{L}^2}
		\\
		&\lesssim
		1
		+ 
		\epsilon \norm{\nabla \partial_t \bff{u}_n}{\bb{L}^2}^2,
	\end{align}
	for any $\epsilon >0$, where in the penultimate step we used~\eqref{eq7}, and in
	the last step we used \eqref{eq5}, Proposition \ref{pro:Delta un}, \ref{pro:nab delta un}, \ref{pro:nab Del sq}, and Young's inequality.
	Inserting the estimates \eqref{equ:dt un inpro 1}, \eqref{equ:dt un inpro 2}, \eqref{equ:dt un inpro 3} and \eqref{equ:dt un inpro 4} into \eqref{equ:dt un-1}, integrating over
	$(0,t)$, and choosing $\epsilon$ sufficiently small yield
	\begin{align*}
		\int_0^t \norm{\nabla \partial_s \bff{u}_n(s)}{\bb{L}^2}^2 \ds
		\lesssim
		1 + \norm{\nabla \bff{u}_n}{\bb{L}^2}^2
		\lesssim
		1,
	\end{align*}
	where in the last step we used Proposition \ref{pro:nab delta un}. This completes the proof.
\end{proof}

As a consequence of Proposition~\ref{pro:Delta un}--Proposition~\ref{pro:dt un H1}, we have the following result.
\begin{corollary}\label{cor:un bou}
	For any $T>0$, let $T^\ast$ be defined by Proposition~\ref{pro:nab delta
	un}. Assume that the initial data~$\bff{u}_0$
	satisfies~$\bff{u}_0\in\bb{H}^r$ for~$r\in\{0,1,2,3\}$. Assume further that
\[
	\norm{\bff{u}_{0n}-\bff{u}_0}{\bb{H}^r} \to 0 \quad \text{as
	$n\to\infty$,}
\]
where $\bff{u}_{0,n}$ is defined in~\eqref{approx}. 
{
Then
\begin{equation}\label{equ:un bou}
	\norm{\bff{u}_n}{L^\infty(0,\overline{T};\bb{H}^r)}
	+
	\norm{\bff{u}_n}{L^4(0,\overline{T};\bb{L}^4)}
	+
	\norm{\bff{u}_n}{L^2(0,\overline{T};\bb{H}^{r+2})}
	+
	\norm{\partial_t\bff{u}_n}{L^2(0,\overline{T};\bb{H}^{r-2})}
	\lesssim 1,
\end{equation}
where
\begin{align}\label{equ:T bar alpha}
	\overline{T}=
	\begin{cases}
		T \quad & \text{if $r=0$},
		\\[1ex]
		T^\ast \quad & \text{if $r>0$},
	\end{cases}
\end{align}
}
\end{corollary}

{
\begin{proof}
	First we recall from~\eqref{equ:un0} that the given assumption yields
	$\norm{\bff{u}_n(0)}{\bb{H}^r}\lesssim \norm{\bff{u}_0}{\bb{H}^r}
	\lesssim 1$. 
Therefore,
Proposition~\ref{pro:Delta un} and Proposition~\ref{pro:dt un H-2} imply
\eqref{equ:un bou} when $r=0$, while Proposition~\ref{pro:nab delta un}, Proposition~\ref{pro:partial un}, Proposition~\ref{pro:dt un H-1} and
inequality~\eqref{eq5} give the result when $r=1$.

Next, Proposition~\ref{pro:Delta un}, Proposition~\ref{pro:nab delta un},
Proposition~\ref{pro:partial un}, Proposition~\ref{pro:Del sqr} and inequality~\eqref{eq5} give the required
estimate for the case $r=2$. Finally,  Proposition~\ref{pro:Delta un}, Proposition~\ref{pro:nab delta un},
Proposition~\ref{pro:partial un}, Proposition~\ref{pro:Del sqr},
Proposition~\ref{pro:nab Del sq}, Proposition~\ref{pro:dt un H1} and inequality~\eqref{eq8} give the result for
the case $r=3$,  completing the proof of the corollary.
\end{proof}
}

\section{Proof of Theorem \ref{the:weakexist}} \label{weaksol}

Let $r\in \{0,1,2,3\}$. It follows from~\eqref{equ:un bou} and the Banach-Alaoglu
theorem that there exists a subsequence of~$\{\bff{u}_n\}$, which is still denoted by
$\{\bff{u}_n\}$, such that
{
\begin{equation}\label{equ:wea con L2}
	\left\{
	\begin{alignedat}{2}
		\bff{u}_n &\rightharpoonup \bff{u} \quad &&\text{weakly* in } L^\infty (0,\overline{T}; \bb{H}^r),
		\\
		\bff{u}_n &\rightharpoonup \bff{u} \quad &&\text{weakly in }
		L^2(0,\overline{T}; \bb{H}^{r+2}),
			\\
		\bff{u}_n &\rightharpoonup \bff{u} \quad &&\text{weakly in }
		L^4(0,\overline{T}; \bb{L}^4),
		\\
		\partial_t \bff{u}_n &\rightharpoonup \partial_t \bff{u} 
		\quad &&\text {weakly in } L^2(0,\overline{T}; \bb{H}^{r-2}),
	\end{alignedat}
	\right.
\end{equation}
where $\overline{T}$ was defined in \eqref{equ:T bar alpha}.
By the Aubin--Lions--Simon lemma (Theorem \ref{aubin}), a further subsequence then satisfies
\begin{equation}\label{equ:unu str L2}
	\bff{u}_n \to \bff{u} \quad \text{strongly in } 
	L^2(0,T; \bb{H}^1).
\end{equation}
}

The next proposition shows the convergence of the nonlinear terms in \eqref{approx}.
\begin{proposition} \label{pro:limit L2}
	Let $T>0$ be arbitrary. Let $\{\bff{\phi}_n\}$ be a sequence in $\bb{V}_n$ such that $\bff{\phi}_n \to \bff{\phi}$ in $\bb{H}^2$. For all $t\in [0,\overline{T}]$, we have
	\begin{align}
		%\lim_{n\to\infty} \int_0^t \inpro{ \nabla \bff{u}_n(s)}{ \nabla \bff{\phi}_n}_{\bb{L}^2} \ds 
		%&=
		%\int_0^t \inpro{ \nabla \bff{u}(s)}{ \nabla \bff{\phi}}_{\bb{L}^2} \ds, \label{equ:limit1}
		%\\
		%\lim_{n\to\infty}\int_0^t \inpro{ \Delta \bff{u}_n(s)}{ \Delta \bff{\phi}_n}_{\bb{L}^2} \ds 
		%&= 
		%\int_0^t \inpro{ \Delta \bff{u}(s)}{ \Delta \bff{\phi}}_{\bb{L}^2} \ds, \label{equ:limit2}
		%\\
		\label{equ:limit3 L2}
		\lim_{n\to\infty} \int_0^t \inpro{ \Pi_n\left((1-|\bff{u}_n(s)|^2)
			\bff{u}_n(s)\right)}{\bff{\phi}_n }_{\bb{L}^2} \ds 
		&= 
		\int_0^t \inpro{ \left((1-|\bff{u}(s)|^2) \bff{u}(s)\right)}{\bff{\phi}
		}_{\bb{L}^2} \ds, 
		\\
		 \label{equ:limit4 L2}
		\lim_{n\to\infty} \int_0^t \inpro{ \Pi_n \left(\bff{u}_n(s)\times \Delta\bff{u}_n(s)\right)}{ \bff{\phi}_n }_{\bb{L}^2} \ds 
		&=
		- \int_0^t \inpro{ \left(\bff{u}(s)\times \nabla \bff{u}(s)\right)}{
			\nabla \bff{\phi} }_{\bb{L}^2} \ds,
		\\
		\label{equ:limit5 L2}
		\lim_{n\to\infty} \int_0^t \inpro{ \Pi_n \left(\Delta \left(|\bff{u}_n(s)|^2 \bff{u}_n(s)\right)\right)}{ \bff{\phi}_n }_{\bb{L}^2} \ds 
		&= 
		\int_0^t \inpro{ \nabla \left(|\bff{u}(s)|^2 \bff{u}(s)\right)}{\nabla \bff{\phi} }_{\bb{L}^2} \ds.
	\end{align}
\end{proposition}

\begin{proof}
	By the definition of $\Pi_n$, in order to prove~\eqref{equ:limit3 L2} it suffices to show
	\begin{equation}\label{equ:uns u L2}
		\lim_{n\to\infty} \int_0^t \inpro{ |\bff{u}_n(s)|^2
			\bff{u}_n(s)}{\bff{\phi}_n }_{\bb{L}^2} \ds
		= 
		\int_0^t \inpro{ |\bff{u}(s)|^2 \bff{u}(s)}{\bff{\phi} }_{\bb{L}^2} \ds.
	\end{equation}
	
	To this end, note that H\"older's inequality implies
	\begin{align*}
		&\left| \int_0^t \inpro{ |\bff{u}_n(s)|^2 \bff{u}_n(s)}{ \bff{\phi}_n}_{\bb{L}^2} \ds 
		- 
		\int_0^t \inpro{ |\bff{u}(s)|^2 \bff{u}(s)}{ \bff{\phi}}_{\bb{L}^2} \ds \right|
		\\
		&\leq 
		\left| \int_0^t \inpro{ |\bff{u}_n(s)|^2 \bff{u}_n(s)}{ \bff{\phi}_n-\bff{\phi} }_{\bb{L}^2} \ds\right| 
		+ 
		\left| \int_0^t \inpro{ |\bff{u}_n(s)|^2 (\bff{u}_n(s)-\bff{u}(s))} {\bff{\phi}}_{\bb{L}^2} \ds \right| 
		\\
		&\quad +
		\left| \int_0^t \inpro{ (|\bff{u}_n(s)|^2 - |\bff{u}(s)|^2) \bff{u}(s)}{ \bff{\phi}}_{\bb{L}^2} \ds \right|
		\\
		&\leq
		\norm{\bff{\phi}_n - \bff{\phi}}{\bb{L}^6} 
		\int_0^t \norm{\bff{u}_n(s)}{\bb{L}^6}^2
		\norm{\bff{u}_n(s)}{\bb{L}^2} \ds
		+
		\norm{\bff{\phi}}{\bb{L}^\infty} 
		\int_0^t \norm{\bff{u}_n(s)}{\bb{L}^4}^2  
		\norm{\bff{u}_n(s) -\bff{u}(s)}{\bb{L}^2}
		\ds
		\\
		&\quad +
		\norm{\bff{\phi}}{\bb{L}^\infty} 
		\int_0^t \norm{\bff{u}_n(s) -\bff{u}(s)}{\bb{L}^2}  \norm{\bff{u}_n(s) +
			\bff{u}(s)}{\bb{L}^4}  \norm{\bff{u}(s)}{\bb{L}^4}
		\ds
		\\
		&\leq
		\norm{\bff{\phi}_n-\bff{\phi}}{\bb{H}^1}
		\norm{\bff{u}_n}{L^2(0,\overline{T};\bb{H}^1)}
		+
		\norm{\bff{\phi}}{\bb{H}^2}
		\norm{\bff{u}_n}{\bb{L}^4(0,\overline{T};\bb{L}^4)}
		\norm{\bff{u}_n-\bff{u}}{\bb{L}^2(0,\overline{T};\bb{L}^2)}
		\\
		&\quad
		+
		\norm{\bff{\phi}}{\bb{H}^2}
		\norm{\bff{u}_n-\bff{u}}{L^2(0,\overline{T};\bb{L}^2)}
		\norm{\bff{u}_n+\bff{u}}{L^4(0,\overline{T};\bb{L}^4)}
		\norm{\bff{u}}{L^4(0,\overline{T};\bb{L}^4)}.
	\end{align*}
	By using the Sobolev embedding $\bb{H}^1 \subset \bb{L}^6$ and $\bb{H}^2 \subset \bb{L}^\infty$, \eqref{equ:un bou},
	and \eqref{equ:unu str L2} we deduce~\eqref{equ:uns u L2}.
	
	Similarly, to show~\eqref{equ:limit4 L2}, it suffices to show
	\begin{align*}
		\lim_{n\to\infty} \int_0^t \inpro{ \left(\bff{u}_n(s)\times \nabla \bff{u}_n(s)\right)}{ \nabla \bff{\phi}_n }_{\bb{L}^2} \ds 
		&=
		\int_0^t \inpro{ \left(\bff{u}(s)\times \nabla \bff{u}(s)\right)}{ \nabla \bff{\phi} }_{\bb{L}^2} \ds.
	\end{align*}
	To this end, note that
	\begin{align*}
		&\left|\int_0^t \inpro{ \bff{u}_n(s) \times \nabla \bff{u}_n(s)}{ \nabla \bff{\phi}_n}_{\bb{L}^2} \ds 
		- 
		\int_0^t \inpro{ \bff{u}(s) \times \nabla\bff{u}(s) }{\nabla \bff{\phi}}_{\bb{L}^2} \ds \right|
		\\
		&\leq 
		\left| \int_0^t \inpro{ \bff{u}_n(s) \times \nabla \bff{u}_n(s)}{ \nabla \bff{\phi}_n - \nabla \bff{\phi} }_{\bb{L}^2} \ds \right| 
		+ 
		\left| \int_0^t \inpro{ (\bff{u}_n(s) -\bff{u}(s)) \times \nabla \bff{u}_n(s)}{ \nabla \bff{\phi}}_{\bb{L}^2} \ds \right|
		\\
		&\quad + 
		\left| \int_0^t \inpro{ \bff{u}(s) \times (\nabla \bff{u}_n(s) -\nabla \bff{u}(s))}{ \nabla \bff{\phi} }_{\bb{L}^2} \ds \right|
		\\
		&\lesssim 
		\norm{\bff{u}_n}{L^2(0,\overline{T};\bb{L}^4)} 
		\norm{\nabla \bff{u}_n}{L^2(0,\overline{T};\bb{L}^4)} 
		\norm{\nabla \bff{\phi}_n - \nabla \bff{\phi}}{\bb{L}^2}
		\\
		&\quad
		+
		\norm{\bff{u}_n -\bff{u}}{L^2(0,\overline{T};\bb{L}^4)} 
		\norm{\nabla \bff{u}_n}{L^2(0,\overline{T};\bb{L}^4)} 
		\norm{\nabla \bff{\phi}}{\bb{L}^2} 
		\\
		&\quad + 
		\norm{\bff{u}}{L^2(0,\overline{T};\bb{L}^4)}
		\norm{\nabla \bff{u}_n-\nabla \bff{u}}{L^2(0,\overline{T};\bb{L}^4)} 
		\norm{\nabla \bff{\phi}}{\bb{L}^2}.
		%\\
		%&\lesssim
		%\norm{\nabla \bff{\phi}_n - \nabla \bff{\phi}}{\bb{L}^2} 
		%+ 
		%\norm{\bff{u}_n -\bff{u}}{L^2(0,T;\bb{H}^2)},
	\end{align*}
	By using the Sobolev embedding $\bb{H}^1 \subset \bb{L}^4$, \eqref{equ:un bou},
	and~\eqref{equ:unu str L2}, we deduce the required convergence.
	
	For the last convergence~\eqref{equ:limit5 L2}, it suffices to show
	\begin{align*}
		\lim_{n\to\infty} \int_0^t \inpro{ \nabla \left(|\bff{u}_n(s)|^2 \bff{u}_n(s)\right)}{\nabla \bff{\phi}_n }_{\bb{L}^2} \ds 
		&= 
		\int_0^t \inpro{ \nabla \left(|\bff{u}(s)|^2 \bff{u}(s)\right)}{\nabla \bff{\phi} }_{\bb{L}^2} \ds.
	\end{align*}
	To this end, note that
	\begin{align*}
		\nonumber
		&\left| \int_0^t \inpro{ \nabla (|\bff{u}_n(s)|^2 \bff{u}_n(s))}{ \nabla \bff{\phi}_n }_{\bb{L}^2} \ds 
		- 
		\int_0^t \inpro{ \nabla (|\bff{u}(s)|^2 \bff{u}(s))}{ \nabla \bff{\phi} }_{\bb{L}^2} \ds \right|
		\\
		\nonumber
		&\leq 
		\left| \int_0^t \inpro{ \nabla(|\bff{u}_n(s)|^2 \bff{u}_n(s))}{ \nabla\bff{\phi}_n - \nabla\bff{\phi} }_{\bb{L}^2} \ds\right| 
		\\
		\nonumber
		&\quad 
		+ \left| \int_0^t \inpro{ \nabla(|\bff{u}_n(s)|^2 \bff{u}_n(s)) - \nabla(|\bff{u}_n(s)|^2 \bff{u}(s))}{ \nabla\bff{\phi}}_{\bb{L}^2} \ds \right|
		\\
		&\quad 
		+\left| \int_0^t \inpro{ \nabla(|\bff{u}_n(s)|^2 \bff{u}(s)) -
			\nabla(|\bff{u}(s)|^2 \bff{u}(s))}{ \nabla \bff{\phi}}_{\bb{L}^2} \ds
		\right|. 
	\end{align*}
	The arguments follow along the line of the previous convergence statements and
	are omitted.
\end{proof}

We are now ready to prove Theorem~\ref{the:weakexist}.

\bigskip
\noindent
\underline{Proof that $\bff{u}$ satisfies \eqref{weakform} and \eqref{equ:wea sol smo}}:
For any $\bff{\phi} \in \bb{H}^2$, take a sequence $\{\bff{\phi}_n\}$ in
$\bb{V}_n$ such that $\bff{\phi}_n \to \bff{\phi}$ in~$\bb{H}^2$. It follows
from~\eqref{approx} that
\begin{align*}
	&\inpro{\bff{u}_n(t)}{\bff{\phi}_n}_{\bb{L}^2} 
	+
	\beta_1 \int_0^t \inpro{ \nabla \bff{u}_n(s)}{ \nabla \bff{\phi}_n}_{\bb{L}^2} \ds 
	+
	\beta_2 \int_0^t \inpro{ \Delta \bff{u}_n(s)}{ \Delta \bff{\phi}_n}_{\bb{L}^2} \ds 
	\\
	&= 
	\inpro{\bff{u}_{0n}}{\bff{\phi}_n}_{\bb{L}^2}
	+
	\beta_3 \int_0^t \inpro{ \Pi_n\left((1-|\bff{u}_n(s)|^2) \bff{u}_n(s)\right)}{\bff{\phi}_n }_{\bb{L}^2} \ds
	\\
	&\quad
	+
	\beta_4 \int_0^t \inpro{ \Pi_n \left(\bff{u}_n(s)\times \Delta\bff{u}_n(s)\right)}{ \bff{\phi}_n }_{\bb{L}^2} \ds 
	-
	\beta_5 \int_0^t \inpro{ \Pi_n \left(\Delta \left(|\bff{u}_n(s)|^2 \bff{u}_n(s)\right)\right)}{ \bff{\phi}_n }_{\bb{L}^2} \ds.
\end{align*}
Hence, letting $n \to \infty$ and using Proposition \ref{pro:limit L2} we deduce~\eqref{weakform}.
{
Noting \eqref{equ:wea con L2}, we have
	\begin{equation*}
		\bff{u} \in L^\infty(0,\overline{T}; \bb{H}^r) \cap L^2(0,\overline{T}; \bb{H}^{r+2})
			\cap L^4(0,\overline{T}; \bb{L}^4) 
			\quad
			\text{and}
			\quad
		\partial_t \bff{u} \in L^2(0, \overline{T}; \bb{H}^{r-2}).
	\end{equation*}
Therefore, applying Theorem \ref{the:E22 to C}, and noting that $[\bb{H}^{r-2},
\bb{H}^{r+2}]_{1/2} \equiv \bb{H}^r$, we obtain \eqref{equ:wea sol smo}.

\bigskip
\noindent
\underline{Proof of \eqref{equ:cont dep u-v L2}}:
Let $\bff{u}$ and $\bff{v}$ be weak solutions to \eqref{equ:formal} with initial data $\bff{u}_0$ and $\bff{v}_0 \in \bb{L}^2$, respectively. Let
$\bff{w}=\bff{u}-\bff{v}$. Then, for all $\bff{\phi}\in \bb{H}^2$,
\begin{align}\label{equ:difference}
	\nonumber
	&\inpro{ \partial_t \bff{w}}{ \bff{\phi}}_{\bb{L}^2} 
	+ 
	\beta_1 \inpro{ \nabla \bff{w}}{ \nabla \bff{\phi} }_{\bb{L}^2} 
	+
	\beta_2 \inpro{ \Delta \bff{w}}{ \Delta \bff{\phi} }_{\bb{L}^2} 
	\\
	\nonumber
	&= 
	\beta_3 \inpro{ \bff{w}}{ \bff{\phi} }_{\bb{L}^2} 
	- 
	\beta_3 \inpro{ |\bff{u}|^2 \bff{u} - |\bff{v}|^2 \bff{v}}{ \bff{\phi} }_{\bb{L}^2} 
	+ 
	\beta_4 \inpro{ \bff{u} \times \nabla \bff{u} 
		-  \bff{v} \times \nabla \bff{v}}{ \nabla \bff{\phi} }_{\bb{L}^2} 
	\\
	&\quad - 
	\beta_5 \inpro{ \nabla (|\bff{u}|^2 \bff{u}- |\bff{v}|^2 \bff{v})}{ \nabla \bff{\phi} }_{\bb{L}^2}.
\end{align}
By using integration by parts for the terms with coefficient~$\beta_4$
and~$\beta_5$, we obtain
\begin{align}\label{equ:dif 2}
	\nonumber
	&\inpro{ \partial_t \bff{w}}{ \bff{\phi}}_{\bb{L}^2} 
	+ 
	\beta_1 \inpro{ \nabla \bff{w}}{ \nabla \bff{\phi} }_{\bb{L}^2} 
	+
	\beta_2 \inpro{ \Delta \bff{w}}{ \Delta \bff{\phi} }_{\bb{L}^2} 
	\\
	\nonumber
	&= 
	\beta_3 \inpro{ \bff{w}}{ \bff{\phi} }_{\bb{L}^2} 
	- 
	\beta_3 \inpro{ |\bff{u}|^2 \bff{u} - |\bff{v}|^2 \bff{v}}{ \bff{\phi} }_{\bb{L}^2} 
	- 
	\beta_4 \inpro{ \bff{u} \times \Delta \bff{u} 
		-  \bff{v} \times \Delta \bff{v}}{ \bff{\phi} }_{\bb{L}^2} 
	\\
	&\quad + 
	\beta_5 \inpro{ |\bff{u}|^2 \bff{u}- |\bff{v}|^2 \bff{v}}{ \Delta \bff{\phi} }_{\bb{L}^2}.
\end{align}
Both equations above will be used at our convenience.
Letting $\bff{\phi}= \bff{w}$ in~\eqref{equ:dif 2}, we have
\begin{align}\label{equ:unique0 weak}
	\nonumber
	\frac{1}{2} \ddt \norm{\bff{w}}{\bb{L}^2}^2 
	+ 
	\beta_2 \norm{\Delta \bff{w}}{\bb{L}^2}^2  
	&\le
	|\beta_1| \norm{\nabla \bff{w}}{\bb{L}^2}^2 
	+
	\beta_3 
	\norm{\bff{w}}{\bb{L}^2}^2 
	+ 
	\beta_3 
	\Big|
	\inpro{ |\bff{u}|^2 \bff{u} - |\bff{v}|^2 \bff{v}}{ \bff{w} }_{\bb{L}^2} 
	\Big|
	\nonumber
	\\
	&\quad +
	\beta_4 
	\Big|
	\inpro{ \bff{u} \times \Delta \bff{u}- \bff{v} \times \Delta \bff{v}}{ \bff{w} }_{\bb{L}^2} 
	\Big|
	\nonumber\\
	&\quad 
	+
	\beta_5 
	\Big|
	\inpro{ |\bff{u}|^2 \bff{u}- |\bff{v}|^2 \bff{v}}{ \Delta \bff{w} }_{\bb{L}^2}
	\Big|.
\end{align}

We will now estimate the inner products on the right-hand side. 
For the first inner product, applying \eqref{equ:D k u2u-v2v} yields
\begin{align} \label{equ:unique2 weak}
	\beta_3 
	\big|
	\inpro{ |\bff{u}|^2 \bff{u} - |\bff{v}|^2 \bff{v}}{ \bff{w} }_{\bb{L}^2} 
	\big|
	&\lesssim 
	\big( \norm{\bff{u}}{\bb{L}^\infty}^2 + \norm{\bff{v}}{\bb{L}^\infty}^2 \big)
	\norm{\bff{w}}{\bb{L}^2}^2.
\end{align}
	For the second inner product on the right-hand side
	of~\eqref{equ:unique0 weak}, we have by using H\"{o}lder's inequality
	and Young's inequality
\begin{align}\label{equ:unique3 weak}
	\beta_4 
	\big|
	\inpro{ \bff{u} \times \Delta \bff{u}- \bff{v} \times \Delta \bff{v}}{ \bff{w} }_{\bb{L}^2} 
	\big|
	&= 
	\beta_4 
	\big|
	\inpro{ \bff{u} \times \Delta \bff{w} + \bff{w} \times \Delta \bff{v}}{\bff{w} }_{\bb{L}^2} 
	\big|
	= 
	\beta_4
	\big|
	\inpro{ \bff{u} \times \Delta \bff{w}}{ \bff{w} }_{\bb{L}^2} 
	\big|
	\nonumber\\
	&\lesssim 
	\norm{\bff{u}}{\bb{L}^\infty}
	\norm{\Delta \bff{w}}{\bb{L}^2}
	\norm{\bff{w}}{\bb{L}^2}
	\nonumber\\
	&\leq
	C 
	\norm{ \bff{u}}{\bb{L}^\infty}^2 \norm{\bff{w}}{\bb{L}^2}^2
	+ 
	\epsilon\norm{\Delta\bff{w}}{\bb{L}^2}^2,
\end{align}
for any $\epsilon>0$.
For the last inner product in~\eqref{equ:unique0 weak}, applying \eqref{equ:D k u2u-v2v}, then using Young's inequality yield
\begin{align}\label{equ:unique4 weak}
	\beta_5 
	\Big|
	\inpro{ |\bff{u}|^2 \bff{u}-  |\bff{v}|^2 \bff{v}}{ \Delta \bff{w} }_{\bb{L}^2} 
	\Big|
	&\leq
	C \big(\norm{\bff{u}}{\bb{L}^\infty}^4 
	+ \norm{\bff{v}}{\bb{L}^\infty}^4 \big) \norm{\bff{w}}{\bb{L}^2}^2
	+ 
	\epsilon \norm{\Delta \bff{w}}{\bb{L}^2}^2
\end{align}
for any $\epsilon > 0$.
Inserting \eqref{equ:unique2 weak}, \eqref{equ:unique3 weak} and
\eqref{equ:unique4 weak} into \eqref{equ:unique0 weak}, and choosing $\epsilon$
sufficiently small, we obtain
\begin{align}\label{equ:ddt unique}
	\nonumber
	\ddt \norm{\bff{w}}{\bb{L}^2}^2
	+
	\norm{\Delta \bff{w}}{\bb{L}^2}^2
	&\lesssim 
	\norm{\nabla \bff{w}}{\bb{L}^2}^2
	+
	\big(\norm{\bff{u}}{\bb{L}^\infty}^2 + \norm{\bff{v}}{\bb{L}^\infty}^2 +
	\norm{\bff{u}}{\bb{L}^\infty}^4 + \norm{\bff{v}}{\bb{L}^\infty}^4 \big)
	\norm{\bff{w}}{\bb{L}^2}^2
	+
	\epsilon \norm{\Delta \bff{w}}{\bb{L}^2}^2
	\\
	\nonumber
	&\lesssim
	\big(1+ \norm{\bff{u}}{\bb{L}^\infty}^2 + \norm{\bff{v}}{\bb{L}^\infty}^2 +
	\norm{\bff{u}}{\bb{L}^\infty}^4 + \norm{\bff{v}}{\bb{L}^\infty}^4 \big)
	\norm{\bff{w}}{\bb{L}^2}^2
	+
	\epsilon \norm{\Delta \bff{w}}{\bb{L}^2}^2
	\\
	&\lesssim
	\big(1+
	\norm{\bff{u}}{\bb{L}^\infty}^4 + \norm{\bff{v}}{\bb{L}^\infty}^4 \big)
	\norm{\bff{w}}{\bb{L}^2}^2
	+
	\epsilon \norm{\Delta \bff{w}}{\bb{L}^2}^2,
\end{align}
where we used \eqref{eq2} and Young's inequality. Choosing $\epsilon$ sufficiently
small, rearranging the above equation, integrating over $(0,t)$, and using
Gronwall's inequality, we infer~\eqref{equ:cont dep u-v L2}. Under
assumption~\eqref{equ:add reg unique}, uniqueness then follows.

\bigskip
\noindent
\underline{Proof of \eqref{equ:uv uv0 Hr}}:

\medskip
\underline{The case $r=1$}:
Taking $\bff{\phi}= -\Delta \bff{w}$ in \eqref{equ:difference}, and integrating
by parts (for the terms with coefficient~$\beta_1$ and~$\beta_2$) we have (after
rearranging the terms)
	\begin{align}\label{equ:cont dep H1}
		\nonumber
		&\frac{1}{2} \ddt \norm{\nabla \bff{w}}{\bb{L}^2}^2
		+
		\beta_2 \norm{\nabla \Delta \bff{w}}{\bb{L}^2}^2
		\\
		\nonumber
		&\leq
		|\beta_1| \norm{\Delta \bff{w}}{\bb{L}^2}^2 
		+
		\beta_3 \norm{\nabla \bff{w}}{\bb{L}^2}^2 
		+
		\beta_3 \big| \inpro{|\bff{u}|^2 \bff{u}- |\bff{v}|^2
		\bff{v}}{\Delta \bff{w}}_{\bb{L}^2} \big|
		\\
		&\quad
		+
		\beta_4 
		\big| \inpro{\bff{u} \times \nabla \bff{u} - \bff{v} \times
		\nabla \bff{v}}{\nabla\Delta \bff{w}}_{\bb{L}^2} \big|
		+
		\beta_5 
		\big| \inpro{\nabla(|\bff{u}|^2 \bff{u}- |\bff{v}|^2
		\bff{v})}{\nabla\Delta \bff{w}}_{\bb{L}^2} \big|.
	\end{align}
	We will estimate the inner products on the right-hand side. For the
	first inner product, it follows successively from~\eqref{equ:unique4 weak},
	\eqref{eq3}, and Young's inequality that
	\begin{align}\label{equ:cont dep H1 1st}
		\beta_3 \big| \inpro{ |\bff{u}|^2 \bff{u}- |\bff{v}|^2 \bff{v}}{\Delta \bff{w}}_{\bb{L}^2} \big|
		&\lesssim
		\big(\norm{\bff{u}}{\bb{L}^\infty}^4 + \norm{\bff{v}}{\bb{L}^\infty}^4 \big) \norm{\bff{w}}{\bb{L}^2}^2
		+
		\norm{\nabla \bff{w}}{\bb{L}^2}^2 
		+
		\epsilon \norm{\nabla \Delta \bff{w}}{\bb{L}^2}^2,
	\end{align}
	for any $\epsilon>0$. For the second inner product, applying \eqref{equ:u cross Dk u} then using Young's inequality yield
	\begin{align}\label{equ:cont dep H1 2nd}
		\beta_4 \big| \inpro{\bff{u} \times \nabla \bff{u} - \bff{v}
		\times \nabla \bff{v}}{\nabla\Delta \bff{w}}_{\bb{L}^2} \big|
		&\lesssim
		\norm{\bff{u}}{\bb{L}^\infty}^2 \norm{\nabla \bff{w}}{\bb{L}^2}^2
		+
		\norm{\bff{w}}{\bb{L}^2}^2 \norm{\nabla \bff{v}}{\bb{L}^\infty}^2
		+
		\epsilon \norm{\nabla \Delta \bff{w}}{\bb{L}^2}^2,
	\end{align}
	for any $\epsilon >0$. For the third inner product on the right-hand
	side of~\eqref{equ:cont dep H1}, applying \eqref{equ:D k u2u-v2v} then using Young's inequality give
	\begin{align}\label{equ:cont dep H1 3rd}
		\beta_5 
		\big| \inpro{\nabla(|\bff{u}|^2 \bff{u}- |\bff{v}|^2
		\bff{v})}{\nabla\Delta \bff{w}}_{\bb{L}^2} \big|
		&\lesssim
		\big( \norm{\bff{u}}{\bb{H}^1}^2 + \norm{\bff{v}}{\bb{H}^1}^2 \big)
		\big( \norm{\bff{u}}{\bb{H}^2}^2 + \norm{\bff{v}}{\bb{H}^2}^2 \big) 
		\norm{\bff{w}}{\bb{H}^1}^2
		\nonumber\\
		&\quad
		+
		\epsilon \norm{\nabla \Delta \bff{w}}{\bb{L}^2}^2,
	\end{align}
	for any $\epsilon >0$.
	Inserting \eqref{equ:cont dep H1 1st}, \eqref{equ:cont dep H1 2nd} and \eqref{equ:cont dep H1 3rd} into \eqref{equ:cont dep H1}, integrating over $(0,t)$, and choosing $\epsilon$ sufficiently small, we obtain
	\begin{align*}
		\norm{\nabla \bff{w}(t)}{\bb{L}^2}^2 
		&+
		\int_0^t \norm{\nabla \Delta \bff{w}(s)}{\bb{L}^2}^2 \ds
		\\
		&\lesssim
		\norm{\nabla \bff{w}(0)}{\bb{L}^2}^2
		+
		\int_0^t 
		\bff{\alpha}(s)
		\norm{\bff{w}(s)}{\bb{L}^2}^2 \ds
		+
		\int_0^t \bff{\alpha}(s) \norm{\nabla \bff{w}(s)}{\bb{L}^2}^2 \ds
		\\
		&\lesssim
		\norm{\nabla \bff{w}(0)}{\bb{L}^2}^2
		+
		\norm{\bff{w}(0)}{\bb{L}^2}^2 
		\int_0^t 
		\bff{\alpha}(s) \exp \left(\int_0^t \bff{\alpha}(\tau) \dtau \right) \ds
		+
		\int_0^t \bff{\alpha}(s) \norm{\nabla \bff{w}(s)}{\bb{L}^2}^2 \ds,
	\end{align*}
	where
	\[
		\bff{\alpha}(s) := 
		\big( \norm{\bff{u}(s)}{\bb{H}^1}^2 + \norm{\bff{v}(s)}{\bb{H}^1}^2 \big)
		\big( \norm{\bff{u}(s)}{\bb{H}^2}^2 + \norm{\bff{v}(s)}{\bb{H}^2}^2 \big)
		+ 
		\norm{\bff{u}(s)}{\bb{H}^3}^2 
		+
		\norm{\bff{v}(s)}{\bb{H}^3}^2 
	\]
	and where in the last step we used~\eqref{equ:cont dep u-v L2} which was
	proved above in this proof.
	Note that for $\bff{u}_0$, $\bff{v}_0 \in \bb{H}^1$, we have $\bff{u}$,
	$\bff{v}\in L^\infty(0,T; \bb{H}^1) \cap L^2(0,T; \bb{H}^3)$, and thus (by Gagliardo--Nirenberg inequality)
	\[
		\int_0^t \bff{\alpha}(s) \ds 
		\lesssim
		\int_0^t \norm{\bff{u}(s)}{\bb{H}^3}^2
		+
		\norm{\bff{v}(s)}{\bb{H}^3}^2 \ds
		\lesssim 1.
	\]
	Gronwall's inequality then yields
	\begin{align*}
		\norm{\nabla \bff{w}(t)}{\bb{L}^2}^2
		\leq
		C \norm{\bff{w}(0)}{\bb{H}^1}^2,
	\end{align*}
which implies~\eqref{equ:uv uv0 Hr}.

\medskip
\underline{The case $r=2$}:
Taking $\bff{\phi}=\partial_t \bff{w}$ in \eqref{equ:dif 2} and using
integration by parts for the term with coefficient~$\beta_5$ give
	\begin{align}\label{equ:cont dep H2}
		\nonumber
		\norm{\partial_t \bff{w}}{\bb{L}^2}^2
		+
		\beta_2 \ddt \norm{\Delta \bff{w}}{\bb{L}^2}^2
		&\leq
		\beta_1 \ddt \norm{\nabla \bff{w}}{\bb{L}^2}^2
		+
		\beta_3 \ddt \norm{\bff{w}}{\bb{L}^2}^2 
		+
		\beta_3 
		\Big|
		\inpro{|\bff{u}|^2 \bff{u}- |\bff{v}|^2 \bff{v}}{\partial_t \bff{w}}_{\bb{L}^2}
		\Big|
		\nonumber\\
		&\quad
		+
		\beta_4 
		\Big|
		\inpro{\bff{u} \times \Delta \bff{u} - \bff{v} \times \Delta \bff{v}}
		{\partial_t \bff{w}}_{\bb{L}^2}
		\Big|
		\nonumber\\
		&\quad
		+
		\beta_5 
		\Big|
		\inpro{\Delta(|\bff{u}|^2 \bff{u}- |\bff{v}|^2
		\bff{v})}{\partial_t \bff{w}}_{\bb{L}^2}
		\Big|.
	\end{align}
Note that when using integration by parts, the integrals on the boundary vanish
due to the boundary property of $\bff{w}$ and~\eqref{equ:nor der v2v}.
	Each inner product on the right-hand side of \eqref{equ:cont dep H2} can
	be estimated as follows. 
	For the first inner product, similarly to~\eqref{equ:unique4 weak} we
	have
	\begin{align}\label{equ:cont dep H2 1st}
		\beta_3 \big| 
		\inpro{|\bff{u}|^2 \bff{u}- |\bff{v}|^2 \bff{v}}{\partial_t \bff{w}}_{\bb{L}^2}
		\big|
		\leq 
		C \big(\norm{\bff{u}}{\bb{L}^\infty}^4 
		+ \norm{\bff{v}}{\bb{L}^\infty}^4 \big) \norm{\bff{w}}{\bb{L}^2}^2
		+ 
		\epsilon \norm{\partial_t \bff{w}}{\bb{L}^2}^2
	\end{align}
	for any $\epsilon>0$. For the second inner product, applying \eqref{equ:u cross Dk u} and Young's inequality yield
	\begin{align}\label{equ:cont dep H2 2nd}
		\beta_4 \big|
		\inpro{\bff{u} \times \Delta \bff{u} - \bff{v} \times \Delta \bff{v}}{\partial_t \bff{w}}_{\bb{L}^2}
		\big|
		&\leq
		C \big( \norm{\bff{u}}{\bb{H}^2}^2 + \norm{\bff{v}}{\bb{H}^2}^2 \big)
		\norm{\bff{w}}{\bb{H}^2}^2
		+
		\epsilon \norm{\partial_t \bff{w}}{\bb{L}^2}^2.
	\end{align}
	For the last inner product in \eqref{equ:cont dep H2}, by \eqref{equ:D k u2u-v2v} and Young's inequality, we have
	\begin{align}\label{equ:cont dep H2 3rd}
		\beta_5 \big| 
		\inpro{\Delta(|\bff{u}|^2 \bff{u}- |\bff{v}|^2 \bff{v})}{\partial_t \bff{w}}_{\bb{L}^2}
		\big| 
		\leq
		C \big( \norm{\bff{u}}{\bb{H}^2}^4 + \norm{\bff{v}}{\bb{H}^2}^4 \big)
		\norm{\bff{w}}{\bb{H}^2}^2
		+
		\epsilon \norm{\partial_t \bff{w}}{\bb{L}^2}^2.
	\end{align}
	Inserting \eqref{equ:cont dep H2 1st}, \eqref{equ:cont dep H2 2nd} and
	\eqref{equ:cont dep H2 3rd} into \eqref{equ:cont dep H2}, integrating
	over $(0,t)$, and choosing $\epsilon>0$ sufficiently small, we obtain
	\begin{align*}
		\int_{0}^{t} 
		\norm{\partial_t \bff{w}(s)}{\bb{L}^2}^2 \ds
		+
		\norm{\Delta \bff{w}(t)}{\bb{L}^2}^2
		&\lesssim
		\norm{\bff{w}(0)}{\bb{H}^2}^2 
		+ 
		\int_0^t \bff{\beta}(s) \norm{\bff{w}(s)}{\bb{L}^2}^2 \ds 
		+
		\int_0^t \bff{\beta}(s) \norm{\Delta \bff{w}(s)}{\bb{L}^2}^2 \ds,
	\end{align*}
	where
	\[
		\bff{\beta}(s) := 1 + \norm{\bff{u}(s)}{\bb{H}^2}^4 + \norm{\bff{v}(s)}{\bb{H}^2}^4.
	\]
	Gronwall's inequality then yields
	\[
		\norm{\Delta \bff{w}(t)}{\bb{L}^2}^2 \leq C
		\norm{\bff{w}(0)}{\bb{H}^2}^2.
	\]
	This, together with \eqref{equ:cont dep u-v L2} shown for $r=0$ and \eqref{equ:uv uv0 Hr} for $r=1$,
	gives the required inequality~\eqref{equ:uv uv0 Hr} for $r=2$.

\bigskip
\underline{The case $r=3$}:
	We now take $\bff{\phi}= -\Delta \partial_t \bff{w}$ in
	\eqref{equ:dif 2}. Using integration by parts for all terms,
	we have
	\begin{align}\label{equ:cont dep H3}
		\nonumber
		&\norm{\partial_t \nabla \bff{w}}{\bb{L}^2}^2
		+
		\beta_1 \ddt \norm{\Delta \bff{w}}{\bb{L}^2}^2
		+
		\beta_2 \ddt \norm{\nabla \Delta \bff{w}}{\bb{L}^2}^2
		\\
		\nonumber
		&=
		\beta_3 \ddt \norm{\nabla \bff{w}}{\bb{L}^2}^2 
		-
		\beta_3 \inpro{\nabla(|\bff{u}|^2 \bff{u}- |\bff{v}|^2 \bff{v})}{\nabla \partial_t \bff{w}}_{\bb{L}^2}
		\\
		\nonumber
		&\quad
		+
		\beta_4 \inpro{\nabla \bff{u} \times \Delta \bff{u} - \nabla \bff{v} \times \Delta \bff{v}}{\nabla \partial_t \bff{w}}_{\bb{L}^2}
		+
		\beta_4 \inpro{\bff{u} \times \nabla \Delta \bff{u} - \bff{v} \times \nabla \Delta \bff{v}}{\nabla \partial_t \bff{w}}_{\bb{L}^2}
		\\
		&\quad
		-
		\beta_5 \inpro{\nabla \Delta(|\bff{u}|^2 \bff{u}- |\bff{v}|^2 \bff{v})}{\nabla \partial_t \bff{w}}_{\bb{L}^2}.
	\end{align}
	Each inner product on the right-hand side of \eqref{equ:cont dep H3} can
	be estimated as follows. 
	For the first inner product, similarly to~\eqref{equ:cont dep H1 3rd} we
	have
	\begin{align}\label{equ:cont dep H3 1st}
		\nonumber
		\beta_3 \big| 
		\inpro{\nabla(|\bff{u}|^2 \bff{u}- |\bff{v}|^2 \bff{v})}{\nabla \partial_t \bff{w}}_{\bb{L}^2}
		\big|
		&\lesssim 
		\big( \norm{\bff{u}}{\bb{H}^1}^2 + \norm{\bff{v}}{\bb{H}^1}^2 \big)
		\big( \norm{\bff{u}}{\bb{H}^2}^2 + \norm{\bff{v}}{\bb{H}^2}^2 \big) 
		\norm{\bff{w}}{\bb{H}^1}^2
		\\
		&\quad +
		\epsilon \norm{\nabla \partial_t \bff{w}}{\bb{L}^2}^2
	\end{align}
	for any $\epsilon>0$.
	For the second inner product, applying H\"older's and Young's inequality
	yields
	\begin{align}\label{equ:cont dep H3 2nd}
		\nonumber
		&\beta_4 \big| 
		\inpro{\nabla \bff{u} \times \Delta \bff{u} - \nabla \bff{v} \times \Delta \bff{v}}{\nabla \partial_t \bff{w}}_{\bb{L}^2} \big|
		\\
		\nonumber
		&\lesssim
		\norm{\nabla \bff{u} \times \Delta \bff{u} - \nabla \bff{v}
		\times \Delta \bff{v}}{\bb{L}^2}^2
		+
		\epsilon
		\norm{\nabla \partial_t \bff{w}}{\bb{L}^2}^2
		\\
		\nonumber
		&\lesssim
		\norm{\nabla \bff{u}}{\bb{L}^\infty}^2
		\norm{\Delta \bff{w}}{\bb{L}^2}^2
		+
		\norm{\nabla \bff{w}}{\bb{L}^\infty}^2
		\norm{\Delta \bff{v}}{\bb{L}^2}^2
		+
		\epsilon
		\norm{\nabla \partial_t \bff{w}}{\bb{L}^2}^2
		\\
		&\lesssim
		\big( \norm{\bff{u}}{\bb{H}^3}^2 + \norm{\bff{v}}{\bb{H}^2}^2 \big)
		\norm{\bff{w}}{\bb{H}^3}^2
		+
		\epsilon \norm{\nabla \partial_t \bff{w}}{\bb{L}^2}^2
	\end{align}
	for any $\epsilon>0$, 
	where in the last step we used the Sobolev
	embedding $\bb{H}^2\subset \bb{L}^\infty$. 
	For the third inner product we have by H\"older's inequality, Young's
	inequality, and~\eqref{equ:u cross Dk u}
	\begin{align}\label{equ:cont dep H3 3rd}
		\nonumber
		&\beta_4 \big| 
		\inpro{\bff{u} \times \nabla \Delta \bff{u} - \bff{v} \times \nabla \Delta \bff{v}}{\nabla \partial_t \bff{w}}_{\bb{L}^2}
		\big|
		\\
		\nonumber
		&\lesssim
		\norm{\bff{u}}{\bb{L}^\infty}^2
		\norm{\nabla \Delta \bff{w}}{\bb{L}^2}^2
		+
		\norm{\bff{w}}{\bb{H}^2}^2
		\norm{\nabla \Delta \bff{v}}{\bb{L}^2}^2
		+
		\epsilon
		\norm{\nabla \partial_t \bff{w}}{\bb{L}^2}^2
		\\
		&\lesssim
		\big( \norm{\bff{u}}{\bb{H}^2}^2 + \norm{\bff{v}}{\bb{H}^3}^2 \big)
		\norm{\bff{w}}{\bb{H}^3}^2
		+
		\epsilon \norm{\nabla \partial_t \bff{w}}{\bb{L}^2}^2
	\end{align}
	for any $\epsilon>0$. For the last inner product in \eqref{equ:cont dep
	H3}, by using H\"older's and Young's inequality, and~\eqref{equ:D k
	u2u-v2v}, we have
	\begin{align}\label{equ:cont dep H3 4th}
		\beta_5 \big| 
		\inpro{\nabla \Delta(|\bff{u}|^2 \bff{u}- |\bff{v}|^2 \bff{v})}{\nabla \partial_t \bff{w}}_{\bb{L}^2}
		\big| 
		&\leq
		C \big( \norm{\bff{u}}{\bb{H}^3}^4 + \norm{\bff{v}}{\bb{H}^3}^4 \big)
		\norm{\bff{w}}{\bb{H}^3}^2
		+
		\epsilon \norm{\nabla \partial_t \bff{w}}{\bb{L}^2}^2.
	\end{align}
	Define
	\[
	\bff{\gamma}(s) := 1+ \norm{\bff{u}(s)}{\bb{H}^3}^4 + \norm{\bff{v}(s)}{\bb{H}^3}^4.
	\]
	Inserting \eqref{equ:cont dep H3 1st}, \eqref{equ:cont dep H3 2nd},
	\eqref{equ:cont dep H3 3rd} and \eqref{equ:cont dep H3 4th} into
	\eqref{equ:cont dep H3}, integrating over $(0,t)$, and choosing
	$\epsilon>0$ sufficiently small, we obtain
	\begin{align*}
		\norm{\nabla \Delta \bff{w}(t)}{\bb{L}^2}^2
		&\lesssim
		\norm{\bff{w}(0)}{\bb{H}^3}^2 
		+ 
		\int_0^t \bff{\gamma}(s) \big(\norm{\bff{w}(s)}{\bb{L}^2}^2 \ds 
		+
		\norm{\nabla \bff{w}(s)}{\bb{L}^2}^2 \ds
		+
		\norm{\nabla \Delta \bff{w}(s)}{\bb{L}^2}^2 \big) \ds
		\\
		&\lesssim
		\norm{\bff{w}(0)}{\bb{H}^3}^2 
		+
		\int_{0}^{t} \gamma(s) 
		\norm{\nabla \Delta \bff{w}(s)}{\bb{L}^2}^2 \ds,
	\end{align*}
	where we used \eqref{eq5}, \eqref{equ:cont dep u-v L2} and~\eqref{equ:uv uv0 Hr} with $r=1$.
	Gronwall's inequality yields
	\[
	\norm{\nabla \Delta \bff{w}(t)}{\bb{L}^2}^2 \leq C \norm{\bff{w}(0)}{\bb{H}^3}^2
	\]
	and the required estimate \eqref{equ:uv uv0 Hr} for $r=3$ then follows.
}

\bigskip
\noindent
\underline{Extension from $[0,T^\ast]$ to $[0,T]$ for $d=3$}:
Recall that $T^\ast \le T$ for $d=3$. We will now show that in this case, we
also have $T^\ast = T$. First, it follows from~\eqref{equ:un bou} that, for $r=1,2,3$,
\begin{equation}\label{equ:del3d}
	\bff{u} \in
	L^\infty(0, T^\ast; \bb{H}^r) \cap L^2(0,T^\ast; \bb{H}^{r+2}).
\end{equation}

Assume that the following estimate holds (which will be shown in Proposition \ref{pro:nab Del 3d} later).
\begin{align}\label{equ:nab u d3}
\norm{\nabla \bff{u}(t)}{\bb{L}^2}^2 
&+ 
\norm{\bff{u}(t)}{\bb{L}^4}^4 
+  
\int_0^t \norm{\nabla \Delta \bff{u}(s)}{\bb{L}^2}^2 \ds   
\nonumber\\
&+ 
\int_0^t \norm{\nabla(|\bff{u}(s)|^2 \bff{u}(s))}{\bb{L}^2}^2 \ds 
+ 
\int_0^t \norm{\bff{u}(s)}{\bb{L}^6}^6 \ds 
\lesssim 
\norm{\bff{u}_0}{\bb{H}^1}^2,
\quad t\in[0,T^\ast].
\end{align}

Then we can repeat the arguments leading to the proof of \eqref{equ:un bou}
with $\bff{u}_n$ replaced by $\bff{u}$ to obtain similar estimates for
$\bff{u}$, with constant depending on $T$.  Proposition \ref{pro:nab Del 3d} and
\eqref{equ:un bou} imply that this weak solution $\bff{u}$ originally defined on
$[0,T^*]$ belongs to~$C([0,T^*]; \bb{H}^r) \cap L^2(0,T^*; \bb{H}^{r+2})$, and
that $\bff{u}(t,\bff{x})$ remains bounded in this norm as $t\to T^\ast$ from the
left. Therefore, the technique of continuation of solutions can be applied and
thus the solution $\bff{u}$ exists on the whole interval $[0,T]$ for any $T>0$.

It remains to prove~\eqref{equ:nab u d3}.

\begin{proposition}\label{pro:nab Del 3d}
	Let $T>0$ be arbitrary and $T^\ast$ be defined as in
	Proposition~\ref{pro:nab delta un}. Let $\bff{u}$ be the unique weak solution of \eqref{llbar}. 
%	For all $t\in [0,T^*]$,
%	\begin{align*}
%		\norm{\nabla \bff{u}(t)}{\bb{L}^2}^2 
%		&+ 
%		\norm{\bff{u}(t)}{\bb{L}^4}^4 
%		+  
%		\int_0^t \norm{\nabla \Delta \bff{u}(s)}{\bb{L}^2}^2 \ds   
%		\\
%		&+ 
%		\int_0^t \norm{\nabla(|\bff{u}(s)|^2 \bff{u}(s))}{\bb{L}^2}^2 \ds 
%		+ 
%		\int_0^t \norm{\bff{u}(s)}{\bb{L}^6}^6 \ds 
%		\lesssim 
%		\norm{\bff{u}_0}{\bb{H}^1}^2,
%	\end{align*}
	Then~\eqref{equ:nab u d3} holds with a constant depending on $T$.
\end{proposition}

\begin{proof}
We aim to choose $\bff{\phi} = \alpha \abs{\bff{u}(t)}^2 \bff{u}(t) $ in
\eqref{weakform}, for some positive constant $\alpha$. Hence, we first consider
the nonlinear terms in the resulting equation with that choice of $\bff{\phi}$.
For the term with coefficient $\beta_1$, we use~\eqref{equ:nab un2} to have
\begin{align}\label{equ:nab nab 3d}
	\inpro{\nabla\bff{u}(s)}{\nabla\big(|\bff{u}(s)|^2\bff{u}(s)\big)}_{\bb{L}^2}
=
2 \norm{\bff{u}(s) \cdot \nabla \bff{u}(s)}{\bb{L}^2}^2
+
\norm{ \abs{\bff{u}(s)} \abs{\nabla \bff{u}(s)}}{\bb{L}^2}^2.
\end{align}
For the term with coefficient $\beta_2$, we use integration by parts to have
\begin{align}\label{equ:Del Del 3d}
	\inpro{\Delta\bff{u}(s)}{\Delta\big(|\bff{u}(s)|^2\bff{u}(s)\big)}_{\bb{L}^2}
=
-
\inpro{\nabla\Delta\bff{u}(s)}{\nabla\big(|\bff{u}(s)|^2\bff{u}(s)\big)}_{\bb{L}^2}.
\end{align}
For the terms involving $\beta_3$ and $\beta_5$, it is straightforward to have
\begin{align*} 
	\inpro{\big( 1-|\bff{u}(s)|^2 \big) \bff{u}(s)}%
	{|\bff{u}(s)|^2 \bff{u}(s)}_{\bb{L}^2}
	&=
	\norm{\bff{u}(s)}{\bb{L}^4}^4
	-
	\norm{\bff{u}(s)}{\bb{L}^6}^6,
	\\
	\inpro{\nabla\big(|\bff{u}(s)|^2 \bff{u}(s) \big)}%
	{\nabla\big(|\bff{u}(s)|^2 \bff{u}(s) \big)}_{\bb{L}^2}
	&=
	\norm{\nabla\big(|\bff{u}(s)|^2 \bff{u}(s) \big)}{\bb{L}^2}^2.
\end{align*}
The term involving $\beta_4$ vanishes. Altogether, we deduce from choosing
$\bff{\phi} = 4\alpha \abs{\bff{u}(t)}^2 \bff{u}(t) $ in~\eqref{weakform} that
\begin{align}\label{equ:nab Del 3d-1}
		\alpha\norm{\bff{u}(t)}{\bb{L}^4}^4 
		&+ 
		8\alpha\beta_1 \int_0^t \norm{\bff{u}(s) \cdot \nabla \bff{u}(s)}{\bb{L}^2}^2 \ds 
		+ 
		4\alpha\beta_1 \int_0^t \norm{ \abs{\bff{u}(s)} \abs{\nabla \bff{u}(s)}}{\bb{L}^2}^2 \ds 
		\nonumber
		\\
		&\quad- 
		4\alpha\beta_2 
		\int_0^t \inpro{\nabla \Delta \bff{u}(s)}{\nabla 
		(\abs{\bff{u}(s)}^2 \bff{u}(s)) }_{\bb{L}^2} \ds 
		\nonumber
		\\
		&=
		4\alpha \inpro{\bff{u}_0}{|\bff{u}(t)|^2 \bff{u}(t)}_{\bb{L}^2}
		+ 
		4\alpha\beta_3 \int_0^t \norm{\bff{u}(s)}{\bb{L}^4}^4 \ds 
		-
		4\alpha\beta_3 \int_0^t \norm{\bff{u}(s)}{\bb{L}^6}^6 \ds 
		\nonumber\\
		&\quad
		-
		4\alpha\beta_5 \int_0^t \norm{\nabla (\abs{\bff{u}(s)}^2
		\bff{u}(s))}{\bb{L}^2}^2 \ds.
\end{align}

Next, we choose $\bff{\phi}= -2\Delta \bff{u}(t)$ in~\eqref{weakform} and use
integration by parts, noting~\eqref{equ:vec Gre} so that the term involving
$\beta_4$ vanishes. We then have, noting~\eqref{equ:nab nab 3d},
	\begin{align}\label{equ:nab Del 3d-2}
		\norm{\nabla \bff{u}(s)}{\bb{L}^2}^2  
		&+ 
		2\beta_1 \int_0^t \norm{\Delta \bff{u}(s)}{\bb{L}^2}^2 \ds 
		+ 
		2\beta_2 \int_0^t \norm{\nabla \Delta \bff{u}(s)}{\bb{L}^2}^2 \ds 
		\nonumber
		\\
		&= 
		2\inpro{\nabla\bff{u}_0}{\nabla \bff{u}(t)}_{\bb{L}^2}
		+ 
		2\beta_3 \int_0^t \norm{\nabla \bff{u}(s)}{\bb{L}^2}^2 \ds 
		\nonumber
		\\
		&\quad -
		4\beta_3 \int_0^t \norm{\bff{u}(s) \cdot \nabla \bff{u}(s)}{\bb{L}^2}^2 \ds 
		-
		2\beta_3 \int_0^t \norm{ \abs{\bff{u}(s)} \abs{\nabla \bff{u}(s)}}{\bb{L}^2}^2 \ds  
		\nonumber
		\\
		&\quad +
		2\beta_5 \int_0^t \inpro{\nabla \Delta \bff{u}(s)}{\nabla
		(\abs{\bff{u}(s)}^2 \bff{u}(s))} \ds.
	\end{align}
Adding \eqref{equ:nab Del 3d-1} and \eqref{equ:nab Del 3d-2} gives
\begin{align}\label{equ:nab Del add 3d}
	& \alpha \norm{\bff{u}(t)}{\bb{L}^4}^4 
	+ 
	\norm{\nabla \bff{u}(t)}{\bb{L}^2}^2 
	+ 
	2\beta_2 \int_0^t \norm{\nabla \Delta \bff{u}(s)}{\bb{L}^2}^2 \ds
	+
	4 \alpha \beta_3 \int_0^t \norm{\bff{u}(s)}{\bb{L}^6}^6 \ds 
	\nonumber
	\\
	&\quad -
	(4\alpha \beta_2 + 2\beta_5) \int_0^t \inpro{\nabla \Delta \bff{u}(s)}{\nabla(|\bff{u}(s)|^2 \bff{u}(s))}_{\bb{L}^2} \ds
	+
	4\alpha \beta_5 \int_0^t \norm{\nabla(|\bff{u}(s)|^2 \bff{u}(s))}{\bb{L}^2}^2 \ds
	\nonumber
	\\
	&= 
	\alpha \inpro{\bff{u}_0}{|\bff{u}(t)|^2 \bff{u}(t)}_{\bb{L}^2}
	+
	\inpro{\nabla\bff{u}_0}{\nabla \bff{u}(t)}_{\bb{L}^2}
	\nonumber
	\\
	&\quad -
	2\beta_1 \int_0^t \norm{\Delta \bff{u}(s)}{\bb{L}^2}^2 \ds 
	+ 
	4 \alpha \beta_3 \int_0^t \norm{\bff{u}(s)}{\bb{L}^4}^4 \ds 
	+ 
	2\beta_3 \int_0^t \norm{\nabla \bff{u}(s)}{\bb{L}^2}^2 \ds 
	\nonumber
	\\
	&\quad - 
	(8 \alpha \beta_1 + 4\beta_3) \int_0^t \norm{\bff{u}(s) \cdot \nabla \bff{u}(s)}{\bb{L}^2}^2 \ds  
	- 
	(4 \alpha \beta_1 + 2\beta_3)  \int_0^t \norm{ \abs{\bff{u}(s)} \abs{\nabla \bff{u}(s)}}{\bb{L}^2}^2 \ds.
\end{align}
Note that if $\alpha= \beta_5/2\beta_2$, then the third, fifth, and sixth terms
on the left-hand side add up to
	\begin{align*}
		2\beta_2 \norm{\nabla \Delta \bff{u}(s)}{\bb{L}^2}^2 
		&- 
		(4\alpha \beta_2 + 2\beta_5) \inpro{\nabla \Delta \bff{u}(s)}{\nabla (\abs{\bff{u}(s)}^2 \bff{u}(s))} 
		+
		4\alpha \beta_5 \norm{\nabla (\abs{\bff{u}(s)}^2 \bff{u}(s))}{\bb{L}^2}^2
		\\
		&= 
		\norm{\sqrt{2\beta_2} \nabla \Delta \bff{u}(s) 
		- 
		\sqrt{4\alpha \beta_5} \nabla (\abs{\bff{u}(s)}^2 \bff{u}(s))}{\bb{L}^2}^2 
		\geq 0.
	\end{align*}
Hence, with this value of $\alpha$ and the use of Young's inequality, \eqref{equ:nab Del add 3d} becomes
	\begin{align*}
		& \alpha \norm{\bff{u}(t)}{\bb{L}^4}^4 
		+ 
		\norm{\nabla \bff{u}(t)}{\bb{L}^2}^2 
		+ 
		4\alpha \beta_3 \int_0^t \norm{\bff{u}(s)}{\bb{L}^6}^6 \ds 
		+ 
		\int_0^t \norm{\sqrt{2\beta_2} \nabla \Delta \bff{u} 
		- 
		\sqrt{4\alpha\beta_5} \nabla (\abs{\bff{u}}^2 \bff{u})}{\bb{L}^2}^2 \ds
		\\
		&\quad \leq 
		C \norm{\bff{u}_0}{\bb{L}^4}^4
		+
		\epsilon \norm{\bff{u}(t)}{\bb{L}^4}^4
		+
		C \norm{\nabla \bff{u}_0}{\bb{L}^2}^2
		+
		\epsilon \norm{\nabla \bff{u}(t)}{\bb{L}^2}^2
		\\
		&\qquad
		+
		2 |\beta_1| \int_0^t \norm{\Delta \bff{u}(s)}{\bb{L}^2}^2 \ds 
		+ 
		4\alpha \beta_3 \int_0^t \norm{\bff{u}(s)}{\bb{L}^4}^4 \ds 
		+ 
		2\beta_3 \int_0^t \norm{\nabla \bff{u}(s)}{\bb{L}^2}^2 \ds 
		\\
		&\qquad 
		+ 
		(8\alpha |\beta_1| + 4\beta_3) \int_0^t \norm{\bff{u}(s) \cdot \nabla \bff{u}(s)}{\bb{L}^2}^2 \ds  
		+ 
		(4\alpha |\beta_1| + 2\beta_3)  \int_0^t \norm{ \abs{\bff{u}(s)}
		\abs{\nabla \bff{u}(s)}}{\bb{L}^2}^2 \ds
		\\
		&\le
		C \norm{\bff{u}_0}{\bb{H}^1}^2
		+
		\epsilon \norm{\bff{u}(t)}{\bb{L}^4}^4
		+
		\epsilon \norm{\nabla \bff{u}(t)}{\bb{L}^2}^2,
\end{align*}
where in the last step we used the Sobolev embedding $\bb{H}^1 \subset \bb{L}^4$
(for $\bff{u}_0$) and Proposition~\ref{pro:Delta un} to bound all the integrals
on the right-hand side. Choosing $\epsilon >0$ sufficiently small, we obtain the
required estimate.
\end{proof}

%\section{H\"{o}lder Continuity in Time of the Solution}\label{sec:Hol}
\section{Proof of Theorem~\ref{the:Hol}}\label{sec:Hol}

\begin{proof}
	For any Banach space $\bb{X}$, since $C^{0,\alpha_2}([0,T];\bb{X}) \subset
	C^{0,\alpha_1}([0,T];\bb{X})$ for $0 < \alpha_1 < \alpha_2$, it
	suffices to prove the theorem for $\alpha=1/2$ and $\beta=1/2-d/8$.

Let $T>0$ and $\tau, t\in [0,T]$ be such that $\tau < t$. Performing integration
by parts on \eqref{weakform} (and noting the regularity of the
solution~$\bff{u}$ given by Theorem \ref{the:weakexist}), we have for any
$\bff{\phi} \in \bb{H}^2$,
\begin{align*}
		\inpro{\bff{u}(t)-\bff{u}(\tau)}{\bff{\phi}}_{\bb{L}^2}
		&-
		\beta_1 \int_\tau^t \inpro{\Delta \bff{u}(s)}{\bff{\phi}}_{\bb{L}^2} \ds
		+
		\beta_2 \int_\tau^t \inpro{\Delta^2 \bff{u}(s)}{\bff{\phi}}_{\bb{L}^2} \ds
		\\
		&=
		\beta_3 \int_\tau^t \inpro{(1-|\bff{u}(s)|^2) \bff{u}(s)}{\bff{\phi}}_{\bb{L}^2} \ds
		-
		\beta_4 \int_\tau^t \inpro{\bff{u}(s) \times \Delta \bff{u}(s)}{\bff{\phi}}_{\bb{L}^2} \ds
		\\
		&\quad +
		\beta_5 \int_\tau^t \inpro{\Delta(|\bff{u}(s)|^2 \bff{u}(s))}{\bff{\phi}}_{\bb{L}^2} \ds.
\end{align*}
Therefore, by H\"{o}lder's inequality,
\begin{align*}
	\big| \inpro{\bff{u}(t)-\bff{u}(\tau)}{\bff{\phi}}_{\bb{L}^2} \big| 
	&\leq
	|\beta_1| \norm{\bff{\phi}}{\bb{L}^2} \int_\tau^t \norm{\Delta \bff{u}(s)}{\bb{L}^2} \ds
	+
	\beta_2 \norm{\bff{\phi}}{\bb{L}^2} \int_\tau^t \norm{\Delta^2 \bff{u}(s)}{\bb{L}^2} \ds
	\\
	&\quad
	+
	\beta_3 \norm{\bff{\phi}}{\bb{L}^2} \int_\tau^t \norm{\bff{u}(s)}{\bb{L}^2}  \ds
	+
	\beta_3 \norm{\bff{\phi}}{\bb{L}^2} \int_\tau^t \norm{\bff{u}(s)}{\bb{L}^6}^3  \ds
	\\
	&\quad
	+
	\beta_4 \norm{\bff{\phi}}{\bb{L}^2} \int_\tau^t \norm{\bff{u}(s)\times \Delta \bff{u}(s)}{\bb{L}^2}  \ds
	\\
	&\quad
	+
	\beta_5 \norm{\bff{\phi}}{\bb{L}^2} \int_\tau^t \norm{\Delta (|\bff{u}(s)|^2 \bff{u}(s))}{\bb{L}^2}  \ds.
\end{align*}
Taking $\bff{\phi}= \bff{u}(t)-\bff{u}(\tau)$, we obtain
\begin{align}\label{equ:ut utau}
	\norm{\bff{u}(t)- \bff{u}(\tau)}{\bb{L}^2} 
	&\lesssim 
	\int_\tau^t  \norm{ \Delta \bff{u}(s)}{\bb{L}^2} \ds
	+ \int_\tau^t \norm{\Delta^2 \bff{u}(s)}{\bb{L}^2} \ds
	+ \int_\tau^t \norm{\bff{u}(s) }{\bb{L}^2} \ds
	\nonumber\\
	&\quad 
	+ \int_\tau^t \norm{\bff{u}(s)}{\bb{L}^6}^3  \ds
	+ \int_\tau^t \norm{ \bff{u}(s) \times \Delta \bff{u}(s)}{\bb{L}^2} \ds
	\nonumber\\
	&\quad
	+ \int_\tau^t \norm{ \Delta (|\bff{u}(s)|^2 \bff{u}(s)) }{\bb{L}^2} \ds.
\end{align}
We will now estimate each term on the right-hand side of~\eqref{equ:ut utau}.
For the linear terms, by H\"{o}lder's inequality and Corollary~\ref{cor:un bou},
\begin{align*}
	\int_\tau^t  \norm{\Delta \bff{u}(s)}{\bb{L}^2} \ds
	&\leq
	|t-\tau|^{\frac{1}{2}} \norm{\Delta \bff{u}}{L^2(0,T;\bb{L}^2)}
	\lesssim
	|t-\tau|^{\frac{1}{2}}
	,
	\\
	\int_\tau^t  \norm{ \Delta^2 \bff{u}(s)}{\bb{L}^2} \ds
	&\leq
	|t-\tau|^{\frac{1}{2}} \norm{\Delta^2 \bff{u}}{L^2(0,T;\bb{L}^2)}
	\lesssim
	|t-\tau|^{\frac{1}{2}},
	\\
	\int_\tau^t \norm{\bff{u}(s) }{\bb{L}^2}  \ds 
	&\leq
	|t-\tau|^{\frac{1}{2}} \norm{\bff{u}}{L^2(0,T;\bb{L}^2)}
	\lesssim
	|t-\tau|^{\frac{1}{2}}.
\end{align*}
For the nonlinear terms on the right-hand side of~\eqref{equ:ut utau}, by
H\"{o}lder's inequality, Corollary~\ref{cor:un bou} and the Sobolev embedding,
\begin{align*}
	\int_\tau^t \norm{\bff{u}(s)}{\bb{L}^6}^3  \ds
	&\leq
	|t-\tau|^{\frac{1}{2}} \norm{\bff{u}}{L^6(0,T;\bb{L}^6)}^3
	\lesssim
	|t-\tau|^{\frac{1}{2}} \norm{\bff{u}}{L^\infty(0,T;\bb{H}^1)}^3
	\\
	&\lesssim
	|t-\tau|^{\frac{1}{2}},
	\\
	\int_\tau^t \norm{ \bff{u}(s) \times \Delta \bff{u}(s)}{\bb{L}^2} \ds
	&\leq
	\int_\tau^t \norm{\bff{u}(s)}{\bb{L}^\infty} \norm{\Delta \bff{u}(s)}{\bb{L}^2} \ds
	\\
	&\leq
	\norm{\Delta \bff{u}}{L^\infty(0,T;\bb{L}^2)}
	|t-\tau|^{\frac{1}{2}} \norm{\bff{u}}{L^2(0, T; \bb{L}^\infty)}
	\\
	&\leq
	|t-\tau|^{\frac{1}{2}} \norm{\bff{u}}{L^2(0, T; \bb{H}^2)}
	\norm{\Delta \bff{u}}{L^\infty(0,T;\bb{L}^2)}
	\\
	&\lesssim
	|t-\tau|^{\frac{1}{2}},
	\\
	\int_\tau^t \norm{ \Delta (|\bff{u}(s)|^2 \bff{u}(s)) }{\bb{L}^2}  \ds 
	&\leq 
	\int_\tau^t \norm{ |\bff{u}(s)|^2 \bff{u}(s) }{\bb{H}^2}  \ds 
	\le
	\int_\tau^t \norm{\bff{u}(s)}{\bb{H}^2}^3  \ds
	\\
	&\leq
	|t-\tau|^{\frac{1}{2}} \norm{\bff{u}}{\bb{L}^6(0,T;\bb{H}^2)}^3
	\leq
	|t-\tau|^{\frac{1}{2}} \norm{\bff{u}}{\bb{L}^\infty(0,T;\bb{H}^2)}^3
	\\
	&\lesssim
	|t-\tau|^{\frac{1}{2}},
\end{align*}
where for the last nonlinear term, we also used \eqref{eq7}.
Altogether, we derive from~\eqref{equ:ut utau} that $\bff{u} \in
C^{0,\alpha}(0,T; \bb{L}^2)$ for $\alpha \in (0,1/2]$.

Finally, by the Gagliardo--Nirenberg inequality (Theorem~\ref{the:Gal Nir} with
$\bff{v}= \bff{u}(t)-\bff{u}(\tau)$, $r=0$, $q=\infty$, $s_1=0$, $s_2=2$),
\begin{align*}
	\norm{ \bff{u}(t)- \bff{u}(\tau) }{\bb{L}^\infty} 
	&\lesssim 
	\norm{\bff{u}(t)-\bff{u}(\tau)}{\bb{L}^2}^{1- \frac{d}{4}}
	\norm{\bff{u}(t)-\bff{u}(\tau)}{\bb{H}^{2}}^{\frac{d}{4}}
	\\
	&\lesssim
	\norm{\bff{u}(t)-\bff{u}(\tau)}{\bb{L}^2}^{1- \frac{d}{4}}
	\norm{\bff{u}}{C([0,T];\bb{H}^2)}^{ \frac{d}{4} }
	\\
	&\lesssim
	|t-\tau|^{\frac{1}{2}- \frac{d}{8}},
\end{align*}
where in the penultimate step we used Theorem~\ref{the:weakexist} and in the
last step we used the previous part of this theorem.
\end{proof}

\section{Appendix}

We collect in this section a few results which were extensively used in this
paper.

\begin{theorem}[Gronwall--Bihari's inequality \cite{Bih56, BoyFab13}]\label{the:bihari}
Let $f$ be a non-decreasing continuous function which is non-negative on $[0,
\infty)$ such that $\int_1^\infty 1/f(x)\,\dx < \infty$. Let $F$ be the
anti-derivative of $-1/f$ which vanishes at $\infty$. Let
$y:[0,\infty)\to[0,\infty)$ be a continuous function and let $g$ be a locally
integrable non-negative function on $[0,\infty)$. Suppose that there exists $y_0>0$
such that for all $t \geq 0$,
	\begin{align*}
		y(t) \leq y_0 + \int_0^t g(s) \,\ds + \int_0^t f(y(s)) \,\ds.
	\end{align*}
Let $T^\ast$ be the unique solution of the equation
	\begin{align*}
		T^\ast = F \left(y_0 + \int_0^{T^\ast} g(s) \,\ds \right).
	\end{align*}
Then for any $T' \in (0, T^\ast)$, we have
	\begin{align} \label{equ:bihari}
		\sup_{0\leq t\leq T'} y(t) \leq F^{-1} 
		\left(F \Big(y_0 + \int_0^{T'} g(s) \,\ds \Big) - T' \right).
	\end{align}
Note that the expression on the right-hand side of \eqref{equ:bihari} tends to
$\infty$ as $T'\to T^\ast$.
\end{theorem}

The following theorem is a special case of a more general result in
\cite{BreMir19}.

\begin{theorem}[Gagliardo--Nirenberg inequalities]\label{the:Gal Nir}
Let $\Omega$ be a bounded domain of $\bb{R}^d$ with Lipschitz boundary, and
let~$\bff{v}: \Omega \to \bb{R}^3$. Then
\begin{align}\label{gagliardo}
	\|\bff{v}\|_{\bb{W}^{r,q}} 
	\leq 
	C \|\bff{v}\|_{\bb{H}^{s_1}}^\theta \|\bff{v}\|_{\bb{H}^{s_2}}^{1-\theta}
\end{align}
for all $\bff{v} \in \bb{H}^{s_2}(\Omega)$, where $s_1,s_2,r$ are non-negative
real numbers satisfying
\begin{align*}
	0\leq s_1<s_2,\quad \theta \in (0,1),\quad 0\leq r <\theta s_1 +(1-\theta)s_2, 
\end{align*}
and $q\in(2,\infty]$ satisfies
\[
	\frac{1}{q}=\frac{1}{2}+\frac{(s_2-s_1)\theta}{d}- \frac{s_2-r}{d}.
\]
Moreover, when $2< q<\infty$, we have
\begin{align*}
	\theta =\frac{2q(s_2-r)-d(q-2)}{2q(s_2-s_1)}.
%\quad \text{and} \quad \frac{2q(s_1-r)}{q-2} < d< \frac{2q(s_2-r)}{q-2}.
\end{align*}
\end{theorem}

\begin{theorem}[Aubin--Lions--Simon lemma \cite{Sim87}]\label{aubin}
Let $X_0\hookrightarrow X\hookrightarrow X_1$ be three Banach spaces such that the inclusion
$X_0 \hookrightarrow X$ is compact and the inclusion $X \hookrightarrow X_1$ is
continuous. For $1\leq p,q\leq \infty$, let
\begin{align*}
	\bb{W}_{p,q} := \{\bff{v} \in L^p(0,T;X_0) : \bff{v}_t \in L^q(0,T;X_1)\}.
\end{align*}
\begin{enumerate}
	\item If $p<\infty$, then $\bb{W}_{p,q}$ is compactly embedded into $L^p(0,T;X)$.
	\item If $p=\infty$ and $q=1$, then $\bb{W}_{p,q}$ is compactly embedded
		into $C([0,T];X)$.
\end{enumerate}
\end{theorem}

{
\begin{theorem}[Theorem II.5.14 in \cite{BoyFab13}]\label{the:E22 to C}
	Let $V$ and $W$ be Hilbert spaces. Then the space
	\[
		\{\bff{v} \in L^2(0,T; V): \partial_t\bff{v} \in L^2(0,T; W)\}
	\]
	is continuously embedded into $C([0,T]; [V,W]_{1/2})$. Here,
	$[V,W]_{1/2}$ is the interpolation space of $V$ and $W$ with order
	$\frac{1}{2}$.
\end{theorem}
}

\section*{Acknowledgements}
The authors gratefully thank the anonymous reviewer for the helpful remarks and suggestions, which greatly improved the clarity of the paper.

The first author is supported by the Australian Government Research Training Program Scholarship awarded at the University of New South Wales, Sydney. The second author is partially supported by the Australian Research Council under grant number DP190101197 and DP200101866.

%%%%%%%%%%%%%%% Bibliography %%%%%%%%%%%%%%%%%%

\newcommand{\noopsort}[1]{}\def\cprime{$'$}
\def\soft#1{\leavevmode\setbox0=\hbox{h}\dimen7=\ht0\advance \dimen7
	by-1ex\relax\if t#1\relax\rlap{\raise.6\dimen7
		\hbox{\kern.3ex\char'47}}#1\relax\else\if T#1\relax
	\rlap{\raise.5\dimen7\hbox{\kern1.3ex\char'47}}#1\relax \else\if
	d#1\relax\rlap{\raise.5\dimen7\hbox{\kern.9ex \char'47}}#1\relax\else\if
	D#1\relax\rlap{\raise.5\dimen7 \hbox{\kern1.4ex\char'47}}#1\relax\else\if
	l#1\relax \rlap{\raise.5\dimen7\hbox{\kern.4ex\char'47}}#1\relax \else\if
	L#1\relax\rlap{\raise.5\dimen7\hbox{\kern.7ex
			\char'47}}#1\relax\else\message{accent \string\soft \space #1 not
		defined!}#1\relax\fi\fi\fi\fi\fi\fi}

\end{document}